\documentclass[a4paper,  dvipsnames, 11pt]{amsart}
\usepackage{preamble}

\tikzstyle{cell}=[fill=white, draw=black, shape=rectangle, minimum width=1cm, minimum height=0.5cm]
\tikzstyle{big cell}=[fill=white, draw=black, shape=rectangle, minimum width=2cm, minimum height=0.5cm]
\tikzstyle{huge cell}=[fill=white, draw=black, shape=rectangle, minimum width=3cm, minimum height=0.5cm]

\tikzstyle{white}=[-, draw=white, line width=5pt]

\begin{document}
\begin{abstract}
	We relativise double categories of relations to stable orthogonal factorisation systems. 
	Furthermore, we present the characterisation of the relative double categories of relations in two ways.
	The first utilises a generalised comprehension scheme, and the second focuses on a specific class of vertical arrows defined solely double-categorically. 
	We organise diverse classes of double categories of relations and correlate them with significant classes of factorisation systems. 
	Our framework embraces double categories of spans and double categories of relations on regular categories, 
	which we meticulously compare to existing work on the characterisations of bicategories and double categories of spans and relations.
\end{abstract}
\maketitle
\tableofcontents
\section{Introduction}
Sets, functions, and relations are fundamental concepts in mathematics.
The category $\Set$ of sets and functions has, on the one hand, been one of the most motivating prototypes of categories,
and on the other hand, the bicategory $\biRel{}$ of sets and binary relations 
has also been studied as a primary entity in category theory.
Afterwards, relations were extended to any category with finite products
by formulating relations between objects as subobjects of their product.
Here, the category must be regular for the composition of relations to be associative and unital, as observed by Lawvere \cite{Law72}.
Subsequently, Succi \cite{Suc75} gave a characterisation of
the bicategory of relations on regular categories using involutions, 
and Carboni, Kasangian, and Street \cite{CKS84} addressed the same problem without such additional structures. 

Relations were further generalised to any category $\one{C}$ with finite products and a 
proper stable factorisation system $(\zero{E}, \zero{M})$ on it by Klein \cite{Kl70}.
He defined the notion of $\zero{M}$-relations in $\one{C}$ by replacing subobjects with $\zero{M}$-subobjects. 
Here, $\zero{M}$-subobjects of an object $A$ are arrows in $\zero{M}$ into $A$\footnote{
$\zero{M}$-subobjects in a category are usually referred to as an equivalence class of morphisms in $\zero{M}$ into an object,
but here we use the term `subobject' to mean the morphism itself.
}.
Later on, Kelly \cite{Kel91} investigated the bicategories of $\zero{M}$-relations $\biRel{\zero{E},\zero{M}}(\one{C})$ in this context.
Then, Pavlovi\'c \cite{Pav95} extended the scope of relativised relations beyond the assumption of properness on the factorisation system
so that they encompass many different situations consistently.
For instance, the non-proper \acl{SOFS} $(\Iso,\Mor)$ whose left class is the class of isomorphisms
and whose right class is the class of all morphisms
yields the bicategory of spans $\biSpan(\one{C})$ on a finitely complete category $\one{C}$.
The bicategories of spans were initially introduced by B\'enabou \cite{Ben67} and 
have been inspected by several authors (see \cite{CKS84,CW87,LWW10} and the references therein).

On the other hand, there are many cases in which one would like to deal with relations and the original category at the same time.
Double categories, devised by Ehresmann \cite{Ehr63}, fulfil this desire by having vertical arrows as arrows in the original category 
and horizontal arrows as relations between objects.
Besides, one can express the interplay between vertical and horizontal arrows 
by dint of various structures on double categories,
such as companions and conjoints \cite{GP04,Shu08}, double limits \cite{GP99}, and cartesian structures \cite{Ale18}.
This remarkable abundance of structures has opened the possibility of studying relations in a double-categorical framework.
Indeed, Lambert \cite{Lam22} defined and characterised the \acl{DCR} $\Rel{}(\one{C})$ on a regular category $\one{C}$, 
in which vertical arrows are morphisms in $\one{C}$, and horizontal arrows are relations in $\one{C}$.
His method to construct an equivalence between an axiomatised \acl{DCR} and one 
in the form of $\Rel{}(\one{C})$ originated from the work of Niefield \cite{Nie12}.
She determined the condition for a double category to admit an oplax/lax adjoint to the double category of spans.
Later, Aleiferi \cite{Ale18} extended the result and
characterised\footnote{
	The characterisation in \cite{Ale18} needs a slight modification as we will point out in \cref{rem:SpanchaAleiferi}.
}
double categories of spans.

The principal objective of this paper is to generalise these results \cite{Ale18,Lam22}
to the cases relative to \acp{SOFS}.
We provide the characterisation of the \acp{DCR} $\Rel{\zero{E}, \zero{M}}(\one{C})$ 
formed by relations defined using a general \acl{SOFS} $(\zero{E}, \zero{M})$ on a finitely complete category $\one{C}$,
as shown in the following theorem.
\begin{theorem*}[\cref{thm:MainThm}]
	The following are equivalent for a double category $\dbl{D}$.
	\begin{enumerate}
		\item
			$\dbl{D}$ is equivalent to $\Rel{\zero{E}, \zero{M}}(\one{C})$
			for some category $\one{C}$ with finite limits and a \acl{SOFS}
			$(\zero{E}, \zero{M})$ on $\one{C}$.
		\item
			$\dbl{D}$ is a cartesian equipment with
			Beck-Chevalley pullbacks and
			$\zero{M}$-comprehension scheme for some stable system
			$\zero{M}$ on $\oneV{\dbl{D}}$.
		\item
			$\dbl{D}$ is a cartesian equipment with
			Beck-Chevalley pullbacks and
			a left-sided $\zero{M}$-comprehension scheme for some stable system
			$\zero{M}$ on $\oneV{\dbl{D}}$.
		\item
			$\Fib(\dbl{D})$ is closed under composition, and
			$\dbl{D}$ is a cartesian equipment with
			Beck-Chevalley pullbacks and
			strong tabulators.
		\item
			$\Fib(\dbl{D})$ is closed under composition, and
			$\dbl{D}$ is a cartesian equipment with
			Beck-Chevalley pullbacks and
			left-sided strong tabulators.
	\end{enumerate}
\end{theorem*}

As we only consider the case of finitely complete categories,
the double category $\dbl{D}$ is assumed to be a cartesian equipment,
meaning it is a cartesian double category with companions and conjoints.
Our goal is to characterise $i)$.
Condition $ii)$ that follows involves the $\zero{M}$-comprehension scheme for a class of vertical arrows $\zero{M}$,
a generalisation of the subobject comprehension scheme in \cite{Lam22}.
Through the $\zero{M}$-comprehension scheme, it is possible to relate the information of horizontal arrows to the vertical category.
Condition $iii)$ is a left-sided version of $ii)$ and is the key to proving the theorem in our general setting.
Here, a left-sided relation refers to one from some object $A$ to the terminal object $1$.
In the proof of the equivalence between $ii)$ and $iii)$,
relations between $A$ and $B$ are transformed into relations between $A\times B$ and $1$.
This manipulation is only possible when the cartesian equipment is `discrete.'
Discreteness is a condition that classifies cartesian bicategories of relations in \cite{CW87},
but we will see that this follows in the existence of so-called Beck-Chevalley pullbacks.
The remaining conditions $iv)$ and $v)$ are formulated with the properties and structures 
intrinsic to double categories.
The class $\Fib(\dbl{D})$ of vertical arrows called fibrations in these conditions, and another class $\Fin(\dbl{D})$ consisting of 
vertical arrows, called final arrows, are defined 
as candidates for the right and left classes of the factorisation system.
Final arrows will turn out to be a double-categorical analogue of the class of extremal epimorphisms in a category.
Besides their crucial role in the characterisation of \aclp{DCR},
it is worth attention that fibrations and final arrows apprehend the structure of double categories in a broader context.
For instance, these two classes in the double category of profunctors $\Prof$ 
become the class of discrete fibrations and final functors,
constituting a factorisation system on $\Catbi$.

Our second aim is to delineate the correspondence between \aclp{DCR} and \aclp{SOFS}.
As presented in \cref{tab:Correspondence}, the desirable properties on the \aclp{DCR} will be layered by the significant classes of factorisation systems.
A property called unit-pureness was an indispensable assumption in the previous literature \cite{Ale18,Lam22},
but the correspondence explicitly explains its effect.
Furthermore, there is a clear correspondence between the properness of factorisation systems and the posetality of \aclp{DCR}.

\begin{table}[h]
	\centering
	\begin{tabular}{|c|c|}
		\hline
		\Acp{SOFS} on finitely complete categories
		& Double categories of relations (\acp{DCR})
		\\
		\hline
		\hline
		\begin{tikzpicture}[scale=0.70,baseline=2ex]
			\node[myrectanglefs] (A)
				at (1, 0)
				{\ac{SOFS}};
			\node[myrectanglefs] (C1)
				at (0, -2)
				{left-proper \ac{SOFS}};
			\node[myrectanglefs] (C2)
				at (5, -2)
				{right-proper \ac{SOFS}};
			\node[myrectanglefs] (D1)
				at (-1, -4)
				{anti-right-proper \ac{SOFS}};
			\node[myrectanglefs, text width=20mm] (D2)
				at (4, -4)
				{proper \ac{SOFS}};
			\node[myrectanglefs] (E)
				at (3, -6)
				{regular \ac{SOFS}};
			\node[myrectanglefs] (F)
				at (-2, -6)
				{$(\Iso, \Mor)$};
			\draw[-] (A) -- (C1);
			\draw[-] (A) -- (C2);
			\draw[-] (C1) -- (D1);
			\draw[-] (C1) -- (D2);
			\draw[-] (C2) -- (D2);
			\draw[-] (D1) -- (E);
			\draw[-] (D2) -- (E);
			\draw[-] (D1) -- (F);
			\addvmargin{1mm}
		\end{tikzpicture}
		&
		\begin{tikzpicture}[scale=0.8,baseline=0]
			\node[myrectangledbl] (A)
				at (1, 0)
				{\Acf{DCR}\\
				\ref{thm:MainThm}};
			\node[myrectangledbl] (C1)
				at (0, -2)
				{unit-pure \ac{DCR}\\
				\ref{thm:unitpuredbl}};
			\node[myrectangledbl] (C2)
				at (5, -2)
				{locally preordered \ac{DCR}\\
				\ref{thm:locorddbl}};
			\node[myrectangledbl] (D1)
				at (-1, -4)
				{unit-pure Cauchy \ac{DCR}\\
				\ref{thm:Cauchydouble}};
			\node[myrectangledbl] (D2)
				at (4, -4)
				{locally posetal \ac{DCR}\\
				\ref{cor:locpos}};
			\node[myrectangledbl] (E)
				at (3, -6)
				{ \ac{DCR} on regular
				 categories\\
				\ref{thm:Regcat}, \ref{thm:RegcatCW}
				\cite{Lam22, CW87}
				};
			\node[myrectangledbl] (F)
				at (-2, -6)
				{Double category of spans\\
				\ref{cor:Spanchacompa}
				\cite{LWW10, Ale18}
				};
			\draw[-] (A) -- (C1);
			\draw[-] (A) -- (C2);
			\draw[-] (C1) -- (D1);
			\draw[-] (C1) -- (D2);
			\draw[-] (C2) -- (D2);
			\draw[-] (D1) -- (E);
			\draw[-] (D2) -- (E);
			\draw[-] (D1) -- (F);
			\addvmargin{1mm}
		\end{tikzpicture}
		\\
		\hline
	\end{tabular}
	\caption{Correspondence between \acfp{SOFS} and \acfp{DCR}}
	\label{tab:Correspondence}
	\end{table}

A noteworthy condition for unit-pure double categories is the Cauchy condition located in the third row of the table.
It states that every adjoint in the horizontal bicategory comes from the vertical category as a companion/conjoint adjunction.
The significance of the Cauchy condition lies in its connotation of the unique choice principle that every functional relation gives rise to a function,
as pointed out by Rosolini \cite{Ros99}.
Classical accounts of the bicategories of relations or spans have defined `functions' as `functional relations',
which obscured the Cauchy condition as a hidden assumption. 
However, by separating vertical arrows from horizontal ones, our framework makes the Cauchy condition explicit
and establishes a correspondence between the Cauchy condition on \aclp{DCR} and a condition on \aclp{SOFS} called anti-right-properness.
In addition, we implement the transformation process
from the unit-pure \acl{DCR} into a Cauchy one. 
It is the extension of the work by Kelly \cite{Kel91} on
converting a category with a proper factorisation system into a regular category.

We also contrast our results with prior work.
Carboni and Walters defined a `bicategory of relations' as a locally posetal cartesian bicategory with discrete objects in \cite{CW87}.
In this setting, they showed that
a `bicategory of relations' is equivalent to the horizontal part of $\Rel{\Regepi, \Mono}(\one{C})$ for some regular category $\one{C}$
if and only if it is `functionally complete.'
On the other hand, as already mentioned, much work has been done on the characterisation of spans.
Lack, Walters, and Wood give a way to characterise the bicategory of spans as a cartesian bicategory in terms of `comonad'
and its `co-Eilenberg-Moore objects' \cite{LWW10},
while the work by Aleiferi \cite{Ale18}, following Niefield \cite{Nie12},
captures spans as a double category using `copointed arrows.'
We explain how these concepts are interpreted in our settings and how their theorems are induced by narrowing down our general results.

\addtocontents{toc}{\protect\setcounter{tocdepth}{1}}
\subsection*{The outline of the paper}
\Cref{sec:Backgrounds} is devoted to the preliminaries of the theory of double categories and orthogonal factorisation systems.

\Cref{sec:axiomatising} is the pivotal part of this paper.
\Cref{sec:BCpb} illustrates how the Beck-Chevalley condition is used to bias relations to one side.
\Cref{sec:comprehensionscheme} deals with comprehension schemes for \aclp{SOFS}
and its left-sided version.
After the two classes of vertical arrows in double categories, called fibrations and final arrows, are formulated,
the main theorem is proved in \cref{subsec:Characterisation}.

\Cref{sec:severalclasses} comprises the discussions on different topics.
\Cref{sec:properness} shows how the left and right properness of \aclp{SOFS} are reflected in double categories.
\Cref{subsec:cauchyreg} discusses the Cauchy condition, especially the process of Cauchisation.
\Cref{sec:Contrast} compares our results with the previous literature \cite{CW87,LWW10,Ale18}.

\Cref{sec:futurework} is devoted to the future work.

\subsection*{Notation}
As already used in the introduction, we adopt the following notation.
We write categories in the bold font as $\one{C}$, $\one{D}$, and $\Set$.
For 2-categories and bicategories, we use 
the calligraphic letters such as $\bi{B}$, $\bi{C}$, and $\CATbi$.
For double categories, we write the initial letter in the blackboard bold font as $\dbl{D}$, $\dbl{E}$, and $\Prof$.
A letter in the Roman font denotes a class of arrows in a category, or that of vertical arrows in a double category,
such as $\zero{M}$, $\Mono$, and $\Fib$.

We write composites of arrows in the diagrammatic order opposite to the classical convention.
For instance, $f\fatsemi g$ denotes the composite $\cdot\to["f"]\cdot\to["g"]\cdot$ in the vertical category of some double category. 

\addtocontents{toc}{\protect\setcounter{tocdepth}{2}}

\section{Background}\label{sec:Backgrounds} 

This section reviews some notions and results that will be used in the remainder of the paper.
In particular, we review the basics of the theory of double categories and
factorisation systems.

\subsection{Double categories}\label{subsec:dblcat}
In this subsection, we take a look at several notions about double categories,
including equipments, tabulators, and cartesian double categories.
We use the notion of Grothendieck fibrations and simply call them fibrations.

By a \textit{(pseudo-)double category} $\dbl{D}$, we mean a pseudo-category in the 2-category $\CATbi$
of locally small categories.
In other words, a double category consists of
two (locally small) categories $\dbl{D}_0$, $\dbl{D}_1$ and 
functors 
\[
\begin{tikzcd}
	\dbl{D}_1\,_\tgt\!\times_{\src}\dbl{D}_1
	\ar[r, "\odot"] 
	&
	\dbl{D}_1
	\ar[r, shift left =2 , "\src"]
	\ar[r, shift right =2 , "\tgt"']
	&
	\dbl{D}_0
	\ar[l, "\Id"{description}]
\end{tikzcd}\text{.}
\]
These data come equipped with natural isomorphisms 
that stand for the associativity law and the unit laws.

Objects and arrows of $\dbl{D}_0$ are called
\emph{objects} and \emph{vertical arrows} of the double category $\dbl{D}$.
We use the notation $f\fatsemi g$ for the composition of $A\to["f"] B \to["g"] C$ in $\dbl{D}_0$.
An object $p$ of $\dbl{D}_1$ whose values of $\src$ and $\tgt$ 
are $A$ and $B$, respectively, is called a \emph{horizontal arrow}\footnote{
Note that in a significant portion of the literature,
vertical arrows are referred to as what we designate as horizontal arrows, and vice versa.
The difference cannot be dismissed since only one class of arrows requires strict associativity and unit laws.}
from $A$ to $B$, 
and written as $p\colon A\sto B$.
We use the notation $p\odot q$, or simply $p q$,
for the composite of $p\colon A\sto B$ and $q\colon B\sto C$ 
in $\dbl{D}_1$.
An arrow $\alpha\colon p\to q$ in $\dbl{D}_1$ is called a \emph{double cell}
(or merely a \emph{cell}) in the double category $\dbl{D}$.
This cell is drawn as below, where $\src(\alpha)=f$ and $\tgt(\alpha)=g$.
\begin{equation}
	\label{dgm:doubleCell}
	\begin{tikzcd}
		A 
		\ar[d, "f"']
		\sar[r,"p"]
		\ar[rd, phantom,"\alpha"]
		&
		B
		\ar[d, "g"]
		\\
		C
		\sar[r,"q"']
		&
		D
	\end{tikzcd}	
\end{equation}
Suppose moreover that $p$ and $q$ are obtained by composing paths of horizontal arrows
$\langle p_1,\ldots,p_n\rangle$ and
$\langle q_1,\ldots,q_m\rangle$,
where we have $p_i\colon A_{i-1}\sto A_{i}$, $q_i\colon C_{i-1}\sto C_i$, $A_0=A$, $A_n=B$, $C_0=C$, and $C_m=D$.
Then, we express $\alpha$ as on the left below.
This includes the case $n=0$ and/or $m=0$, as is on the right below. As usual, by the composition of a $0$-length path, we mean the horizontal identity.
\[
	\begin{tikzcd}
		A_0
		\ar[d, "f"']
		\sar[r,"p_1"]
		\doublecell[rrrd]{\alpha}
			&
			A_1
			\sar[r,"p_2"]
				&
				\cdots
				\sar[r,"p_n"]
					&
					A_n
					\ar[d, "g"]
		\\
		C_0
		\sar[r,"q_1"']
			&
			C_1
			\sar[r,"q_2"']
				&
				\cdots
				\sar[r,"q_m"']
					&
					C_m
	\end{tikzcd}
	\hspace{7ex}
	\begin{tikzcd}[column sep=small]
			&
			A_0
			\ar[ld, "f"']
			\ar[rd, "g"]
				&
		\\
		C_0
		\sar[r,"q_1"']
		\doublecell[rr, shift left=3ex]{\alpha}
			&
			C_1
			\sar[r,"q_2"']
				&
				C_2
	\end{tikzcd}
\]
The natural isomorphisms that stand for the associativity law and the unit laws describe the associativity and the unitality of the horizontal 
composition of horizontal arrows and that of cells.

Interchanging the roles of $\src$ and $\tgt$ in a double category $\dbl{D}$,
we obtain another double category.
We call it the \emph{horizontal opposite} of $\dbl{D}$
and write it as $\dbl{D}^\hop$.
Sending the data of $\dbl{D}$ by the 2-functor $(-)^\op\colon \CATbi^\co\to\CATbi$,
we get another double category.
We call it the \emph{vertical opposite} of $\dbl{D}$ and write it as $\dbl{D}^\vop$.

The following is a way to construct a bicategory from a double category, which is fundamental and has been discussed in several contexts.
See, for example, \cite{Gra20} for more discussion.
\begin{definition}
	For a double category $\dbl{D}$,
	the \emph{horizontal bicategory} $\biH{\dbl{D}}$ is a bicategory
	whose objects are objects in $\dbl{D}$,
	1-cells are horizontal arrows, 
	and 2-cells are cells in $\dbl{D}$ of the form on the left below, which are called \emph{horizontal cells} in $\dbl{D}$.

	Similarly, the \emph{vertical 2-category} $\biV{\dbl{D}}$ is a 2-category
	defined in the same way, 
	in which the composition of 1-cells is strict and coincides with that of $\dbl{D}_0$.
	We call a 2-cell $\beta\colon f\Rightarrow g$ in this 2-category a \emph{vertical cell} in $\dbl{D}$, which is a cell in $\dbl{D}$ of the form shown on the right.
	\[
		\begin{tikzcd}
			\cdot
			\sar[r, "p"]
			\ar[d, equal]
			\doublecell[rd]{\alpha}
				&
				\cdot
				\ar[d, equal]
			\\
			\cdot
			\sar[r, "q"']
				&
				\cdot
		\end{tikzcd}
		\hspace{3ex}
		,
		\hspace{3ex}
		\begin{tikzcd}
			\cdot
			\ar[d, "f"', bend right, shift right=0.2ex]
			\ar[d, "g", bend left, shift left=0.2ex]
			\doublecell[d]{\beta}
			\\
			\cdot
		\end{tikzcd}
	\]

	Note that $\biH{\dbl{D}}(X, Y)$ is exactly the fibre above $(X, Y)\in\dbl{D}_0\times\dbl{D}_0$
	of the functor $\langle\src,\tgt\rangle\colon \dbl{D}_1\to\dbl{D}_0\times\dbl{D}_0$.
\end{definition}
For consistency with this terminology, 
we write $\oneV{\dbl{D}}$ for the category $\dbl{D}_0$.
In other words, $\oneV{\dbl{D}}$ is the category underlying $\biV{\dbl{D}}$.

The above definition gives a way to form a bicategory of objects and horizontal arrows, while there also exists a slightly unfamiliar way
to construct a bicategory from a double category, consisting of vertical arrows and cells.
\begin{definition}
	For a double category $\dbl{D}$,
	the \emph{horizontal bicategory of cells}, $\bi{C}_h(\dbl{D})$, is a bicategory
	defined as follows.
	Objects are vertical arrows in $\dbl{D}$, and a 1-cell $g \to h$ is a triple $(p, \alpha, q)$ that forms a cell in $\dbl{D}$ of the following form.
	\[
		\begin{tikzcd}
			\cdot
			\sar[r, "p"]
			\ar[d, "g"']
			\doublecell[rd]{\alpha}
				&
				\cdot
				\ar[d, "h"]
			\\
			\cdot
			\sar[r, "q"']
				&
				\cdot
		\end{tikzcd}
	\]
	A 2-cell $(p, \alpha, q)\Rightarrow (p', \alpha', q')\colon g \to h$ is a pair $(\gamma,\delta)$ of horizontal cells in $\dbl{D}$ satisfying the following.
	\[
		\begin{tikzcd}
			\cdot
			\sar[r, "p"]
			\ar[d, "g"']
			\doublecell[rd]{\alpha}
				&
				\cdot
				\ar[d, "h"]
			\\
			\cdot
			\sar[r, "q"]
			\ar[d, equal]
			\doublecell[rd]{\delta}
				&
				\cdot
				\ar[d, equal]
			\\
			\cdot
			\sar[r, "q'"']
				&
				\cdot
		\end{tikzcd}
		\hspace{2ex}
		=
		\hspace{2ex}
		\begin{tikzcd}
			\cdot
			\sar[r, "p"]
			\ar[d, equal]
			\doublecell[rd]{\gamma}
				&
				\cdot
				\ar[d, equal]
			\\
			\cdot
			\sar[r, "p'"]
			\ar[d, "g"']
			\doublecell[rd]{\alpha'}
				&
				\cdot
				\ar[d, "h"]
			\\
			\cdot
			\sar[r, "q'"']
				&
				\cdot
		\end{tikzcd}
	\]
\end{definition}

The following are our first examples of \textit{double categories of relations} defined in \cref{def:RelEM}, the main objective of this paper.
\begin{example}
	The double category $\Rel{}(\Set)$ consists of the following data.
	\begin{itemize}
		\item%
			The vertical category $\oneV{\Rel{}(\Set)}$ is the category $\Set$ of sets and functions.
		\item%
			The horizontal arrows are binary relations between sets;
			i.e., a horizontal arrow $p\colon A\sto B$ is a subset of $A\times B$.
			A cell of the form \cref{dgm:doubleCell} exists if and only if 
			for any $a\in A$ and $b\in B$ such that $(a,b) \in p$,
			we have $(f(a),g(b))\in q$.
			There is at most one cell framed by a pair of vertical arrows and a pair of horizontal arrows. 
		\item%
			The composite $pq$ of relations $p\colon A\sto B$ and $q\colon B\sto C$ is defined as the following relation.
			For $a\in A$ and $c\in C$, we have $(a,c)\in pq$ if and only if there exists $b\in B$ such that $(a,b)\in p$ and $(b,c)\in q$.
			The unit horizontal arrow $\Id_A$ is defined by the diagonal $\{\,(a,a)\mid a\in A\,\}$.
	\end{itemize}
	This double category is introduced in \cite[\S 3.4]{GP99}
	and is generalised by replacing $\Set$ with any regular category in \cite{Lam22},
	which we characterise in \cref{thm:Regcat,thm:RegcatCW}. 
\end{example}

\begin{example}\label{ex:Span}
	Let $\one{C}$ be a category with finite limits.\footnote{%
		For the definition of $\Span(\one{C})$, not all finite limits are necessary; in fact, only pullbacks are sufficient.
		However, they are necessary for the double category to be an example of a \textit{double category of relations} in our sense.
	}
	The double category $\Span(\one{C})$ consists of the following data.
	\begin{itemize}
		\item%
			The vertical category $\oneV{\Span(\one{C})}=\Span(\one{C})_0$ is precisely the same as $\one{C}$.
			Therefore, objects and vertical arrows in $\Span(\one{C})$ are the same as objects and arrows in $\one{C}$.
		\item%
			$\langle\src,\tgt\rangle\colon\Span(\one{C})_1\to\one{C}\times\one{C}$ is defined by the following pullback.
			\[
				\begin{tikzcd}
					\Span(\one{C})_1
					\ar[d, "{\langle\src,\tgt\rangle}"']
					\ar[r]
					\pullback[rd]
						&
						\one{C}^{\rightarrow}
						\ar[d, "\cod^{\one{C}}"]
					\\
					\one{C}\times\one{C}
					\ar[r, "\times"']
						&
						\one{C}
				\end{tikzcd}
			\]
			Unpacking this definition, one finds that a horizontal arrow $R:A\sto B$ is a \textit{span} from objects $A$ to $B$; i.e.,
			a pair $(l_R,r_R)$ of arrows in $\one{C}$ with the same domain $\abs{R}$ such that the codomains of $l_R$ and $r_R$ are $A$ and $B$, respectively.
			A cell of the form on the left below is precisely an arrow $\alpha\colon\abs{R}\to\abs{S}$ between the apexes of $R$ and $S$ that makes the diagram on
			the right below commute.
			\[
				\begin{tikzcd}
					A
					\sar[r, "R"]
					\ar[d, "f"']
					\doublecell[rd]{\alpha}
						&
						B
						\ar[d, "g"]
					\\
					C
					\sar[r, "S"']
						&
						D
				\end{tikzcd}
				\hspace{5ex}
				\begin{tikzcd}
					A
					\ar[d, "f"']
						&
						\abs{R}
						\ar[l, "l_R"']
						\ar[r, "r_R"]
						\ar[d, "\alpha"description]
							&
							B
						\ar[d, "g"]
					\\
					C
						&
						\abs{S}
						\ar[l, "l_S"]
						\ar[r, "r_S"']
							&
							D
				\end{tikzcd}
			\]
		\item%
			For a composable pair of spans $R\colon A\sto B$ and $S\colon B\sto C$, the composite $R\odot S$ is defined as the following pullback.
			\[
				\begin{tikzcd}[row sep=small]
						&
							&
							\abs{R\odot S}
							\ar[ld, ""']
							\ar[lldd, bend right, "l_{R\odot S}"']
							\ar[rrdd, bend left, "r_{R\odot S}"]
							\ar[rd, ""]
							\pullbackA[dd]
								&
									&
					\\
						&
						\abs{R}
						\ar[ld, "l_R"]
						\ar[rd, "r_R"]
							&
								&
								\abs{S}
								\ar[ld, "l_S"']
								\ar[rd, "r_S"']
									&
					\\
					A
						&
							&
							B
								&
									&
									C
				\end{tikzcd} 
			\]
			The unit horizontal arrow $\Id_A$ is defined by the pair $(\id_A,\id_A)$,
			and the unit cell $\Id_f$ is defined by $f$ itself.
			This composition is associative and unital up to isomorphism by 
			the universal property of pullbacks.
			The composition of cells is defined by the universal property of pullbacks,
			which makes it associative and unital.
	\end{itemize}
		The double category $\Span(\one{C})$ is introduced in \cite[\S 3.2]{GP99} 
		in the case $\one{C}=\Set$,
		and general cases are discussed in \cite{Nie12}.
		We call these \emph{double categories of spans}.
\end{example}

	Now we introduce morphisms between double categories called \emph{lax double functors} \cite{GP99,Gra20}.
	The definition is given in an `unbiased' way compared to the usual definition of lax functors in the literature.\footnote{%
	This type of definition is closely related to the observation in \cite[Example 3.5]{CS10} that a lax double functor is just a functor
	between \textit{virtual double categories}, restricted to double categories. 
	}
Given two double categories $\dbl{D}$ and $\dbl{E}$, a lax double functor $F\colon\dbl{D}\to\dbl{E}$ consists of
the following data.
\begin{itemize}
	\item%
		Two functors $F_0\colon\dbl{D}_0\to\dbl{E}_0$ and $F_1\colon\dbl{D}_1\to\dbl{E}_1$ commuting with $\src$ and $\tgt$.
		We often omit the subscripts $0,1$ for brevity.
	\item%
		For any path of horizontal arrows $\langle p_1,\ldots,p_n\rangle$, a \textit{coherence} horizontal cell
		\[
			\begin{tikzcd}
				F(A_0)
				\ar[d, equal]
				\sar[r,"F(p_1)"]
				\doublecell[rrrd]{\gamma_{\langle p_1,\ldots,p_n\rangle}}
					&
					F(A_1)
					\sar[r,"F(p_2)"]
						&
						\cdots
						\sar[r,"F(p_n)"]
							&
							F(A_n)
							\ar[d, equal]
				\\
				F(A_0)
				\sar[rrr,"F(p_1\odot p_2\odot\cdots \odot p_n)"']
					&
						&
							&
							F(A_n)
			\end{tikzcd}
		\]
		satisfying the following two conditions.
		\begin{itemize}
			\item
				Suppose we are given a path $\langle\alpha_1,\ldots,\alpha_n\rangle$ of cells, seen as arrows in $\bi{C}_h(\dbl{D})$, and
				let $\langle p_1,\ldots,p_n\rangle$ and $\langle q_1,\ldots,q_n\rangle$ be paths of horizontal arrows that are on the above and the below
				of this sequence of cells, respectively. Then the pair $(\gamma_{\langle p_1,\ldots,p_n\rangle}, \gamma_{\langle q_1,\ldots,q_n\rangle})$
				forms a 2-cell $ F(\alpha_1)\odot\cdots\odot F(\alpha_n)\Rightarrow F(\alpha_1\odot\cdots\odot\alpha_n)$ in $\bi{C}_h(\dbl{D})$.
			\item
				The coherence horizontal cells are closed under composition; i.e.,
				we have $\gamma_{\langle p_1\rangle}=\Id_{p_1}$ and the composite of the following cells is identified with
				$\gamma_{\langle p^1_1,\ldots,p^1_{m_1},p^2_1,\ldots,p^n_1,\ldots,p^n_{m_n}\rangle}$ through the image of the coherence
				between the composite of the path $\langle p^1_1,\ldots,p^1_{m_1},p^2_1,\ldots,p^n_1,\ldots,p^n_{m_n}\rangle$
				and that of $\langle q_1,\ldots,q_{n}\rangle$ under $F$.
				Here, we put $q_i= p^i_1\odot\cdots\odot p^i_{m_i}$ for each $i$.
				\[
					\begin{tikzcd}
						\cdot
						\doublecell[rrd]{\gamma_{p^1_1,\ldots,p^1_{m_1}}}
						\sar[r, "F(p^1_1)"]
						\ar[d, equal]
							&
							\cdots
							\sar[r, "F(p^1_{m_1})"]
								&
								\cdot
								\ar[d, equal]
								\sar[r, "F(p^2_1)"]
								\doublecell[rrd]{\gamma_{p^2_1,\ldots,p^2_{m_2}}}
									&
									\cdot
										&
										\cdots
										\doublecell[rrd]{\gamma_{p^{n-1}_1,\ldots,p^{n-1}_{m_{n-1}}}}
										\doublecell[d]{\cdots}
											&
											\cdot
											\sar[r, "F(p^{n-1}_{m_{n-1}})"]
												&
												\cdot
												\doublecell[rrd]{\gamma_{p^n_1,\ldots,p^n_{m_n}}}
												\ar[d, equal]
												\sar[r, "F(p^n_1)"]
													&
													\cdots
													\sar[r, "F(p^n_{m_n})"]
														&
														\cdot
														\ar[d, equal]
						\\
						\cdot
						\doublecell[rrrrrrrrd]{\gamma_{q_1,\ldots,q_n}}
						\sar[rr, "F(q_1)"']
						\ar[d, equal]
							&
								&
								\cdot
								\sar[rr, "F(q_2)"']
									&
										&
										\cdots
										\sar[rr,"F(q_{n-1})"']
											&
												&
												\cdot
												\sar[rr, "F(q_n)"']
													&
														&
														\cdot
														\ar[d, equal]
						\\
						\cdot
						\sar[rrrrrrrr, "F( q_1\odot\cdots\odot q_{n})"']
							&
								&
									&
										&
											&
												&
													&
														&
														\cdot
					\end{tikzcd}
				\]
		\end{itemize}
\end{itemize}
The last condition for coherence cells ensures that to define a lax double functor, we only have to
obtain $\gamma_{\langle p_1,\ldots,p_n\rangle}$ for the case $n=2$ or $n=0$ and check suitable conditions,
which results in the classical definition of lax functors in \cite{GP99,Gra20}.

An \emph{oplax double functor} $F\colon\dbl{D}\to\dbl{E}$ is a lax double functor $\dbl{D}^\vop\to\dbl{E}^\vop$.
A \emph{pseudo-double functor} (or double functor) between two double categories
is a lax functor whose coherence cell $\gamma_{\langle p_1,\ldots,p_n\rangle}$ is invertible in $\dbl{E}_1$
for each path $\langle p_1,\ldots,p_n\rangle$ of horizontal arrows.
In other words, a double functor is just an internal pseudo-functor in $\CATbi$.

Suppose we have parallel lax double functors $F,G\colon\dbl{D}\to\dbl{E}$.
A \emph{vertical transformation} $\theta\colon F\to G$ is a pair $(\{\theta_A\}_{A\in\dbl{D}_0}, \{\theta_p\}_{p\in\dbl{D}_1})$
of natural transformations $F_0\Rightarrow G_0$ and $F_1\Rightarrow G_1$ compatible with $\src$ and $\tgt$ such that
for each path $\langle p_1,\ldots,p_n\rangle$ of horizontal arrows in $\dbl{D}$, the following equation holds.
\[
	\begin{tikzcd}
		F(A_0)
		\ar[d, equal]
		\sar[r,"F(p_1)"]
		\doublecell[rrrd]{\gamma_{\langle p_1,\ldots,p_n\rangle}}
			&
			F(A_1)
			\sar[r,"F(p_2)"]
				&
				\cdots
				\sar[r,"F(p_n)"]
					&
					F(A_n)
					\ar[d, equal]
		\\
		F(A_0)
		\sar[rrr,"F(p_1\odot p_2\odot\cdots \odot p_n)"']
		\doublecell[rrrd]{\theta_{p_1\odot p_2\odot\cdots \odot p_n}}
		\ar[d, "\theta_{A_0}"']
			&
				&
					&
					F(A_n)
					\ar[d, "\theta_{A_n}"]
		\\
		G(A_0)
		\sar[rrr,"G(p_1\odot p_2\odot\cdots \odot p_n)"']
			&
				&
					&
					G(A_n)
	\end{tikzcd}
	=
	\begin{tikzcd}
		F(A_0)
		\ar[d, "\theta_{A_0}"']
		\sar[r,"F(p_1)"]
		\doublecell[rd]{\theta_{p_1}}
			&
			F(A_1)
			\sar[r,"F(p_2)"]
			\ar[d, "\theta_{A_1}"']
			\doublecell[rd]{\theta_{p_2}}
				&
				\cdots
				\sar[r,"F(p_n)"]
				\doublecell[d]{\cdots}
				\doublecell[rd]{\theta_{p_n}}
					&
					F(A_n)
					\ar[d, "\theta_{A_n}"]
		\\
		G(A_0)
		\ar[d, equal]
		\sar[r,"G(p_1)"]
		\doublecell[rrrd]{\gamma'_{\langle p_1,\ldots,p_n\rangle}}
			&
			G(A_1)
			\sar[r,"G(p_2)"]
				&
				\cdots
				\sar[r,"G(p_n)"]
					&
					G(A_n)
					\ar[d, equal]
		\\
		G(A_0)
		\sar[rrr,"G(p_1\odot p_2\odot\cdots \odot p_n)"']
			&
				&
					&
					G(A_n)
	\end{tikzcd}
\]

We write $\DblCat$ for the 2-category of double categories,
pseudo-double functors, and vertical transformations.

For each 2-cell $\alpha\colon f\Rightarrow g\colon A\to B$ in $\biV{\dbl{D}}$, we obtain another 2-cell $F\alpha\colon Ff\Rightarrow Fg\colon FA\to FB$
in $\biV{\dbl{E}}$ as the following composite.
\[
	\begin{tikzcd}[column sep=small]
			&
			FA
			\ar[ld, equal]
			\ar[rd, equal]
				&
		\\
		FA
		\ar[d, "Ff"']
		\sar[rr, "F(\Id_A)"]
		\doublecell[rrd]{F\alpha}
		\doublecell[rr, shift left=3.5ex]{\gamma_{\langle\rangle}}
			&
				&
				FA
				\ar[d, "Fg"]
		\\
		FB
		\sar[rr, "F(\Id_B)"]
		\doublecell[rr, shift right=3ex]{\gamma^{-1}_{\langle\rangle}}
			&
				&
				FB
		\\
			&
			FB
			\ar[lu, equal]
			\ar[ru, equal]
				&
	\end{tikzcd}
\]
A straightforward discussion shows that this assignment extends to a 2-functor $\biV{F}\colon\biV{\dbl{D}}\to\biV{\dbl{E}}$.
Moreover, a vertical natural transformation $F\Rightarrow G$ also restricts to a 2-natural transformation $\biV{F}\Rightarrow\biV{G}$ and
they define a 2-functor $\bi{V}\colon\DblCat\to 2\Catbi$.

The 2-category $\DblCat$ has finite products,
whose data (objects, vertical arrows, horizontal arrows, and cells)
are pairs of data for each component.
The terminal double category is denoted as $\mathbbm{1}$,
and the product double category of $\dbl{D}$ and $\dbl{E}$ is
denoted as $\dbl{D}\times\dbl{E}$.

We move on to the illustration of several structures on double categories.

\begin{proposition}[{\cite[Theorem 4.1]{Shu08}}]
	\label{prop:repradj}
	Let $\dbl{D}$ be a double category,
	$f\colon X\to Y$ be a vertical arrow in $\dbl{D}$,
	and $p\colon X\sto Y$ and $q\colon Y\sto X$ be horizontal arrows.
	Then, the (structural) 2-out-of-3 condition holds for the following three data;
	i.e., given any two of the three pieces of data, the other is uniquely determined under a suitable ternary relation.
	\begin{enumerate}
		\item
			\underline{Companion.}
			A pair $(\alpha, \beta)$ satisfying the following.
			\begin{equation}
				\label{eq:companion}
					\begin{tikzcd}
							&
							X
							\ar[ld, equal]
							\sar[r, "p"]
							\ar[d, "f"]
							\doublecell[ld, shift left=2ex]{\beta}
								&
								Y
								\ar[ld, equal]
								\doublecell[ld, shift right=2ex]{\alpha}
						\\
						X
						\sar[r, "p"']
							&
							Y
								&
					\end{tikzcd}
					=
					\begin{tikzcd}
						X
						\sar[r, "p"]
						\ar[d, equal]
						\doublecell[rd]{\verteq}
							&
							Y
							\ar[d, equal]
						\\
						X
						\sar[r, "p"']
							&
							Y
					\end{tikzcd}
					\hspace{1ex}
					,
					\hspace{1ex}
					\begin{tikzcd}[column sep = small]
							&
							X
							\ar[rd, "f"]
							\ar[ld, equal]
								&
						\\
						X
						\ar[rd, "f"']
						\sar[rr, "p"{description, near start}]
						\doublecell[rr,  shift left=2.5ex]{\beta}
						\doublecell[rr, shift right=2.5ex]{\alpha}
							&
								&
								Y
								\ar[ld, equal]
						\\
							&
							Y
								&
					\end{tikzcd}
					=
					\begin{tikzcd}
						X
						\ar[d, "f"']
						\sar[r, equal]
						\doublecell[rd]{\Id_f}
							&
							X
							\ar[d, "f"]
						\\
						Y
						\sar[r, equal]
							&
							Y
					\end{tikzcd}
			\end{equation}
			If $f$ and $p$ come equipped with these structures,
			we say that $p$ is a \emph{companion} of $f$.
		\item
			\underline{Conjoint.}
			A pair $(\gamma, \delta)$ satisfying the following.
			\begin{equation}
					\begin{tikzcd}
						Y
						\sar[r, "q"]
						\ar[rd, equal]
						\doublecell[rd, shift left=2ex]{\gamma}
							&
							X
							\ar[d, "f"]
							\ar[rd, equal]
							\doublecell[rd, shift right=2ex]{\delta}
								&
						\\
							&
							Y
							\sar[r, "q"']
								&
								X
					\end{tikzcd}
					=
					\begin{tikzcd}
						Y
						\sar[r, "q"]
						\ar[d, equal]
						\doublecell[rd]{\verteq}
							&
							X
							\ar[d, equal]
						\\
						Y
						\sar[r, "q"']
							&
							X
					\end{tikzcd}
					\hspace{1ex}
					,
					\hspace{1ex}
					\begin{tikzcd}[column sep = small]
							&
							X
							\ar[ld, "f"']
							\ar[rd, equal]
								&
						\\
						Y
						\ar[rd, equal]
						\sar[rr, "q"{description, near start}]
						\doublecell[rr,  shift left=2.5ex]{\delta}
						\doublecell[rr, shift right=2.5ex]{\gamma}
							&
								&
								X
								\ar[ld, "f"]
						\\
							&
							Y
								&
					\end{tikzcd}
					=
					\begin{tikzcd}
						X
						\ar[d, "f"']
						\sar[r, equal]
						\doublecell[rd]{\Id_f}
							&
							X
							\ar[d, "f"]
						\\
						Y
						\sar[r, equal]
							&
							Y
					\end{tikzcd}
			\end{equation}
			If $f$ and $q$ come equipped with these structures,
			we say that $q$ is a \emph{conjoint} of $f$.
		\item
			\underline{Adjoint in $\biH{\dbl{D}}$.}
			A pair $(\eta, \varepsilon)$ satisfying the following.
			\begin{equation}
				\begin{tikzcd}[column sep=small]
						&
						Y
						\sar[r, "q"]
						\ar[ld, equal]
						\doublecell[d]{\verteq}
							&
							X
							\ar[ld, equal]
							\ar[rd, equal]
								&
					\\
					Y
					\sar[r, "q"]
					\doublecell[rr, shift right=3ex]{\varepsilon}
						&
						X
						\sar[r, "p"{description, near start, inner sep=0.1mm}]
						\doublecell[rr, shift left=3ex]{\eta}
							&
							Y
							\sar[r, "q"']
							\doublecell[d]{\verteq}
								&
								X
								\ar[ld, equal]
					\\
						&
						Y
						\sar[r, "q"']
						\ar[lu, equal]
						\ar[ru, equal]
							&
							X
								&
				\end{tikzcd}
				=
				\begin{tikzcd}
					Y
					\sar[r, "q"]
					\ar[d, equal]
					\doublecell[rd]{\verteq}
						&
						X
						\ar[d, equal]
					\\
					Y
					\sar[r, "q"']
						&
						X
				\end{tikzcd}
				\hspace{1ex}
				,
				\hspace{1ex}
				\begin{tikzcd}[column sep=small]
						&
						X
						\sar[r, "p"]
						\ar[ld, equal]
						\ar[rd, equal]
							&
							Y
							\ar[rd, equal]
							\doublecell[d]{\verteq}
								&
					\\
					X
					\sar[r, "p"']
					\doublecell[rr, shift left=3ex]{\eta}
						&
						Y
						\sar[r, "q"{description, near start, inner sep=0.1mm}]
						\doublecell[rr, shift right=3ex]{\varepsilon}
							&
							X
							\sar[r, "p"]
								&
								Y
					\\
						&
						X
						\ar[lu, equal]
						\doublecell[u]{\verteq}
						\sar[r, "p"']
							&
							Y
							\ar[lu, equal]
							\ar[ru, equal]
								&
				\end{tikzcd}
				=
				\begin{tikzcd}
					X
					\sar[r, "p"]
					\ar[d, equal]
					\doublecell[rd]{\verteq}
						&
						Y
						\ar[d, equal]
					\\
					X
					\sar[r, "p"']
						&
						Y
				\end{tikzcd}
			\end{equation}
	\end{enumerate}
	In particular, a vertical arrow with companion and conjoint produces
	an adjoint in $\biH{\dbl{D}}$.
	We call such an adjoint a \emph{representable} adjoint.
\end{proposition}

\begin{definition}
	\label{def:equip}
	A double category $\dbl{D}$ is an \emph{equipment}
	if the functor $\langle \src,\tgt\rangle\colon \dbl{D}_1\to \dbl{D}_0\times\dbl{D}_0$
	is a fibration.
\end{definition}
Equipments are also known as `framed bicategories' \cite{Shu08} 
and `fibrant double categories' \cite{Ale18}.

\begin{remark}
	\label{rem:equip}
	A double category $\dbl{D}$ is an equipment if and only if
	$\langle \src,\tgt\rangle$ is an opfibration,
	hence a bifibration. 
	Also, being an equipment is equivalent to the condition that
	for every vertical arrow $f\colon X\to Y$, there are horizontal arrows $p\colon X\sto Y$ and $q\colon Y\sto X$,
	equipped with two (hence all) of the data listed in \cref{prop:repradj};
	see \cite[Theorem 4.1]{Shu08}. 
	Under this correspondence,
	$\alpha$ in \Cref{eq:companion} is the cartesian lifting of $(f\colon X\to Y,\id\colon Y\to Y)$,
	and likewise for other cells.
	The companion and conjoint of $f\colon X\to Y$ are written as $f_!$ and $f^*$.
	Note that the vertical composition of two cartesian cells is cartesian,
	and the vertical composition of two opcartesian cells is opcartesian.

	By a \emph{cartesian (resp. opcartesian) cell},
	we mean a cartesian (resp. opcartesian) morphism of 
	the bifibration $\langle \src,\tgt\rangle$.
	From a horizontal arrow $p\colon Y\sto Z$ and vertical arrows $f\colon W \to Y$ and $g\colon X\to Z$,
	the cartesian lift of $(f,g)$ to $p$ in the bifibration gives the cartesian cell as the cell on the left below.
	\[
		\begin{tikzcd}
			W
			\ar[d, "f"']
			\sar[r, "{p(f,g)}"]
			\doublecell[rd]{\cart}
				&
				X
				\ar[d, "g"]
			\\
			Y
			\sar[r, "p"']
				&
				Z
		\end{tikzcd}
		\hspace{2ex}
		,
		\hspace{2ex}
		\begin{tikzcd}
			W
			\ar[d, "f"']
			\sar[r, "q"]
			\doublecell[rd]{\opcart}
				&
				X
				\ar[d, "g"]
			\\
			Y
			\sar[r, "{\Ext(q;f,g)}"']
				&
				Z
		\end{tikzcd} 
		\hspace{2ex}
		,
		\hspace{2ex}
		\begin{tikzcd}[column sep=small]
			W
			\ar[rd, "f"']
			\sar[rr, "{Z(f,g)}"]
			\doublecell[rr,shift right = 2ex]{\cart} 			
			&
			&
				X
				\ar[ld, "g"]
			\\
				&
				Z
		\end{tikzcd}
		\hspace{1ex}
		,
		\hspace{1ex}
		\begin{tikzcd}[column sep=small]
			&
			X 
			\ar[ld, "f"']
			\ar[rd, "g"]
			&
			\\
			Y 
			\sar[rr, "{\Ext(X;f,g)}"']
			\doublecell[rr,shift left = 2ex]{\opcart}
			&
			&
			Z
		\end{tikzcd}
	\]
	Here the cartesian cell is unique up to invertible horizontal cells,
	so we just write $\cart$ for the cartesian cell and call the horizontal arrow $p(f,g)$ the \emph{restriction} of $p$ along $f$ and $g$.
	Taking the horizontal dual, the opcartesian cell is unique up to invertible horizontal cells,
	so we just write $\opcart$ for the opcartesian cell and call the horizontal arrow $\Ext(q;f,g)$ the \emph{extension} of $q$ along $f$ and $g$.
	In particular, as presented in the right half of the above diagrams,
	the restriction of $\Id_Z$ through $f$ and $g$ is written as $Z(f,g)$,
	and the extension of $\Id_X$ through $f$ and $g$ is written as $\Ext(X;f,g)$ for brevity.

	The restriction $p(f,g)$ and the extension $\Ext(q;f,g)$ are realised as $f_!pg^*$ and $f^*qg_!$,
	respectively, using the companion and conjoint,
	and the cartesian cell and the opcartesian cell are realised as below.
	\[
		\begin{tikzcd}
			W
			\ar[rd, "f"']
			\sar[r, "f_!"]
			\doublecell[rd, shift left=2.0ex]{\alpha}
				&
				Y
				\doublecell[rd]{\verteq}
				\sar[r,"p"]
				\ar[d, equal]
					&
					Z
					\ar[d, equal]
					\sar[r, "g^*"]
						&
						X
						\doublecell[ld, shift right=2.0ex]{\gamma}
						\ar[ld, "g"]
			\\
				&
				Y
				\sar[r, "p"']
					&
					Z
						&
		\end{tikzcd}
		\hspace{2ex}
		,
		\hspace{2ex}
		\begin{tikzcd}
			&
			W
			\ar[ld, "f"']
			\sar[r, "q"]
			\ar[d, equal]
			\doublecell[rd]{\verteq}
				&
				X
				\ar[rd, "g"]
				\ar[d, equal]
			\\
			Y
			\sar[r, "f^*"']
			\doublecell[r, shift left=2.0ex,xshift =1ex]{\delta}
				&
				W 
				\sar[r, "q"']
					&
					X
					\sar[r, "g_!"']
					\doublecell[r, shift left=2.0ex,xshift=-1ex]{\beta}
						&
						Z
		\end{tikzcd}
	\]
	Put it another way,
	if we are given a cartesian cell $\alpha$ and $\gamma$ as above,
	then the above cell is the restriction of $p$ through $f$ and $g$.
	Since the $\alpha$ and $\gamma$ are cartesian cells and the $\beta$ and $\delta$ are opcartesian cells,
	we just write $\cart$ and $\opcart$ for them as well.
	For a comprehensive treatment on equipments, see \cite[\S 4]{Shu08}.
\end{remark}

Here, we introduce a convenient lemma that we shall refer to as the `sandwich lemma' throughout this paper.

\begin{lemma}
	\label{lem:Sandwich}
	Let $\dbl{D}$ be an equipment.
	Given a sequence of horizontally composable cells
	\begin{equation}
		\label{dgm:cartopcartcart}
	\begin{tikzcd}
		\cdot 
		\sar[r] 
		\ar[d] 
		\doublecell{\cart}
		&
		\cdot 
		\sar[r] 
		\ar[d] 
		\doublecell{\opcart}
		& 
		\cdot 
		\sar[r] 
		\ar[d] 
		\doublecell{\cart}
		& 
		\cdot 
		\ar[d] 
		\\
		\cdot 
		\sar[r] 
		& 
		\cdot 
		\sar[r] 
		& 
		\cdot 
		\sar[r] 
		& 
		\cdot 
	\end{tikzcd}
	\end{equation}
	with the opcartesian cell sandwiched between two cartesian cells,
	the composition of these cells is cartesian.
	The same thing holds when swapping the roles of `cartesian' and `opcartesian'.
\end{lemma}
\begin{proof}
	By \cref{rem:equip},
	we can rewrite the diagram \Cref{dgm:cartopcartcart} as follows,
	in which the names of the cells correspond to that in \cref{prop:repradj}.
	\[
		\begin{tikzcd}
			\cdot
			\sar[r,"f_!"]
			\ar[rd,"f"']
			\doublecell[rd,  shift left=2.0ex]{\alpha}
				&
				\cdot
				\sar[r]
				\doublecell[rd]{\verteq}
				\ar[d, equal]
					&
					\cdot
					\sar[r]
					\ar[d, equal]
						&
						\cdot
						\sar[r]
						\ar[d, equal]
						\ar[ld,]
						\doublecell[ld, shift right=2.0ex]{\gamma}
						\doublecell[ld,  shift left=2.0ex]{\delta}
						\doublecell[rd]{\verteq}
							&
							\cdot
							\ar[rd,]
							\sar[r,]
							\doublecell[rd,  shift left=2.0ex]{\alpha}
							\doublecell[rd, shift right=2.0ex]{\beta}
							\ar[d, equal]
								&
								\cdot
								\doublecell[rd]{\verteq}
								\sar[r]
								\ar[d, equal]
									&
									\cdot
									\ar[d, equal]
									\sar[r,"k^*"]
										&
										\cdot
										\doublecell[ld, shift right=2.0ex]{\gamma}
										\ar[ld,"k"]
			\\
				&
				\cdot
				\sar[r]
					&
					\cdot
					\sar[r]
						&
						\cdot
						\sar[r]
							&
							\cdot
							\sar[r]
								&
								\cdot
								\sar[r]
									&
									\cdot
										&
		\end{tikzcd}
	\]	
	Because of the equalities described in \cref{prop:repradj},
	the middle sequence of square cells are all identities.
	Again by \cref{prop:repradj}, this implies that the composition of the cells in the diagram is cartesian.
	Considering the same statement for the vertical opposite of $\dbl{D}$, we obtain the dual.
\end{proof}

Now we define cartesian double categories as cartesian objects in the 2-category $\DblCat$. 
\begin{definition}[{\cite[Definition 4.2.2]{Law15},\cite[Definition 4.2.1]{Ale18}}]
	\label{def:cartesian}
	A double category $\dbl{D}$ is \emph{cartesian}
	if $\Delta \colon \dbl{D}\to\dbl{D}\times\dbl{D}$ and $!\colon \dbl{D}\to\mathbbm{1}$
	have the right adjoints $\times\colon \dbl{D}\times\dbl{D}\to\dbl{D}$ and 
	$1\colon \mathbbm{1}\to\dbl{D}$ in $\DblCat$.
	We say $1$ is the (vertical) terminal object in $\dbl{D}$.
\end{definition}

\begin{remark}
	If $\dbl{D}$ has the terminal object $1$,
	it serves as the terminal object in $\dbl{D}_0$, and the horizontal identity $\Id_1$ is the terminal object in $\dbl{D}_1$.
	In particular, for each horizontal arrow $p\colon A\sto B$ in $\dbl{D}$, there exists a unique cell of the following form.
	\[
		\begin{tikzcd}[column sep=small]
			A
			\sar[rr, "p"]
			\ar[rd, "!"']
			\doublecell[rr, shift right=3ex]{!}
				&
					&
					B
					\ar[ld, "!"]
			\\
				&
				1
					&
		\end{tikzcd}
		\vspace{-2.5ex}
	\]
\end{remark}

\begin{remark}
	\label{rem:localProduct} 
	Comprehensive discussions on cartesian double categories can 
	be found in \cite{Ale18}.
	Let us write $\bi{E}quip$ for the full sub-2-category of $\DblCat$ consisting of equipments
	and $\bi{F}ib$ for the 2-category of fibrations and morphisms of fibrations.
	Since we have the forgetful 2-functors $\oneV{\mhyphen}\colon\bi{E}quip\to\Catbi$ 
	and $\bi{E}quip\to\bi{F}ib$ that assigns $\langle \src,\tgt\rangle$ in \Cref{def:equip} to each equipment,
	$\oneV{\dbl{D}}$ has finite products 
	and $\biH{\dbl{D}}$ locally has finite products, if $\dbl{D}$ is cartesian \cite[Theorem 4.2.2]{Ale18}.
	Conversely, for an equipment $\dbl{D}$,
	if $\oneV{\dbl{D}}$ has finite products and $\biH{\dbl{D}}$ locally has finite products,
	then the functors above have lax right adjoints, but not pseudo in general \cite[Theorem 4.1.2, Proposition 4.2.3]{Ale18}.
\end{remark}
\begin{remark}
	\label{rem:biVPreservesCartObj}
	Suppose we are given parallel horizontal arrows $p,q\colon A \sto B$. Then the product $p\land q\colon A\sto B$ and
	the terminal object $\top\colon A \sto B$
	in the hom-category of $\biH{\dbl{D}}$ are defined by the following restrictions in $\dbl{D}$.
	\[
		\begin{tikzcd}
			A
			\sar[r, "p\land q"]
			\ar[d, "\Delta"']
			\doublecell[rd]{\cart}
				&
				B
				\ar[d, "\Delta"]
			\\
			A\times A
			\sar[r, "p\times q"']
				&
				B\times B
		\end{tikzcd}
		\hspace{3ex}
		,
		\hspace{3ex}
		\begin{tikzcd}[column sep=small]
			A
			\sar[rr, "\top"]
			\ar[rd, "!"']
			\doublecell[rr, shift right=3ex]{\cart}
				&
					&
					B
					\ar[ld, "!"]
			\\
				&
				1
					&
		\end{tikzcd}
	\]
	Moreover, we have already constructed another 2-functor $\bi{V}\colon\DblCat\to 2\Catbi$ that preserves finite products.
	Therefore, for each cartesian double category $\dbl{D}$, the vertical 2-category $\biV{\dbl{D}}$ has finite products.
\end{remark}

\begin{remark}
	A cartesian double category $\dbl{D}$ is naturally seen as a pseudo-monoid object in $\DblCat$.
	Since every invertible vertical transformation between equipments gives rise to a \textit{horizontal transformation} (in the sense of \cite{GP99}\footnote{%
		This horizontal transformation is weaker than what we obtain through merely rephrasing vertical to horizontal.
		Indeed, our definition of vertical transformation needs to be natural when restricted to the vertical category, while
		the horizontal transformation here shall be pseudo-natural when restricted to the horizontal bicategory.
		For more description of weaker versions of the notions for double categories, see \cite{Ver11}.
	})
	through taking companions for each component, $\biH{\dbl{D}}$ and $\bi{C}_h(\dbl{D})$ have structures like those for monoidal bicategories.
	See \cite[Theorem 5.1]{Shu10} for more details.
\end{remark}

\begin{definition}[{cf.\ \cite{GP99}}]
	Let $\dbl{D}$ be a double category.
	A \emph{tabulator} of a horizontal arrow $p\colon X\sto Y$
	is the representing object of the functor 
	$\dbl{D}_1(\Id(-),p)\colon \dbl{D}_0^\op\to\Set$
	equipped with the universal cell as follows.
	\begin{equation}
		\label{eq:tabu}
		\begin{tikzcd}[column sep=small]
				&
				\top p
				\ar[ld, "l"']
				\ar[rd, "r"]
					&
			\\
			X
			\sar[rr, "p"']
			\doublecell[rr, shift left=3ex]{\tau}
				&
					&
					Y
		\end{tikzcd}
	\end{equation}	
	In this case, we say that $\tau$ is a \emph{tabulating cell} of $p$. 
	Therefore, $\dbl{D}$ \emph{has tabulators} if $\Id\colon \dbl{D}_0\to\dbl{D}
	_1$ has a right adjoint $\top\colon \dbl{D}_1\to\dbl{D}_0$.
	A tabulator \cref{eq:tabu} is \emph{strong}
	if $\tau$ is an opcartesian cell in $\dbl{D}$.
	A tabulator \cref{eq:tabu} is \emph{left-sided} 
	if $Y$ is the terminal object,
	and $\dbl{D}$ \emph{has left-sided tabulators} if 
	every horizontal arrow to the terminal object has a tabulator, which is left-sided.
\end{definition}

By a `tabulator' of a horizontal arrow $p$ in the definition,
we occasionally mean the span $(l,r)$ or the vertical arrow 
$\langle l,r\rangle\colon \top p\to X\times Y$
instead of their domain $\top p$.

\begin{remark}
	A tabulator defined in this way only has the one-dimensional
	universality.
	This is an example of one-dimensional limits in terms of \cite{GP99}.
	To modify this,
	we may extend $\Id\colon\dbl{D}_0\to\dbl{D}_1$ to a pseudo-functor $\dbl{D}\to\dbl{H}(\dbl{D})$,
	where $\dbl{H}(\dbl{D})$ is the double category of horizontal arrows in $\dbl{D}$.
	Important though this kind of universality is, 
	it is not necessary to consider it for our aim.
	See \cite[5.3]{GP99} for more details.
\end{remark}

\begin{example}
	\label{eg:Prof}
	Let $\Prof$ be the double category of small categories,
	functors, and profunctors. Here, we define
	a profunctor $p\colon A\sto B$ as a functor $A^\op\times B\to \Set$.
	The tabulator of a profunctor $p\colon A\sto B$ is called the `category of elements'
	in \cite{GP99}
	and `two-sided discrete fibration' in \cite{LR20}.
	One can readily check that this tabulator is strong.
\end{example}
\begin{example}
	Two-sided split fibrations give an `augmented virtual double category'
	$\spFib$ (see \cite[Example 2.11]{Kou22}).
	The notion of (strong) tabulator is defined similarly in this generalised setting,
	and a tabulator of a horizontal arrow (=two-sided split fibration)
	is itself seen as a pair of functors. It is not strong in general.
\end{example}

\subsection{Orthogonal factorisation systems}\label{subsec:OFS}

In this subsection, we review some basic notations surrounding \textit{orthogonal factorisation systems} \cite{FK72}.

\begin{definition}
	An \ac{OFS} on a category $\one{C}$ consists of 
	a pair $(\zero{E},\zero{M})$ of classes of arrows in $\one{C}$
	that satisfies the following conditions:
	\begin{enumerate}
		\item
			$\zero{E}$ and $\zero{M}$ are closed under composition and contain isomorphisms.
		\item
			Every arrow $e\colon X\to Y$ in $\zero{E}$ is left orthogonal to every arrow $m\colon A\to B$ in $\zero{M}$; that is, 
			every commutative square
			\[
				\begin{tikzcd}
					X 
					\ar[r]
					\ar[d,"e"',two heads]
						&
						A
						\ar[d,"m",tail]
					\\
					Y
					\ar[r]
					\ar[ru,"f",dotted]
						&
						B
				\end{tikzcd}	
			\]
			has a unique diagonal filler $f\colon Y\to A$ that makes two triangles commute.
		\item
			Every arrow $f$ in $\one{C}$ factors as $f= e\fatsemi m$ where $e$ belongs to $\zero{E}$ and $m$ belongs to $\zero{M}$.
	\end{enumerate}
	$\zero{E}$ and $\zero{M}$ are called \emph{the left class} and \emph{the right class} of the \ac{OFS}, respectively.
	For a factorisation of $f$ as $f=e\fatsemi m$, we say $m$ is the $\zero{M}$-\emph{image} of $f$ if $e\in\zero{E}$ and $m\in\zero{M}$.
	
	An \ac{OFS} $(\zero{E},\zero{M})$ is a \emph{\acf{SOFS}} if $\zero{E}$ is stable under pullback.
	We also say it is \emph{right-proper} if $\zero{M}$ is a subclass of the class of all monomorphisms,
	\emph{left-proper} if $\zero{E}$ is a subclass of the class of epimorphisms,
	and \emph{proper} if it is both right-proper and left-proper.
	Furthermore, we call it \emph{anti-right-proper} 
	if the class of monomorphisms is a subclass of $\zero{M}$.
	We only treat \aclp{OFS},
	and so we omit the adjective `orthogonal' in the sequel.
\end{definition}

An inspiring example of an \acl{OFS} is $(\Regepi,\Mono)$ in a regular category,
where $\Mono$ is the class of monomorphisms and $\Regepi$ is the class of regular epimorphisms.
It is well known that a category is regular if and only if it has a \acl{SOFS} $(\Regepi,\Mono)$,
see \cite[Scholium 1.3.5]{Joh02} for example.
We call this type of \aclp{OFS} \emph{regular}.

We often draw an arrow in the left class as $\thto{}$,
and an arrow in the right class as $\tto{}$.
For a class $\zero{S}$ of arrows, 
we write $\lclass{\zero{S}}$ for the class of arrows left orthogonal to all members in $\zero{S}$,
and $\rclass{\zero{S}}$ for the class of arrows right orthogonal to all members in $\zero{S}$.
Note that if a category $\one{C}$ has equalisers, then anti-right-properness implies left-properness.

\begin{definition}[{\cite{HNT20}}]
	\label{defn:StableSystem}
	Let $\one{C}$ be a category.
	A \textit{stable system} on $\one{C}$ is a
	class $\zero{M}$ of arrows
	in $\one{C}$ that is stable under composition and \textit{pullback},
	that is, $\zero{M}$ satisfying the following.
	\begin{enumerate}
		\item
			Every isomorphism is in $\zero{M}$.
		\item
			Given a composable pair of arrows $\cdot\to["f"]\cdot\to["g"]\cdot$,
			the composite $f\fatsemi g$ is in $\zero{M}$
			whenever so are both $f$ and $g$.
		\item
			Given a pullback square below in $\one{C}$,
			if $f$ is in $\zero{M}$ then so is $f'$.
			\begin{equation}
				\begin{tikzcd}
					\cdot
					\ar[r]
					\ar[d, "f'"']
					\pullback
						&
						\cdot
						\ar[d, "f"]
					\\
					\cdot
					\ar[r]
						&
						\cdot
				\end{tikzcd}
			\end{equation}
	\end{enumerate}
\end{definition}

We end this section with two classical results on \aclp{OFS}.
The first one, which connects \aclp{OFS} on a category $\one{C}$ to subcategories of the arrow category $\one{C}^\rightarrow$,
appears in \cite{IK86} and is a special case of \cite[Theorem 9]{BG16}.
We also refer to \cite{BK17,HJ03} for other characterisations of \aclp{OFS}
purely in terms of bifibrations (or hyperdoctrines)
rather than stable systems as extra structures.
The second argues that the stability conditions for an \ac{OFS} can be specified in terms of the theory of fibrations.
This result are according to \cite{CJKP97,HJ03}.

\begin{proposition}
	\label{prop:OFSfromStableSystem}
	Let $\one{C}$ be a category with pullbacks
	and a stable system $\zero{M}$.
	By $\one{M}$, we mean the full subcategory of $\one{C}^{\rightarrow}$
	consisting of arrows in $\zero{M}$.
	Note that $\one{M}\to[hook]\one{C}^{\rightarrow}\to\one{C}$
	is a subfibration of the codomain fibration $\cod^\one{C}\colon\one{C}^\rightarrow\to\one{C}$,
	for which we also write $\one{M}$.
	The following are equivalent.
	\begin{enumerate}
		\item
			$(\lclass{\zero{M}}, \zero{M})$ is an orthogonal factorisation system on $\one{C}$.
		\item
			$\one{M}\to[hook]\cod^\one{C}$ has a left adjoint in
			$\Catbi/\one{C}$.
	\end{enumerate}
\end{proposition}
\begin{proof}
	\begin{proof}[i) $\Rightarrow$ ii)]
	Let $C$ be an object in $\one{C}$.
	We show that, for each $f\in\one{C}/C$, there exists
	$\im(f)$ in the fibre $\one{M}_C$ of $\one{M}$ over $C$ and an arrow $\eta_f\colon f\to\im(f)$ in $\one{C}/C$
	such that, for each $m\in\one{M}$ and $(u, v)\colon f\to m$ in
	$\one{C}^\rightarrow$, there exists a unique
	$\bar{u}$ that satisfies the following condition.
	\begin{equation}
		\label{dgm:imIsLeftAdjoint}
		\parbox{\dimexpr\linewidth-4em}{
		$(\bar{u},v)\colon \im(f)\to m$ is an arrow in $\one{C}^\rightarrow$,
		and the following diagram commutes in $\one{C}^\rightarrow$:
		\[
			\begin{tikzcd}[ampersand replacement=\&]
				f
				\ar[r, "\eta_f"]
				\ar[rd, "\mmbox{(u,v)}"']
					\&
					\im(f)
					\ar[d, "\mmbox{(\bar{u},v)}"]
				\\
					\&
					m
			\end{tikzcd}
			\hspace{10ex}
		\]
		}
	\end{equation}
	Let $f\in\one{C}/C$ be an arrow.
	Take a factorisation $f=\eta_f\fatsemi \im(f)$ with
	$\eta_f\in\lclass{\zero{M}}$ and $\im(f)\in\zero{M}$.
	For each $m\in\zero{M}$ and an arrow $(u,v)\colon f\to m$ in $\one{C}^\rightarrow$,
	the square below commutes, and
	since $\eta_f\in\lclass{\zero{M}}$,
	there exists a unique diagonal filler $\bar{u}$.
	\begin{equation}
		\begin{tikzcd}
			\cdot
			\ar[rr, "u"]
			\ar[d,"\eta_f"', two heads]
				&
					&
					\cdot
					\ar[d, "m", tail]
			\\
			\cdot
			\ar[r, "\im(f)"', tail]
			\ar[rru, "\bar{u}", dotted]
				&
				\cdot
				\ar[r,"v"']
					&
					\cdot
		\end{tikzcd}
	\end{equation}
	But this is precisely the same as \cref{dgm:imIsLeftAdjoint}.
	\end{proof}
	\begin{proof}[ii) $\Rightarrow$ i)]
		Let $\im\colon \one{C}^\rightarrow\to \one{M}$ be the left adjoint of
		the inclusion $i\colon \one{M}\to[hook]\one{C}^\rightarrow$ in $\Catbi/\one{C}$.
		We write $\eta\colon \id_{\one{C}^\rightarrow}\to \im\fatsemi i$
		for the unit of this adjoint.
		It suffices to show the following.
		\begin{equation}
			\label{eqn:RightClassOfUnits}
			\rclass{\{\eta_f\,|\,f\in\one{C}^\rightarrow\}}=\zero{M}
		\end{equation}
		Indeed, if this holds,
		we obtain an orthogonal factorisation system
		$(\lclass{\zero{M}},\zero{M})$ since $f=\eta_f\fatsemi \im_f$ gives
		the factorisation of each map $f\in\one{C}^\rightarrow$.

			Firstly, we show $\zero{M}\subseteq\rclass{\{\eta_f\,|\,f\in\one{C}^\rightarrow\}}$.
			Fix $m\colon X\to Y\in\zero{M}$ and $f\colon A\to B$ in $\one{C}$, and we show $m$ is right orthogonal to $\eta_f$.
			In other words, it suffices to show that, for each $(u, v)\colon \eta_f \to m$, there exists a unique diagonal filler $k\colon \dom(\im(f))\to X$, satisfying
			$\eta_f\fatsemi  k=u$ and $k\fatsemi m=v$.
			However, we need only to verify this 
			for the case when we have $\dom(\im(f))=Y$ and $v=\id$
			since $\zero{M}$ is stable under pullback.
			To put it another way, we need only to show that
			the function $\one{C}/Y(\id_Y, m)\to \one{C}/Y(\eta_f, m)$ obtained
			by precomposing the unique map $\eta_f\to \id_Y$ in $\one{C}/Y$
			is a bijection.
			We observe that this bijection is obtained as the following composite of bijections.
			\footnotesize
			\[
				\begin{split}
					\one{C}/Y(\id_Y, m)
					\cong
					\left(\one{C}/B\right)\!\Big/ {\im(f)}\!%
					\left(\btk\im(f)\ar[d, "\id"']\\\im(f)\etk,%
					\btk m\fatsemi \im(f)\ar[d, "m"']\\\im(f)\etk\right)
					\cong
					\one{M}_B/{\im(f)}\!%
					\left(\btk\im(f)\ar[d, "\id"']\\\im(f)\etk,%
					\btk m\fatsemi \im(f)\ar[d, "m"']\\\im(f)\etk\right)\!
					\\
					\cong 
					\left(\one{C}/B\right)\!\Big/{\im(f)}\!%
					\left(\btk f\ar[d, "\eta_f"']\\\im(f)\etk,%
					\btk m\fatsemi \im(f)\ar[d, "m"']\\\im(f)\etk\right)
					\cong
					\one{C}/Y(\eta_Y, m)
				\end{split}
			\]\normalsize
			The first and the last bijections follow from
			the canonical equivalence $\one{C}/Y\cong(\one{C}/B)/\im(f)$.
			The second bijection follows from the fully-faithful inclusion
			$\one{M}_B\to[hook]\one{C}/B$.
			The third is the transport for the sliced adjunction
			(see, for example, \cite{nLa})
			$\one{M}_B/\im(f)\adjoint(\one{C}/B)/i\im(f)$
			induced from the adjunction
			$i\colon \one{M}_B\adjoint\one{C}/B\colon \im$,
			which sends a map
			$u\colon \langle\im(f), \id\rangle\to\langle m\fatsemi \im(f), m\rangle$
			to $\eta_f\fatsemi u$.

			Secondly, we show the converse
			$\rclass{\{\eta_f\,|\,f\in\one{C}^\rightarrow\}}\subseteq\zero{M}$.
			Let $g$ be an arrow right orthogonal to units.
			In particular, there is a unique diagonal filler
			$k\colon \dom(\im(g))\to\dom(g)$ for
			the map $(\id, \im(g))\colon \eta_g\to g$ in $\one{C}^{\rightarrow}$,
			which gives a retraction of $\eta_g$.
			On the other hand, $k\fatsemi \eta_g$ is a diagonal filler
			for the map $(\eta_g, \im(g))\colon \eta_g\to\im(g)$,
			as shown below:
			\[
				\begin{tikzcd}[column sep=small, row sep=small]
					\cdot
					\ar[rr, "\eta_g"]
					\ar[dd, "\eta_g"']
					\ar[rd, equal]
						&
							&
							\cdot
							\ar[dd, "\im(g)"]
					\\
						&
						\cdot
						\ar[ur, "\eta_g"]
						\ar[rd, "g"description]
							&
					\\
					\cdot
					\ar[ur, "k"]
					\ar[rr, "\im(g)"']
						&
							&
							\cdot
				\end{tikzcd}
			\]
			Therefore,
			$k\fatsemi \eta_g$ coincides with $\id$ for the uniqueness
			of diagonal filler.
			\vspace{-4ex}
			\color{white}
	\end{proof}
			\color{black}
\end{proof}

\begin{fact}[{\cite[2.12]{CJKP97}, \cite[Corollary 3.9]{HJ03}}]
	Let $(\zero{E}, \zero{M})$ be an orthogonal factorisation system
	and write $\one{M}$ and $\im$ for the corresponding subfibration
	and the reflection discussed
	in the above proposition. The following are equivalent.
	\begin{enumerate}
		\item
			$(\zero{E}, \zero{M})$ is a \acl{SOFS}.
		\item
			The reflection $\im\colon \cod^{\one{C}}\to\one{M}$ is a fibred functor.
		\item
			$\zero{M}$ is strongly \acs{BC} as a bifibration in the sense of \cite{Shu08}.
	\end{enumerate}
\end{fact}
Here, a bifibration is defined to be \textit{strongly \acs{BC}} if
the Beck-Chevalley condition holds for all pullback squares in the base category as a bifibration,
according to \cite[Definition 13.21]{Shu08}.
This condition provides theoretical support for our construction of \textit{double categories of relations}, as we will see later in \Cref{lem:FrConstruction}.
The paper distinguishes the notion of strongly \acs{BC} bifibrations from a weaker notion of \acs{BC} bifibrations
called \textit{weakly \acs{BC}} bifibrations,
in which the Beck-Chevalley condition holds only for pullback squares one of whose legs is a product projection.

\section{Axiomatising double categories of relations}\label{sec:axiomatising}
\subsection{Beck-Chevalley pullbacks}\label{sec:BCpb}

To begin with, we recall from \cite{WW08} the notion of the \textit{Beck-Chevalley condition}, interpreted in terms of double categories.
\begin{definition}[{cf.\ \cite[2.4]{WW08}}]
	\label{defn:BeckChevalley}
	Let $\dbl{D}$ be a double category.
	A \emph{diamond cell} in $\dbl{D}$ is a quadruple of vertical arrows together with a vertical cell $\alpha$ of the form on the left below.
	A diamond cell is called an \emph{identity diamond cell} if the vertical cell is the identity cell.
	We say a diamond cell
	satisfies the
	\emph{Beck-Chevalley condition}
	if there exists a horizontal arrow $p\colon B\sto C$
	and $\alpha$ factors as an opcartesian cell followed by a cartesian cell
	as shown in the right below:
	\begin{equation}
		\begin{tikzcd}[column sep=small]
				&
				A
				\ar[ld, "g"']
				\ar[rd, "f"]
				\doublecell[dd]{\alpha}
					&
			\\
			B
			\ar[rd, "h"']
				&
					&
					C
					\ar[ld, "k"]
			\\
				&
				D
					&
		\end{tikzcd}
		\hspace{1ex}
		=
		\hspace{1ex}
		\begin{tikzcd}[column sep=small]
				&
				A
				\ar[ld, "g"']
				\ar[rd, "f"]
					&
			\\
			B
			\sar[rr, "p"{near end, description, inner sep=0.1mm}]
			\ar[rd, "h"']
			\doublecell[rr, shift left=2.2ex]{\opcart}
			\doublecell[rr, shift right=2.2ex]{\cart}
				&
					&
					C
					\ar[ld, "k"]
			\\
				&
				D
					&
		\end{tikzcd}
	\end{equation}
	Although this condition is defined for a diamond cell,
	we often abuse the terminology and say that a cell $\alpha$ satisfies the Beck-Chevalley condition when the quadruple of vertical arrows is evidently recognised from the context.

\end{definition}

We explain how this concept coincides with that in \cite{WW08,LWW10}.
Take a cartesian bicategory $\bi{B}$ whose subbicategory of maps $\bi{M}\subset\bi{B}$ is a 2-category,
and we can construct an equipment $\dbl{D}$ satisfying $\bi{B}=\biH{\dbl{D}}$ and $\bi{M}=\biV{\dbl{D}}$ (see \cite{Ver11,Ale18} for more details).
Then the notion of the Beck-Chevalley condition for this equipment is the same as what is defined in \cite{WW08}, and further investigated in \cite{LWW10}.

However, the following definitions for \textit{pullback squares} slightly differ from those in the literature. This is because they deal with \textit{bipullbacks} in
the bicategory $\bi{M}=\biV{\dbl{D}}$, while we consider pullbacks in the category $\oneV{\dbl{D}}$.

\begin{definition}
	A \emph{Beck-Chevalley pullback square} in $\dbl{D}$ is a pullback square in $\oneV{\dbl{D}}$
	as presented on the left below for which the two identity diamond cells placed in both directions
	as in the diagrams in the middle and on the right below
	satisfy the Beck-Chevalley condition.
	\[
	\begin{tikzcd}
		P 
		\ar[r, "s"]
		\ar[d, "t"']
		\pullback[rd]
		&
		A
		\ar[d, "f"]
		\\
		B
		\ar[r, "g"']
		&
		C
	\end{tikzcd}	
	\hspace{2ex}
	,
	\hspace{2ex}
	\begin{tikzcd}[column sep=small]
		&
		P
		\ar[ld, "s"']
		\ar[rd, "t"]
		\doublecell[dd]{\Id}
			&
	\\
	A
	\ar[rd, "f"']
		&
			&
			B
			\ar[ld, "g"]
	\\
		&
		C
			&
	\end{tikzcd}
	\hspace{2ex}
	,
	\hspace{2ex}
	\begin{tikzcd}[column sep=small]
		&
		P
		\ar[ld, "t"']
		\ar[rd, "s"]
		\doublecell[dd]{\Id}
			&
	\\
	B
	\ar[rd, "g"']
		&
			&
			A
			\ar[ld, "f"]
	\\
		&
		C
			&
	\end{tikzcd}
	\]
	We say a double category $\dbl{D}$ \emph{has the Beck-Chevalley pullbacks}
	if the vertical category $\oneV{\dbl{D}}$ has pullbacks and
	their pullback squares are all Beck-Chevalley pullback squares.
\end{definition}

\begin{definition}
	We say a vertical arrow $f\colon A\to B$ is a \emph{cover} if the identity diamond cell on the left below satisfies the Beck-Chevalley condition.
	Dually, $f$ is an \emph{inclusion} if the identity diamond cell on the right below satisfies the Beck-Chevalley condition.
	\begin{equation}
		\label{eqn:CoverInclusion}
		\begin{tikzcd}[column sep=small]
				&
				A
				\ar[ld, "f"']
				\ar[rd, "f"]
				\doublecell[dd]{=}
					&
			\\
			B
			\ar[rd, equal]
				&
					&
					B
					\ar[ld, equal]
			\\
				&
				B
					&
		\end{tikzcd}
		\hspace{2ex}
		,
		\hspace{2ex}
		\begin{tikzcd}[column sep=small]
				&
				A
				\doublecell[dd]{=}
					&
			\\
			A
			\ar[ru, equal]
				&
					&
					A
					\ar[lu, equal]
			\\
				&
				B
				\al[lu, "f"]
				\al[ru, "f"']
					&
		\end{tikzcd}
	\end{equation}
	We let $\Cov(\dbl{D})$ and $\Inc(\dbl{D})$ denote the class of covers and inclusions in $\dbl{D}$, respectively.
\end{definition}
\begin{remark}
	In other words, $f\colon A\to B$ is an inclusion if the restriction $B(f,f)$ is isomorphic to the horizontal identity $\Id_A$,
	and $f\colon A\to B$ is a cover if the extension $\Ext(A;f,f)$ is isomorphic to the horizontal identity $\Id_B$. 
	With inclusions and covers, we gain a better command of
	the diagrammatic calculation of cartesian and opcartesian cells via the sandwich lemma \Cref{lem:Sandwich}.
	For example, the following cells are all cartesian,
	where $\hto$ and $\thto$ denote an inclusion and a cover, respectively.
	\[
		\begin{tikzcd}
			\cdot
			\ar[d, hook]
			\sar[r]
			\doublecell{\opcart}
				&
				\cdot
				\ar[d]
			\sar[r]
			\doublecell{\cart}
				&
				\cdot
				\ar[d]
				\\
				\cdot
				\sar[r]
				&
				\cdot
				\sar[r]
				&
				\cdot
		\end{tikzcd}
		\hspace{2ex}
		,
		\hspace{2ex}
		\begin{tikzcd}
			\cdot
			\ar[d]
			\sar[r]
			\doublecell{\cart}
				&
				\cdot
				\ar[d, two heads]
				\sar[r]
				\doublecell{\cart} 
					&
					\cdot
					\ar[d]
					\\
					\cdot
					\sar[r]
					&
					\cdot
					\sar[r]
					&
					\cdot
		\end{tikzcd}
		\hspace{2ex}
		,
		\hspace{2ex}
		\begin{tikzcd}
			&
			\cdot
			\ar[ld]
			\ar[rd]
			\doublecell[d]{\cart}
				&
			\\
			\cdot
			\sar[r]
			\ar[d, hook]
			\doublecell{\opcart}
				&
			\cdot
			\sar[r]
			\ar[d, hook]
			\doublecell{\opcart}
				&
			\cdot
			\ar[d, hook]
			\\
			\cdot
			\sar[r]
			&
			\cdot
			\sar[r]
			&
			\cdot
		\end{tikzcd}
	\]
\end{remark}
\begin{example}
	In $\Span(\one{C})$, all morphisms are covers and
	inclusions are precisely monomorphisms.

	In $\Prof$ as in 
	\cref{eg:Prof}, an inclusion is precisely a fully faithful functor,
	while a cover is precisely a functor that 
	is called a \textit{Cauchy dense} functor in \cite{Day77}
	or an \textit{absolutely dense} functor in \cite{EV02}.
	This is also the same as \textit{lax epimorphism} in the 2-category $\Catbi$
	in the sense of \cite{AESV01}.
\end{example}

\begin{lemma}\label{lem:monoincl}
	In a double category with Beck-Chevalley pullbacks,
	every monomorphism in its vertical category is an inclusion.
\end{lemma}
\begin{proof}
	The commutative square on the right in \Cref{eqn:CoverInclusion} is a pullback square
	if $f$ is a monomorphism.
\end{proof}

As a special case of Beck-Chevalley pullbacks,
we shall introduce the following notion, which is of interest in the context.

\begin{definition}[{cf.\ \cite[Definition 3.10]{LWW10}, \cite[Definition 2.1]{CW87}, \cite[Definition 3.1]{WW08}}]
	\label{defn:discrete}
	Let $\dbl{D}$ be a cartesian double category, and
	$X\in\dbl{D}$ be an object.
	Then
	we say $X$ is \emph{discrete}
	if the following pullback squares in $\oneV{\dbl{D}}$
	are all Beck-Chevalley pullback squares.
	\begin{equation}\label{dgm:Discrete}
		\begin{tikzcd}[column sep=tiny]
				&
				X
				\ar[ld, equal]
				\ar[rd, equal]
				\pullbackA[dd]
					&
			\\
			X
			\ar[rd, "\Delta"']
				&
					&
					X
					\ar[ld, "\Delta"]
			\\
				&
				X\times X
					&
		\end{tikzcd}
		\hspace{1ex}
		,
		\hspace{1ex}
		\begin{tikzcd}[column sep=tiny]
					&
					X
					\ar[ld, "\Delta"']
					\ar[rd, "\Delta"]
					\pullbackA[dd]
						&
				\\
				X\times X
				\ar[rd, "\id_X\times\Delta"']
					&
						&
						X\times X
						\ar[ld, "\Delta\times\id_X"]
				\\
					&
					X\times X\times X
						&
		\end{tikzcd}
		.
	\end{equation}
	Equivalently, we say $X$ is discrete if $\Delta$ is an inclusion and
	the pullback square formed by
	pulling back $\Delta\times\id_X$ along $\id_X\times\Delta$
	is a Beck-Chevalley pullback square.

	We say a cartesian equipment is \emph{discrete} if
	every object is discrete.
\end{definition}

Since all the three identity cells in \cref{dgm:Discrete} are made of pullbacks squares,
a cartesian equipment with Beck-Chevalley pullbacks is discrete.

\begin{remark}
	In \cite{WW08}, a \textit{Frobenius object} is defined as an object $X$ such that the latter two cells in \cref{dgm:Discrete} 
	satisfy the Beck-Chevalley condition, while in \cite{LWW10}, an object satisfying the Beck-Chevalley for the first cell in \cref{dgm:Discrete} is
	called a \textit{separable} object.
	In other words, a separable object is an object $A$ whose diagonal
	$\Delta\colon A\to A\times A$ is an inclusion. Moreover,
	since the horizontal identity $\Id_{A\times A}$ is isomorphic to
	$\Id_{A}\times\Id_{A}$, $A$ is separable if and only if the projection
	for the local product
	$\Id_{A}\land\Id_{A}\to\Id_{A}$ is invertible.
\end{remark}

Allegories were defined in \cite{FS90,Joh02} to construct an abstract framework of the category of relations.
The definition involves a condition called the \textit{modular law}, 
which states that for any morphisms $Q\colon A\sto B, R\colon B\sto C$ and $S\colon A\sto C$,
the following inequality holds:
\begin{equation}
	\label{eqn:ModularLawAllegory}
	QR\land S \leq Q(R\land Q^\circ S)
\end{equation}
where $Q^\circ$ is the involution of $Q$.
If $Q$ is of the form $f^*$ for some map $f$,
then the opposite inequality follows from the counit of the adjunction $(f^*)^\circ\dashv f^*$.
The equality for this case is sometimes called the \textit{Frobenius condition}
as an instance of Frobenius reciprocity for adjoint functors.

\begin{definition}
	Let $\dbl{D}$ be a cartesian equipment.
	We say $\dbl{D}$ satisfies the \emph{modular law for vertical arrows}
	if, for any vertical arrow $f$ and any two horizontal arrows $R, S$ of the form
	\[
		\begin{tikzcd}
			A
			\sar[r, "R"]
			\ar[d, "f"']
				&
				B
			\\
			X
			\sar[r, "S"]
				&
				B
		\end{tikzcd}\text{,}
	\]
	the cell on the left below
	factors as an opcartesian cell followed by a cartesian cell
	as in the right below.
	\begin{equation}\label{eqn:ModularLaw}
		\begin{tikzcd}
			A
			\sar[r, "R\land f_!S"]
			\ar[d, "\mmbox{\langle\id, f\rangle}"']
			\doublecell{\cart}
				&
				B
				\ar[d, "\Delta"]
			\\
			A\times X
			\sar[r, "R\times S"]
			\ar[d, "f\times\id"']
			\doublecell{\opcart}
				&
				B\times B
				\ar[d, equal]
			\\
			X\times X
			\sar[r, "f^*R\times S"']
				&
				B\times B
		\end{tikzcd}
		\hspace{1ex}
		=
		\hspace{1ex}
		\begin{tikzcd}
			A
			\sar[r, "R\land f_!S"]
			\ar[d, "f"']
			\doublecell{\opcart}
				&
				B
				\ar[d, equal]
			\\
			X
			\sar[r]
			\ar[d, "\Delta"']
			\doublecell{\cart}
				&
				B
				\ar[d, "\Delta"]
			\\
			X\times X
			\sar[r, "f^*R\times S"']
				&
				B\times B
		\end{tikzcd}
	\end{equation}
\end{definition}

Since the local product $\wedge$ in $\biH{\dbl{D}}(A,B)$ is given by the restriction along the diagonals,
the top horizontal arrow of the cartesian cell on the right is
$f^*R\wedge S$.
On the other hand, the bottom horizontal arrow of the opcartesian cell on the right
is $f^*(R\wedge f_!S)$.
The modular law for vertical arrows, therefore, states the existence of an isomorphism $f^*R\wedge S \cong f^*(R\wedge f_!S)$ for any $f$, $R$, and $S$,
which deserves its name.

\begin{proposition}\label{prop:BCModularLaw}
	A cartesian equipment $\dbl{D}$ satisfies the modular law for vertical arrows
	if $\dbl{D}$ has Beck-Chevalley pullbacks.
\end{proposition}
\begin{proof}
	Since a square
	$f\fatsemi \Delta = \langle\id, f\rangle\fatsemi (f\times\id)$ is a pullback,
	the right-hand composition in \cref{eqn:ModularLaw} is achieved as follows.
	\[
		\begin{tikzcd}
				&
				A
				\sar[rrr]
				\ar[dr, "\mmbox{\langle\id, f\rangle}"description]
				\ar[dl, "f"']
					&
					\!
					\doublecell{\cart}
						&
							&
							B
							\ar[dl, "\Delta"description]
							\ar[dr, equal]
								&
			\\
			X
			\ar[dr, "\Delta"']
			\sar[rr, ""]
			\ar[rr,
					shift left=3ex,
					phantom, "\opcart"{description, inner sep=0mm},
					start anchor={[xshift=0ex, yshift=0ex]center},
					end anchor={[xshift=-2ex, yshift=0ex]center},
					]
			\ar[rr,
					shift right=3ex,
					phantom, "\cart"{description, inner sep=0mm},
					start anchor={[xshift=0ex, yshift=0ex]center},
					end anchor={[xshift=-2ex, yshift=0ex]center},
					]
				&
					&
					A\times X
					\sar[r, "R\times S"]
					\ar[dl, "f\times\id"description]
					\doublecell{\opcart}
						&
						B\times B
						\ar[rd, equal]
						\sar[rr, "\Delta^*"near end]
					\ar[rr,
							shift left=3ex,
							phantom, "\opcart"{description, inner sep=0mm},
							start anchor={[xshift=2ex, yshift=0ex]center},
							end anchor={[xshift=0ex, yshift=0ex]center},
							]
					\ar[rr,
							shift right=3ex,
							phantom, "\cart"{description, inner sep=0mm},
							start anchor={[xshift=2ex, yshift=0ex]center},
							end anchor={[xshift=0ex, yshift=0ex]center},
							]
							&
							\!
								&
								B
								\ar[dl, "\Delta"]
			\\
				&
				X\times X
				\sar[rrr]
					&
						&
						\!
							&
							B\times B
								&
		\end{tikzcd}
	\]
	This diagram gives the desired factorisation
	by \cref{lem:Sandwich}.
\end{proof}

\begin{definition}
	Let $\dbl{D}$ be a cartesian equipment, and $X$ be an object of $\dbl{D}$.
	We say $X$ is \emph{self-dual} if there exist horizontal arrows $\eta_X\colon 1\sto X\times X$ and
	$\varepsilon_X\colon X\times X \sto 1$ equipped with horizontal isomorphisms
	$(\eta_X\times\Id_X)(\Id_X\times\varepsilon_X)\cong\Id_X$ and
	$(\Id_X\times\eta_X)(\varepsilon_X\times\Id_X)\cong\Id_X$.

	Given a vertical arrow $f \colon  X \to Y$, we say $f$ is \emph{self-dual}
	if both $X$ and $Y$ are self-dual and two cells
	\[
		\begin{tikzcd}
			1
			\sar[r, "\eta_X"]
			\ar[d, equal]
			\doublecell[rd]{\eta_f}
				&
				X\times X
				\ar[d, "f\times f"]
			\\
			1
			\sar[r, "\eta_Y"']
				&
				Y\times Y
		\end{tikzcd}
		\hspace{2ex}
		,
		\hspace{2ex}
		\begin{tikzcd}
			X\times X
			\sar[r, "\varepsilon_X"]
			\ar[d, "f\times f"']
			\doublecell[rd]{\varepsilon_f}
				&
				1
				\ar[d, equal]
			\\
			Y\times Y
			\sar[r, "\varepsilon_Y"']
				&
				1
		\end{tikzcd}
	\]
	exist and satisfy the following equalities.
	\[
		\begin{tikzcd}
			\hspace{2ex}X\hspace{2ex}
			\sar[r, "\Id_X\times\eta_X"]
			\ar[d, "f"']
			\doublecell[rd]{\Id_f\times\eta_f}
				&
				X\times X\times X
				\sar[r, "\varepsilon_X\times\Id_X"]
				\ar[d, "f\times f\times f"description]
				\doublecell[rd]{\varepsilon_f\times\Id_f}
					&
					\hspace{2ex}X\hspace{2ex}
					\ar[d, "f"]
			\\
			\hspace{2ex}Y\hspace{2ex}
			\sar[r, "\Id_Y\times\eta_Y"']
			\doublecell[rrd]{\vertcong}
			\ar[d, equal]
				&
				Y\times Y\times Y
				\sar[r, "\varepsilon_Y\times\Id_Y"']
					&
					\hspace{2ex}Y\hspace{2ex}
					\ar[d, equal]
			\\
			\hspace{2ex}Y\hspace{2ex}
			\sar[rr, "\Id_Y"']
				&
					&
					\hspace{2ex}Y\hspace{2ex}
		\end{tikzcd}
		=
		\begin{tikzcd}
			\hspace{2ex}X\hspace{2ex}
			\sar[r, "\Id_X\times\eta_X"]
			\ar[d, equal]
			\doublecell[rrd]{\vertcong}
				&
				X\times X\times X
				\sar[r, "\varepsilon_X\times\Id_X"]
					&
					\hspace{2ex}X\hspace{2ex}
					\ar[d, equal]
			\\
			\hspace{2ex}X\hspace{2ex}
			\sar[rr, "\Id_X"]
			\ar[d, "f"']
			\doublecell[rrd]{\Id_f}
				&
					&
					\hspace{2ex}X\hspace{2ex}
					\ar[d, "f"]
			\\
			\hspace{2ex}Y\hspace{2ex}
			\sar[rr, "\Id_Y"']
				&
					&
					\hspace{2ex}Y\hspace{2ex}
		\end{tikzcd}
	\]
	\[
		\begin{tikzcd}
			\hspace{2ex}X\hspace{2ex}
			\sar[r, "\eta_X\times\Id_X"]
			\ar[d, "f"']
			\doublecell[rd]{\eta_f\times\Id_f}
				&
				X\times X\times X
				\sar[r, "\Id_X\times\varepsilon_X"]
				\ar[d, "f\times f\times f"description]
				\doublecell[rd]{\Id_f\times\varepsilon_f}
					&
					\hspace{2ex}X\hspace{2ex}
					\ar[d, "f"]
			\\
			\hspace{2ex}Y\hspace{2ex}
			\sar[r, "\eta_Y\times\Id_Y"']
			\doublecell[rrd]{\vertcong}
			\ar[d, equal]
				&
				Y\times Y\times Y
				\sar[r, "\Id_Y\times\varepsilon_Y"']
					&
					\hspace{2ex}Y\hspace{2ex}
					\ar[d, equal]
			\\
			\hspace{2ex}Y\hspace{2ex}
			\sar[rr, "\Id_Y"']
				&
					&
					\hspace{2ex}Y\hspace{2ex}
		\end{tikzcd}
		=
		\begin{tikzcd}
			\hspace{2ex}X\hspace{2ex}
			\sar[r, "\eta_X\times\Id_X"]
			\ar[d, equal]
			\doublecell[rrd]{\vertcong}
				&
				X\times X\times X
				\sar[r, "\Id_X\times\varepsilon_X"]
					&
					\hspace{2ex}X\hspace{2ex}
					\ar[d, equal]
			\\
			\hspace{2ex}X\hspace{2ex}
			\sar[rr, "\Id_X"]
			\ar[d, "f"']
			\doublecell[rrd]{\Id_f}
				&
					&
					\hspace{2ex}X\hspace{2ex}
					\ar[d, "f"]
			\\
			\hspace{2ex}Y\hspace{2ex}
			\sar[rr, "\Id_Y"']
				&
					&
					\hspace{2ex}Y\hspace{2ex}
		\end{tikzcd}
	\]
	In other words, there exist invertible 2-cells $(\Id_f\times\eta_f)(\varepsilon_f\times\Id_f)\cong\Id_{f}$ and
	$(\eta_f\times\Id_f)(\Id_f\times\varepsilon_f)\cong\Id_{f}$ in $\bi{C}_h(\dbl{D})$
	that are compatible with
	those for $X$ and $Y$.
\end{definition}

\begin{remark}\label{rem:DagcpctCorrespondence}
	For each triple of self-dual objects $X, Y, Z$,
	the mapping 
	$p\mapsto (\Id_Z\times p)(\varepsilon_Z\times\Id_Y)$
	defines an equivalence
	$\biH{\dbl{D}}(X,Z\times Y)\simeq\biH{\dbl{D}}(Z\times X,Y)$.
	The pseudo-inverse of this equivalence is given by
	the mapping $q\mapsto (\eta_Z\times\Id_X)(\Id_Z\times q)$.
	Dually, the mapping $p\mapsto (\Id_X\times\eta_Z)(p\times\Id_Z)$
	defines another equivalence
	$\biH{\dbl{D}}(X\times Z, Y)\simeq\biH{\dbl{D}}(X,Y\times Z)$.

	The composition of these two equivalences in particular cases gives an equivalence
	$\biH{\dbl{D}}(X,Y)\simeq\biH{\dbl{D}}(Y,X)$,
	which we write $p \mapsto p^\dagger=(\Id_Y\times\eta_X)(\Id_Y\times p\times\Id_X)(\varepsilon_Y\times\Id_X)$.

	In the same way, for self-dual vertical arrows $f, g, h$,
	we obtain the following three kinds of equivalences for cells,
	which are compatible with those equivalences for horizontal arrows.
	\[
	\begin{split}
		\alpha\mapsto \alpha^\kappa:=(\Id_h\times\alpha)(\varepsilon_h\times\Id_g)
		\colon \ 
		&
		\bi{C}_h(\dbl{D})(f, h\times g)
		\ \simeq\ 
		\bi{C}_h(\dbl{D})(h\times f,g)
		\\
		\alpha\mapsto \alpha^\lambda:=(\Id_f\times\eta_h)(\alpha\times\Id_h)
		\colon \ 
		&
		\bi{C}_h(\dbl{D})(f\times h, g)
		\ \simeq\ 
		\bi{C}_h(\dbl{D})(f,g\times h)
		\\
		\alpha\mapsto \alpha^\dagger:=(\Id_g\times\eta_f)(\Id_g\times\alpha\times\Id_f)(\varepsilon_g\times\Id_f)
		\colon \ 
		&
		\bi{C}_h(\dbl{D})(f, g)
		\ \simeq\ 
		\bi{C}_h(\dbl{D})(g, f)
	\end{split}
	\]
\end{remark}

\begin{proposition}\label{prop:discselfdual}
	In a discrete cartesian equipment
	$\dbl{D}$, every object and every vertical arrow are self-dual.
\end{proposition}
\begin{proof}
	Let $X\in\dbl{D}$ be an object. Define $\eta_X\colon 1\sto X\times X$ and $\varepsilon_X\colon X\times X\sto 1$ as the following extensions.
	\[
		\begin{tikzcd}[column sep=small]
				&
				X
				\ar[ld, "!"']
				\ar[rd, "\Delta"]
					&
			\\
			\hspace{3ex}1\hspace{3ex}
			\doublecell[rr, shift left=3ex]{\opcart}
			\sar[rr, "\eta_X"']
				&
					&
					X\times X
		\end{tikzcd}
		\hspace{3ex}
		\begin{tikzcd}[column sep=small]
				&
				X
				\ar[ld, "\Delta"']
				\ar[rd, "!"]
					&
			\\
			X\times X
			\doublecell[rr, shift left=3ex]{\opcart}
			\sar[rr, "\varepsilon_X"']
				&
					&
					\hspace{3ex}1\hspace{3ex}
		\end{tikzcd}
	\]
	We now show that these horizontal arrows give the self-dual of $X$. Consider the following cells.
	\begin{equation}\label{dgm:selfdual}
		\begin{tikzcd}[column sep = small]
				&
					&
					X
					\ar[ld,"\Delta"']
					\ar[lldd,"\id"', bend right = 45]
					\ar[rd,"\Delta"]
					\ar[rrdd,"\id", bend left = 45]
						&
							&
			\\
				&
				X\times X
				\ar[ld,"\id\times !"']
				\ar[rd,"\id\times\Delta"description]
				\doublecell[rr, shift left=3ex]{\opcart}
				\doublecell[rr, shift right=3ex]{\cart}
				\sar[rr]
					&
						&
						X\times X
						\ar[ld,"\Delta\times\id"description]
						\ar[rd,"!\times\id"]
								&
			\\
			\hspace{3.5ex}X\hspace{3.5ex}
			\sar[rr, "\Id\times\eta_X"']
			\doublecell[rr, shift left=2.5ex]{\opcart}
				&
					&
					X\times X\times X
					\sar[rr,"\varepsilon_X\times\Id"']
					\doublecell[rr, shift left=2.5ex]{\opcart}
						&
							& 
							\hspace{3.5ex}X\hspace{3.5ex}
		\end{tikzcd}
	\end{equation}
	The left and right cells are opcartesian
	because $X\times -$ and $-\times X$ define pseudo-functors.
	The Beck-Chevalley condition for the pullback square in the middle follows from the discreteness of $X$.
	\cref{lem:Sandwich} implies that this establishes an extension of the identity vertical arrows on $X$,
	thereby inducing the horizontal isomorphism $(\Id_X\times\eta_X)(\varepsilon_X\times\Id_X)\cong\Id_X$.
	The other isomorphism for them to be the data of the self-dual is verified similarly.

	Given a vertical arrow $f\colon X\to Y$, one can obtain $\eta_f$ through the following equality, and $\varepsilon_f$ can be determined likewise.
	\[
		\begin{tikzcd}[column sep=small]
				&
				X
				\ar[ld, "!"']
				\ar[rd, "\Delta"]
					&
			\\
			\hspace{2ex}1\hspace{2ex}
			\sar[rr, "\eta_X"']
			\doublecell[rr, shift left=2.5ex]{\opcart}
			\ar[d, equal]
			\doublecell[rrd]{\eta_f}
				&
						&
						X\times X
						\ar[d, "f\times f"]
			\\
			\hspace{2ex}1\hspace{2ex}
			\sar[rr, "\eta_Y"']
				&
						&
						Y\times Y
		\end{tikzcd}
		=
		\begin{tikzcd}[column sep=small]
				&
				X
				\ar[d, "f"description]
				\ar[ld, "!"']
				\ar[rd, "\Delta"]
					&
			\\
			\hspace{2ex}1\hspace{2ex}
			\ar[d, equal]
				&
				Y
				\ar[ld, "!"']
				\ar[rd, "\Delta"]
						&
						X\times X
						\ar[d, "f\times f"]
			\\
			\hspace{2ex}1\hspace{2ex}
			\sar[rr, "\eta_Y"']
			\doublecell[rr, shift left=2.5ex]{\opcart}
				&
						&
						Y\times Y
		\end{tikzcd}
	\]
	Composing vertically the opcartesian cell \cref{dgm:selfdual}
	with the horizontal composite of $\Id_f\times\eta_f$ and $\varepsilon_f\times\Id_f$,
	one can check the desired equality for the self-duality of $f$ by a diagrammatic calculation.
\end{proof}

In a discrete cartesian equipment, $\eta$ and $\varepsilon$ are defined by the universal properties of opcartesian cells,
ensuring their compatibility with vertical compositions as below.
\begin{lemma}\label{lem:propofDiscreteDagger}
	Let $\dbl{D}$ be a discrete cartesian equipment, and consider the structure $(\eta, \varepsilon)$ of self-duals established in \Cref{prop:discselfdual}.
	We have the following.
	\begin{enumerate}
		\item
			$\eta_{\id_X}=\id_{\eta_X}$ and $\varepsilon_{\id_X}=\id_{\varepsilon_X}$
			are valid for each object $X$.
		\item
			Given a compatible pair of vertical arrows $f\colon A\to B$ and $g\colon B\to C$, the following hold.
			\[
				\begin{tikzcd}
					1
					\sar[r, "\eta_A"]
					\ar[d, equal]
					\doublecell[rdd,xshift=-1ex]{\eta_{f\fatsemi g}}
						&
						A\times A
						\ar[d, "f\times f"]
					\\
					1
					\ar[d, equal]
						&
						B\times B
						\ar[d, "g\times g"]
					\\
					1
					\sar[r, "\eta_C"']
						&
						C\times C
				\end{tikzcd}
				=
				\begin{tikzcd}
					1
					\sar[r, "\eta_A"]
					\ar[d, equal]
					\doublecell[rd]{\eta_{f}}
						&
						A\times A
						\ar[d, "f\times f"]
					\\
					1
					\sar[r, "\eta_B"]
					\ar[d, equal]
					\doublecell[rd]{\eta_{g}}
						&
						B\times B
						\ar[d, "g\times g"]
					\\
					1
					\sar[r, "\eta_C"']
						&
						C\times C
				\end{tikzcd}
				\hspace{2ex}
				,
				\hspace{2ex}
				\begin{tikzcd}
					A\times A
					\sar[r, "\varepsilon_A"]
					\ar[d, "f\times f"']
					\doublecell[rdd,xshift=1ex]{\varepsilon_{f\fatsemi g}}
						&
						1
						\ar[d, equal]
					\\
					B\times B
					\ar[d, "g\times g"']
						&
						1
						\ar[d, equal]
					\\
					C\times C
					\sar[r, "\varepsilon_C"']
						&
						1
				\end{tikzcd}
				=
				\begin{tikzcd}
					A\times A
					\ar[d, "f\times f"']
					\sar[r, "\varepsilon_A"]
					\doublecell[rd]{\varepsilon_{f}}
						&
						1
						\ar[d, equal]
					\\
					B\times B
					\ar[d, "g\times g"']
					\sar[r, "\varepsilon_B"]
					\doublecell[rd]{\varepsilon_{g}}
						&
						1
						\ar[d, equal]
					\\
					C\times C
					\sar[r, "\varepsilon_C"']
						&
						1
				\end{tikzcd}
			\]
	\end{enumerate}
\end{lemma}

\begin{proposition}
	\label{prop:pseudoFunctorialityDagger}
	In a discrete cartesian equipment, the equivalence $(-)^\dagger$ described in \cref{rem:DagcpctCorrespondence} extends to
	an equivalence  $\dagger\colon \dbl{D}\to\dbl{D}^\hop$ in $\DblCat$.
\end{proposition}
\begin{proof}
	We only show $(-)^\dagger$ extends to a pseudo-functor $\dagger\colon \dbl{D}\to\dbl{D}^\hop$.
	This is because one can similarly verify that the assignment
	$p\mapsto (\eta\times\Id)(\Id\times p\times\Id)(\Id\times\varepsilon)$
	also extends to a pseudo-functor,
	and it gives the pseudo-inverse of $\dagger$.

	$\dagger$ is the identity on the vertical category $\oneV{\dbl{D}}$.
	Firstly, considering $i)$ and $ii)$ of \cref{lem:propofDiscreteDagger}, we obtain the vertical functoriality of $\dagger$.
	It remains to show the horizontal functoriality.

	For each $f$, the cell $\Id_f^\dagger$ is isomorphic to $\Id_f$ through the canonical isomorphism for the self-duality of $f$.
	Let $\alpha\colon f\to g$ and $\beta\colon g\to h$ be 1-cells in $\bi{C}_h(\dbl{D})$.
	Consider the following invertible 2-cells in this bicategory.
	\small
	\begin{equation}
	\begin{split}
		&
		\begin{tikzcd}[ampersand replacement=\&, column sep=8ex]
			h
			\ar[r, "\beta^\dagger\odot\alpha^\dagger"]
			\&
			f
		\end{tikzcd}
		\\
		\isoRightarrow
		&
		\begin{tikzcd}[ampersand replacement=\&, column sep=8ex]
			h
			\ar[r,
				"
				\Id_h
				\times
				\eta_g
				"
			]
			\&
			h
			g
			g
			\ar[rr,
				"
				\Id_h
				\times
				\beta
				\times
				\Id_g
				"
			]
			\&
			\&
			h
			h
			g
			\ar[r,
				"
				\varepsilon_h
				\times
				\Id_g
				"
			]
			\&
			g
			\ar[r,
				"
				\Id_g
				\times
				\eta_f
				"
			]
			\&
			g
			f
			f
			\ar[rr,
				"
				\Id_g
				\times
				\alpha
				\times
				\Id_f
				"
			]
			\&
			\&
			g
			g
			f
			\ar[r,
				"
				\varepsilon_g
				\times
				\Id_f
				"
			]
			\&
			f
		\end{tikzcd}
		\\
		\isoRightarrow
		&
		\begin{tikzcd}[ampersand replacement=\&, column sep=8ex]
			h
			\ar[r,
				"\Id_h
				\times
				\eta_f
				"
			]
			\&
			h
			f
			f
			\ar[r,
				"
				\Id_h
				\times
				\alpha
				\times
				\Id_f
				"
			]
			\&
			h
			g
			f
			\ar[rr,
				"
				\Id_h
				\times
				\eta_g
				\times
				\Id_g
				\times
				\Id_f
				"
			]
			\&
			\&
			h
			g
			g
			g
			f
			\ar[rr,
				"
				\Id_h
				\times
				\Id_g
				\times
				\varepsilon_g
				\times
				\Id_f
				"
			]
			\&
			\&
			h
			g
			f
			\ar[r,
				"
				\Id_h
				\times
				\beta
				\times
				\Id_f
				"
			]
			\&
			h
			h
			f
			\ar[r,
				"
				\varepsilon_h
				\times
				\Id_f
				"
			]
			\&
			f
		\end{tikzcd}
		\\
		\isoRightarrow
		&
		\begin{tikzcd}[ampersand replacement=\&, column sep=8ex]
			h
			\ar[r,
				"
				\Id_h
				\times
				\eta_f
				"
			]
			\&
			h
			f
			f
			\ar[rr,
				"
				\Id_h
				\times
				(\alpha\odot\beta)
				\times
				\Id_f
				"
			]
			\&
			\&
			h
			h
			f
			\ar[r,
				"
				\varepsilon_h
				\times
				\Id_f
				"
			]
			\&
			f
		\end{tikzcd}
		\hspace{2ex}
		=
		\hspace{2ex}
		\begin{tikzcd}[ampersand replacement=\&,column sep=8ex]
			h
			\ar[r, "{(\alpha\odot\beta)^\dagger}"]
			\&
			f
		\end{tikzcd}
	\end{split}
	\end{equation}
	\normalsize
	The first 2-cell is given by definition,
	and the third one is through the canonical isomorphism for the self-duality of $g$.
	The second arises from the pseudo-functoriality of $-\times-\times-\times-\times-\colon \bi{C}_h(\dbl{D})^5\to\bi{C}_h(\dbl{D})$
	by considering cells obtained from the coherence of the pseudo-monoid structure of $\dbl{D}$ in $\DblCat$.
	We omit the verification of the coherence condition required to show that this isomorphism indeed gives the structure of a pseudo-functor. 
\end{proof}

In particular, $\dagger$ induces an equivalence of categories $\biV{\dbl{D}}(A,B)\simeq\biV{\dbl{D}}(A,B)^\op$ for each pair of objects $A$ and $B$. 
Moreover, this 2-category $\biV{\dbl{D}}$ is essentially a 1-category: 
\begin{lemma}[{cf.\ \cite[Proposition 3.13]{LWW10}}]
	Let $\dbl{D}$ be a discrete cartesian equipment. Then the vertical 2-category $\biV{\dbl{D}}$ is \textit{locally essentially discrete}; i.e.,
	for each pair of objects $A, B$, the hom-category $\biV{\dbl{D}}(A, B)$ is equivalent to a discrete category.
\end{lemma}
\begin{proof}
	It suffices to show $\biV{\dbl{D}}$ is locally preordered.
	For each pair of parallel 2-cells $\alpha,\beta\colon g\Rightarrow f\colon A\to B$ in $\biV{\dbl{D}}$,
	consider the following vertical cell shown on the left below. 
	The vertical arrow $\Delta$ being an inclusion assures that this uniquely factors through
	$\Delta$ as follows.
	\[
		\begin{tikzcd}
				&
				A
				\ar[d, "\Delta"]
				\ar[ld, bend right, "g"']
				\ar[rd, bend left, "f"]
					&
			\\
			B
				&
				A\times A
				\ar[d, "g\times g"', bend right=55]
				\ar[d, "f\times f", bend left=55]
				\doublecell[r]{=}
				\doublecell[l]{=}
				\doublecell[d]{\alpha\times\beta}
					&
					B
			\\
				&
				B\times B
				\al[lu, bend left, "\Delta"]
				\al[ru, bend right, "\Delta"']
					&
		\end{tikzcd}
		\hspace{2ex}
		=
		\hspace{2ex}
		\begin{tikzcd}
			A
			\ar[d, bend right, "g"']
			\doublecell[d]{\theta}
			\ar[d, bend left, "f"]
			\\
			B
			\ar[d, "\Delta"]
			\\
			B\times B
		\end{tikzcd}
	\]
	Here, the vertical cell $\alpha\times\beta$ is the image of the pair $(\alpha,\beta)$ under the 2-functor
	$\biV{-\times-}\colon\biV{\dbl{D}}\times\biV{\dbl{D}}\to\biV{\dbl{D}}$.
	By postcomposing the projections, we obtain $\alpha=\theta=\beta$. This follows from the fact that
	$\bi{V}$ preserves cartesian objects, as we observed in
	\cref{rem:biVPreservesCartObj}.
\end{proof}

Due to the preservation of companions by any pseudo-functor, the following corollary is deduced.
\begin{corollary}\label{prop:ConjointDagger}
	Let $\dbl{D}$ be a
	discrete cartesian equipment,
	and $f\colon A \to X$ be a vertical arrow.
	Then, the conjoint $f^*$ is isomorphic to $(f_!)^\dagger$.
\end{corollary}

\begin{remark}
	\label{rem:involution}
	By \Cref{prop:ConjointDagger}, we have $(f^*g_!)^\dagger\cong g^*f_!$ for any pair of vertical arrows $f\colon A\to X$ and $g\colon A\to Y$
	in a discrete cartesian equipment $\dbl{D}$.
	If $\dbl{D}$ furthermore has strong tabulators, we can deduce the original form of the \textit{modular law} \Cref{eqn:ModularLawAllegory} from the one for vertical arrows,
	although the inequality in \Cref{eqn:ModularLawAllegory} is replaced by mere existence of a cell.
	This is because every horizontal arrow $Q\colon A\sto B$ is presented as $Q\cong f_!g^*$ for the legs of a strong tabulator of $Q$ and 
	we have
	\begin{align*}
		QR\land S &\cong f_!g^*R\land S \\
		&\cong f_!(g^*R\land f^*S) & \text{(by \Cref{prop:BCModularLaw})}\\
		&\Rightarrow f_!(g^*R\land g^*g_!f^*S) & \text{(by the unit $\eta$ for the representable adjoint)}\\
		&\cong f_!g^*(R\land g_!f^*S) & \text{(since restriction preserves local products)}\\
		&\cong Q(R\land Q^\dagger S).
	\end{align*}
	From this observation, 
	we see that the horizontal bicategories
	of the locally posetal discrete equipments with strong tabulators
	as we will define in 
	\Cref{defn:locallyPreordered} are allegories.
\end{remark}

Using this pseudo-functor, we can deal with cells and horizontal arrows symmetrically. On the other hand,
the other equivalences $(-)^\lambda$ and $(-)^\kappa$ described in \cref{rem:DagcpctCorrespondence}
enable us to `tilt' cells and horizontal arrows to the left and right:
\begin{lemma}\label{lem:GeneralUnbiasing}
	Let $\dbl{D}$ be a
	discrete cartesian equipment and fix two vertical arrows $f\colon A\to A'$ and $g\colon B\to B'$.
	Let $(p\colon A\sto B, q\colon A\times B \sto 1,r\colon B\sto A)$ and $(p'\colon A'\sto B', q'\colon A'\times B' \sto 1,r'\colon B'\sto A')$ be
	triples of corresponding horizontal arrows under the equivalence in \cref{rem:DagcpctCorrespondence}. 
	Then, there are bijective correspondences among cells of the following forms.
	\[
		\begin{tikzcd}
			A
			\sar[r, "p"]
			\ar[d, "f"']
			\doublecell[rd]{\alpha}
				&
				B
				\ar[d, "g"]
			\\
			A'
			\sar[r, "p'"']
				&
				B'
		\end{tikzcd}
		\vline
		\,
		\vline
		\begin{tikzcd}
			A\times B
			\sar[r, "q"]
			\ar[d, "f\times g"']
			\doublecell[rd]{\beta}
				&
				1
				\ar[d, equal]
			\\
			A'\times B'
			\sar[r, "q'"']
				&
				1
		\end{tikzcd}
		\vline
		\,
		\vline
		\begin{tikzcd}
			B
			\sar[r, "r"]
			\ar[d, "g"']
			\doublecell[rd]{\gamma}
				&
				A
				\ar[d, "f"]
			\\
			B'
			\sar[r, "r'"']
				&
				A'
		\end{tikzcd}
	\]
	Moreover, they are vertically functorial.
	This means that for additional vertical arrows $f'\colon A'\to A''$ and $g'\colon B'\to B''$, and a triple 
	$(p''\colon A''\sto B'', q''\colon A''\times B'' \sto 1,r''\colon B''\sto A'')$,
	the correspondences are compatible with the vertical composition of cells.
	In particular, $\alpha$ is (op)cartesian if and only if
	$\gamma$ is as well.
\end{lemma}

\begin{proof}
	We see the correspondence between $\alpha$ and $\beta$, and the other correspondences is similarly verified.
	We have the adjoint equivalence $(-)^\kappa\colon\bi{C}_h(\dbl{D})(f, g)\adjointleft\bi{C}_h(\dbl{D})(f\times g,\id_1)\lon (-)_\kappa$,
	where $(-)_\kappa$ is the pseudo-inverse of $(-)^\kappa$ defined in \cref{rem:DagcpctCorrespondence}.
	Let $\varphi\colon \id\cong (-^\kappa)_\kappa$ and $\psi\colon (-_\kappa)^\kappa\cong \id$ 
	be the invertible unit and counit.

	For pairs $(p,q)$ and $(p',q')$ of corresponding horizontal arrows,
	fix invertible horizontal cells $\zeta\colon q\cong p^\kappa$ and $\zeta'\colon q' \cong p'^\kappa$.
	Consider the following assignments.
	\[
		\begin{tikzcd}
			A
			\sar[r, "p"]
			\ar[d, "f"']
			\doublecell[rd]{\alpha}
				&
				B
				\ar[d, "g"]
			\\
			A'
			\sar[r, "p'"']
				&
			B'
		\end{tikzcd}
		\mapsto
		\begin{tikzcd}[column sep=large]
			A\times B
			\sar[r, "q"]
			\ar[d, equal]
			\doublecell[rd]{\zeta}
				&
				1
				\ar[d, equal]
			\\
			A\times B
			\sar[r,"p^\kappa"',near end]
			\ar[d, "f\times g"']
			\doublecell[rd]{\alpha^\kappa}
				&
				1
				\ar[d, equal]
			\\
			A'\times B'
			\sar[r, "{p'}^\kappa", near end]
			\ar[d, equal]
			\doublecell[rd]{{\zeta'}\inv}
				&
				1
				\ar[d, equal]
			\\
			A'\times B'
			\sar[r, "q'"']
				&
				1
		\end{tikzcd}
		\hspace{2ex}
		,
		\hspace{2ex}
		\begin{tikzcd}
			A\times B
			\sar[r, "q"]
			\ar[d, "f\times g"']
			\doublecell[rd]{\beta}
				&
				1
				\ar[d, equal]
			\\
			A'\times B'
			\sar[r, "q'"']
				&
				1
		\end{tikzcd}
		\mapsto
		\begin{tikzcd}[column sep=large]
			A
			\sar[r, "p"]
			\ar[d, equal]
			\doublecell[rd]{\varphi_p\fatsemi {\zeta_\kappa}\inv}
				&
			B
				\ar[d, equal]
			\\
			A
			\sar[r,"q_\kappa"',near end]
			\ar[d, "g"']
			\doublecell[rd]{\beta_\kappa}
				&
				B
				\ar[d, "f"]
			\\
			A'
			\sar[r, "q'_\kappa", near end]
			\ar[d, equal]
			\doublecell[rd]{\zeta'_\kappa\fatsemi \varphi_{p'}\inv}
				&
				B'
				\ar[d, equal]
			\\
			A'
			\sar[r, "p'"']
				&
			B'
		\end{tikzcd}
	\]
	Both assignments are vertically functorial by the 
	vertical functoriality of $(-)^\kappa$ and $(-)_\kappa$ in \cref{lem:propofDiscreteDagger}.
	To show that they are mutually inverse, consider the following equations in
	the category $\dbl{D}_1$ consisting of horizontal arrows and cells in $\dbl{D}$.
	Starting from $\alpha$, we have
	\[
	\begin{split}
		&
		\begin{tikzcd}[ampersand replacement=\&, column sep=5ex]
			p
			\ar[r,"\varphi_p"]
			\&
			(p^\kappa)_\kappa
			\ar[r,"\zeta_\kappa\inv"]
			\&
			q_\kappa
			\ar[rrr,"{\left(\zeta\fatsemi\alpha^\kappa\fatsemi(\zeta')\inv\right)_\kappa}"]
			\&
			\&
			\&
			{q'}_\kappa
			\ar[r,"(\zeta'_\kappa)\inv"]
			\&
			({p'}^\kappa)_\kappa
			\ar[r,"(\varphi_{p'})\inv"]
			\&
			p'
		\end{tikzcd}
		\\
		=&
		\begin{tikzcd}[ampersand replacement=\&, column sep=5ex]
			p
			\ar[r,"\varphi_p"]
			\&
			(p^\kappa)_\kappa
			\ar[r,"(\alpha^\kappa)_\kappa"]
			\&
			({p'}^\kappa)_\kappa
			\ar[r,"(\varphi_{p'})\inv"]
			\&
			p'
		\end{tikzcd}
		\qquad =\,
		\begin{tikzcd}[ampersand replacement=\&, column sep=5ex]
			p
			\ar[r,"\alpha"]
			\&
			p'
		\end{tikzcd}.
	\end{split}
	\]
	The last equality follows from the naturality of $\varphi$.

	On the other hand, starting from $\beta$,
	we have
	\[
	\begin{split}
		&
		\begin{tikzcd}[ampersand replacement=\&, column sep=5ex]
			q
			\ar[r,"\zeta"]
			\&
			p^\kappa
			\ar[rrrr,"{\left(\varphi_p\fatsemi\zeta_\kappa\inv\fatsemi\beta_\kappa\fatsemi\zeta'_\kappa\fatsemi{\varphi_{p'}}\inv\right)^\kappa}"]
			\&
			\&
			\&
			\&
			{p'}^\kappa
			\ar[r,"{\zeta'}\inv"]
			\&
			q'
		\end{tikzcd}
		\\
		=&
		\begin{tikzcd}[ampersand replacement=\&, column sep=5ex]
			q
			\ar[r,"\zeta"]
			\&
			p^\kappa
			\ar[r,"{(\psi_{p^\kappa})\inv}"]
			\&
			{\left((p^\kappa)_\kappa\right)^\kappa}
			\ar[rrr,"{\left(\left(\zeta\inv\fatsemi\beta\fatsemi\zeta'\right)_\kappa\right)^\kappa}"]
			\&\&\&
			\left(({p'}^\kappa)_\kappa\right)^\kappa
			\ar[r,"{(\psi_{{p'}^\kappa})^\kappa}"]
			\&
			{p'}^\kappa
			\ar[r,"{\zeta'}\inv"]
			\&
			q'
		\end{tikzcd}
		\\
		=&
		\begin{tikzcd}[ampersand replacement=\&, column sep=5ex]
			q
			\ar[r,"\zeta"]
			\&
			p^\kappa
			\ar[r,"\zeta\inv"]
			\&
			q
			\ar[r,"\beta"]
			\&
			q'
			\ar[r,"\zeta'"]
			\&
			{p'}^\kappa
			\ar[r,"{\zeta'}\inv"]
			\&
			q'
		\end{tikzcd}
		\qquad
		=
		\,
		\begin{tikzcd}[ampersand replacement=\&, column sep=5ex]
			q
			\ar[r,"\beta"]
			\&
			q'
		\end{tikzcd}.
	\end{split}
	\]	
	In the first equality, we turn to 
	the triangle identity for $\varphi$ and $\psi$,
	yielding $(\varphi_p)^\kappa=(\psi_{p^\kappa})\inv$.
	The second equality follows from the naturality of $\psi$. 
\end{proof}
Composing the opcartesian cell defining $\varepsilon$ in \cref{prop:discselfdual}, we obtain: 
\begin{lemma}\label{lem:UnbiasingLemma}
	Let $\dbl{D}$ be a discrete cartesian equipment,
	and $f, g\colon  A\to X$ be parallel vertical arrows.
	Suppose that $p\colon X \sto Y$ and $\bar{p}\colon X\times Y \sto 1$ correspond to
	each other as described in \cref{rem:DagcpctCorrespondence}.
	Then,
	cells of the form on the left bijectively correspond to
	cells of the form on the right
	\[
		\begin{tikzcd}[column sep=small]
				&
				A
				\ar[ld, "f"']
				\ar[rd, "g"]
					&
			\\
			\hspace{3ex}X\hspace{3ex}
			\doublecell[rr, shift left=3ex]{\alpha}
			\sar[rr, "p"']
				&
					&
					\hspace{3ex}Y\hspace{3ex}
		\end{tikzcd}
		\hspace{3ex}
		\vline\ \vline
		\hspace{3ex}
		\begin{tikzcd}[column sep=small]
				&
				A
				\ar[ld, "\mmbox{\langle f,g\rangle}"']
				\ar[rd, "!"]
					&
			\\
			X\times Y
			\doublecell[rr, shift left=3ex]{\bar\alpha}
			\sar[rr, "\bar{p}"']
				&
					&
					\hspace{3ex}1\hspace{3ex}
		\end{tikzcd},
	\]
	through the following function.
	\[
		\begin{tikzcd}[column sep=small]
				&
				A
				\ar[ld, "f"']
				\ar[rd, "g"]
					&
			\\
			\hspace{3ex}X\hspace{3ex}
			\doublecell[rr, shift left=3ex]{\alpha}
			\sar[rr, "p"']
				&
					&
					\hspace{3ex}Y\hspace{3ex}
		\end{tikzcd}
		\hspace{3ex}
		\mapsto
		\hspace{3ex}
		\begin{tikzcd}[column sep = small]
				&
					&
					A
					\ar[ld,"\Delta"']
					\ar[rd,"!"]
						&
							&
			\\
				&
				A\times A
				\ar[ld,"f\times g"']
				\ar[rd,"g\times g"description]
				\doublecell[rrrd]{\varepsilon_g}
				\doublecell[rr, shift left=3ex]{\opcart}
				\sar[rr, "\varepsilon_X"']
					&
						&
						\hspace{2.0ex}1\hspace{2.0ex}
						\ar[rd,equal]
								&
			\\
			X\times Y
			\sar[rr, "p\times\Id"']
			\doublecell[rr, shift left=2.5ex]{\alpha\times \Id_g}
				&
					&
					Y\times Y
					\sar[rr,"\varepsilon_Y"']
						&
							& 
							\hspace{2.0ex}1\hspace{2.0ex}
		\end{tikzcd}
	\]
	In particular, $\alpha$ is tabulating / opcartesian
	if and only if $\bar\alpha$ is as well.
\end{lemma}
\begin{proof}
	The bijective correspondence follows from \cref{lem:GeneralUnbiasing}
	and the universal property of the opcartesian cell,
	and by tracking the bijective correspondence, one can also check that
	$\alpha$ is tabulating if and only if
	$\bar\alpha$ is as well.

	It remains to show that 
	$\alpha$ is opcartesian if and only if
	$\bar\alpha$ is as well.
	Suppose that $\bar\alpha$ is opcartesian.
	We can observe the following bijective correspondences.
	Specifically, the leftmost correspondence arises from the above function, while the rightmost one stems from the one we have observed in \cref{lem:GeneralUnbiasing}.
	The one in the middle results from the opcartesian cell $\bar{\alpha}$.
	Therefore, $\zeta_1$ bijectively corresponds to $\zeta_2$,
	which arises from the precomposition of $\alpha$
	because of the vertical functoriality of the correspondences in \cref{lem:GeneralUnbiasing}.
	It follows that $\alpha$ is opcartesian.
	\[
		\begin{tikzcd}[column sep=small]
				&
				A
				\ar[ld, "f"']
				\ar[rd, "g"]
					&
			\\
			X
			\ar[d, "h"']
				&
					&
					Y
					\ar[d, "k"]
			\\
			C
			\doublecell[rr, shift left=6ex]{\zeta_1}
			\sar[rr, "q"']
				&
					&
					D
		\end{tikzcd}
		\hspace{2ex}
		\vline
		\,
		\vline
		\hspace{2ex}
		\begin{tikzcd}[column sep=small]
				&
				A
				\ar[ld, "{\langle f, g\rangle}"']
				\ar[rd,"!"]
					&
			\\
			X\times Y
			\ar[d, "h\times k"']
				&
					&
					1
					\ar[d, equal]
			\\
			C\times D
			\doublecell[rr, shift left=6ex]{\bar\zeta_1}
			\sar[rr, "\bar{q}"']
				&
					&
					1
		\end{tikzcd}
		\hspace{2ex}
		\vline
		\,
		\vline
		\hspace{2ex}
		\begin{tikzcd}[column sep=small]
			X\times Y
			\sar[rr, "\bar{p}"]
			\ar[d, "h\times k"']
				&
					&
					1
					\ar[d, equal]
			\\
			C\times D
			\doublecell[rr, shift left=3.5ex]{\bar\zeta_2}
			\sar[rr, "\bar{q}"']
				&
					&
					1
		\end{tikzcd}
		\hspace{2ex}
		\vline
		\,
		\vline
		\hspace{2ex}
		\begin{tikzcd}[column sep=small]
			X
			\sar[rr, "p"]
			\ar[d, "h"']
				&
					&
					Y
					\ar[d, "k"]
			\\
			C
			\doublecell[rr, shift left=3.5ex]{\zeta_2}
			\sar[rr, "q"']
				&
					&
					D
		\end{tikzcd}
	\]

	To establish the converse, suppose that $\alpha$ is opcartesian.
	Consider $\beta$ and $q$ as the following extensions, and assume they are mapped from 
	$\tilde{\beta}$ and $\tilde{q}$
	by the previously mentioned function.
	\[
		\begin{tikzcd}[column sep=small]
				&
				A
				\ar[ld, "f"']
				\ar[rd, "g"]
					&
			\\
			\hspace{3ex}X\hspace{3ex}
			\doublecell[rr, shift left=3ex]{\tilde{\beta}}
			\sar[rr, "\tilde{q}"']
				&
					&
					\hspace{3ex}Y\hspace{3ex}
		\end{tikzcd}
		\hspace{3ex}
		\vline\ \vline
		\hspace{3ex}
		\begin{tikzcd}[column sep=small]
				&
				A
				\ar[ld, "\mmbox{\langle f,g\rangle}"']
				\ar[rd, "!"]
					&
			\\
			X\times Y
			\doublecell[rr, shift left=3ex]{\beta}
			\sar[rr, "q"']
				&
					&
					\hspace{3ex}1\hspace{3ex}
		\end{tikzcd}
	\]
	Since $\beta$ is opcartesian, $\tilde{\beta}$ is as well by the argument above. 
	Therefore, $\alpha$ and $\tilde{\beta}$ are isomorphic through
	the unique cell $\gamma\colon \tilde{q}\Rightarrow p$ obtained by the universal property of the opcartesian cell $\tilde{\beta}$.
	In light of the vertical functoriality of the correspondence discussed in \cref{lem:GeneralUnbiasing},
	$\gamma$ corresponds to a horizontal invertible cell $\bar\gamma\colon q\Rightarrow\bar{p}$,
	which gives rise to a factorisation of $\bar\alpha$ through $\beta$.
	This shows $\bar\alpha$ is also opcartesian.
\end{proof}

\subsection{\texorpdfstring{$\zero{M}$}{M}-comprehension schemes}\label{sec:comprehensionscheme}

We shall initiate our discourse by drawing an archetype of the double category of relations.
As indicated in the introduction, we will define relations relative to a stable orthogonal factorisation system $(\zero{E}, \zero{M})$ in a category $\one{C}$.

\begin{definition}\label{def:RelEM}
	Let $\one{C}$ be a category with finite limits and 
	$(\zero{E}, \zero{M})$ be a stable orthogonal factorisation system in $\one{C}$.
	The double category $\Rel{\zero{E}, \zero{M}}(\one{C})$ is defined as follows.
	\begin{itemize}
		\item The vertical category $\Rel{\zero{E}, \zero{M}}(\one{C})_0$ is precisely the same as $\one{C}$.
		Therefore, objects and vertical arrows in $\Rel{\zero{E}, \zero{M}}(\one{C})$ are the same as objects and arrows in $\one{C}$. 
		\item $\langle\src,\tgt\rangle\colon\Rel{\zero{E},\zero{M}}(\one{C})_1\to\one{C}\times\one{C}$ is defined by the following pullback.
		\[
			\begin{tikzcd}
				\Rel{\zero{E},\zero{M}}(\one{C})_1
				\ar[d, "{\langle\src,\tgt\rangle}"']
				\ar[r]
				\pullback[rd]
					&
					\one{M}
					\ar[d]
				\\
				\one{C}\times\one{C}
				\ar[r, "\times"']
					&
					\one{C}
			\end{tikzcd}
		\]
		Here the functor $\one{M}\to\one{C}$ on the right 
		is is the reflective subfibration of the codomain fibration $\one{C}^{\rightarrow}\to\one{C}$
	defined in \cref{prop:OFSfromStableSystem}.
	The concrete description is given as follows.
	Horizontal arrows are $\zero{M}$-relations in $\one{C}$,
	where a $\zero{M}$-relation $R\colon A\sto B$ is a morphism $\langle l,r\rangle\colon R\to A\times B$ 
	in $\zero{M}$.
	A cell of the form on the left below 
	is an arrow $\alpha\colon R\to S$ that makes the diagram on the right below commute.
	\begin{equation}
		\label{eq:cellofrel}
		\begin{tikzcd}
			A
			\sar[r, "R"]
			\ar[d, "f"']
				&
				B
			\ar[d, "g"]
			\\
			C
			\sar[r, "S"']
			\doublecell[r,shift left=4ex]{\alpha}
				&
				D
		\end{tikzcd}
		\qquad
		\begin{tikzcd}
			R
			\ar[r, tail,"{\langle l, r\rangle}"]
			\ar[d, "\alpha"']
				&
				A\times B
			\ar[d, "f\times g"]
			\\
			S
			\ar[r, tail ,"{\langle l', r'\rangle}"']
				&
				C\times D
		\end{tikzcd}
	\end{equation}
	\item The horizontal composition of $R\colon A\sto B$ and $S\colon B\sto C$ is given by
	the $\zero{M}$-image of the arrow $P\to A\times C$,
	where $P$ is the pullback of $R$ and $S$ over $B$
	and the arrow $P\to A\times C$ is induced by the arrows $P\to A$ and $P\to C$ in the following diagram. 
	\begin{equation}\label{eq:relcompo}
		\begin{tikzcd}[column sep=small, row sep=small]
			&
			&
		P
			\ar[ld]
			\ar[rd]
			\ar[dd,phantom,"\mbox{\rotatebox[origin=c]{45}{$\lrcorner$}}"{marking,near start}]
			&
			&
			\\
			&
		R
			\ar[ld]
			\ar[rd]
			&
			&
		S
			\ar[ld]
			\ar[rd]
			&
			\\	
		A
			&
			&		
		B
			&
			&	
		C
		\end{tikzcd}		
	\end{equation}
	The unit on $A$ is the $\one{M}$-image of the diagonal $\Delta_A=\langle\id_A,\id_A\rangle\colon A\to A\times A$.
	The composition of cells is defined using the universal property of pullbacks and the orthogonality of the factorisation system.
	\end{itemize}
\end{definition}

\begin{remark} 
	\label{rem:leftclassrelation}
In the literature (e.g., \cite{Kl70,Kel91,Pav95}),
many authors have studied the bicategories of $\zero{M}$-relations,
to which we bring the vertical part to obtain double categories. 
It is occasionally useful to regard $\zero{M}$-relations as equivalence relations of spans in $\one{C}$,
where the equivalence relation is defined via the left class $\zero{E}$ of the factorisation system,
as shown in \cite{HNST22}.
Related to this, the paper \cite{Kaw73} discusses the category of equivalence classes of spans where 
the equivalence relation is defined by a class of epimorphisms, 
which is not necessarily a left class of a factorisation system.
In this paper, we do not pursue this direction.
\end{remark}

\begin{remark}\label{rem:stabasso}
	The stability of a factorisation system $(\zero{E}, \zero{M})$ ensures the 
	associativity and unitality of the horizontal composition of $\zero{M}$-relations.
	To see this, observe that for any composable sequence of horizontal arrows whose data are given by
	 $\langle l_i,r_i\rangle \colon A_{i-1}\sto A_i\quad (i=1,\dots,n)$,
	the composite is given by the $\zero{M}$-image of the arrow $Q\to A_0\times A_n$ in $\one{C}$,
	where $Q$ is the multiple pullback of $R_1,\dots,R_n$ over $A_1,\dots,A_{n-1}$,
	(or equivalently, the composition in the double category of spans $\Span(\one{C})$),
	not depending on the order of composition.
	For instance, take the case $n=3$.
	In the following diagrams,
	every square is a pullback square in $\one{C}$.
	\[
		\begin{tikzcd}[column sep=small, row sep=small]
			&
			&
			&
		Q
			\ar[ld]
			\ar[rrdd]
			\ar[dd,phantom,"\mbox{\rotatebox[origin=c]{45}{$\lrcorner$}}"{marking,near start}]
			&
			&
			\\
			&
			&
		P
		\ar[dd,phantom,"\mbox{\rotatebox[origin=c]{45}{$\lrcorner$}}"{marking,near start}]
			\ar[ld]
			\ar[rd]
			&
			&
			&
			\\
			&
		R_1
			\ar[ld]
			\ar[rd]
			&
			&
		R_2
			\ar[ld]
			\ar[rd]
			&
			&
		R_3
			\ar[ld]
			\ar[rd]
			&
			&
			\\	
		A_0
			&
			&		
		A_1
			&
			&	
		A_2
			&
			&
		A_3
		\end{tikzcd}
		,
		\hspace{2ex}
		\begin{tikzcd}[column sep=small, row sep=small]
			&
			Q 
			\ar[ld]
			\ar[rd,two heads]
			\ar[dd,phantom,"\mbox{\rotatebox[origin=c]{45}{$\lrcorner$}}"{marking,near start}]
			&
			&
			\\
		P 
			\ar[dr,two heads]
			&
			&
		S 
			\ar[ld]
			\ar[rd]
			\ar[dd,phantom,"\mbox{\rotatebox[origin=c]{45}{$\lrcorner$}}"{marking,near start}]
			&
			\\
			&
		R_1R_2
			\ar[ld]
			\ar[rd]
			&
			&
		R_3
		\ar[ld]
		\ar[rd] 
		\\
		A_0 
			&
			&
		A_2
			&
			&
		A_3
		\end{tikzcd}
	\]
	The composite of $R_1$ and $R_2$ is the $\zero{M}$-image of $P\to A_0\times A_2$,
	where $P$ is the pullback of $R_1$ and $R_2$ over $A_1$,
	and the composite $(R_1R_2)R_3$ is given by the $\zero{M}$-image of the arrow $S\to A_0\times A_3$,
	where $S$ is the pullback of $R_1R_2$ and $R_3$ over $A_2$.
	By the pullback lemma, $Q$ exhibits itself as a pullback of $P$ and $S$ over $R_1R_2$,
	and hence the arrow $Q\to S$ is in $\zero{E}$ by the stability of the factorisation system.
	By the uniqueness of factorisation,
	the composite $(R_1R_2)R_3$ is the same, up to isomorphism, as the $\zero{M}$-image of the arrow $Q\to A_0\times A_3$ is defined from the diagrams.
	The unitality of the horizontal composition is also proved in a similar way.
	For the horizontal composition of cells,
	the associativity and the unitality are easily verified from the uniqueness part of the orthogonality of the factorisation system.
\end{remark}

\begin{remark}
	\label{lem:FrConstruction}
	These double categories should be better understood through the $\dbl{F}\mathrm{r}$-construction
	established in \cite[Theorem 14.4]{Shu08}.
	The counterpart of the $\dbl{F}\mathrm{r}$-construction in the context of
	proarrow equipments can be found in \cite[Proposition 4.2.5]{Law15}.

	The theorem in \cite{Shu08} states as follows.
	Let $\one{B}$ be a category with finite limits
	and $F\colon \one{E}\to\one{B}$ be a bifibration with fibred products.
	Suppose, moreover, that $F$ satisfies the following.
	\begin{itemize}
		\item%
			The binary product functor $-\times-\colon\one{E}\times\one{E}\to\one{E}$
			induced by finite products in the base category and fibred products 
			preserves opcartesian arrows.
		\item%
			$F$ is \textit{strongly \textsc{bc}}:
			for each pullback square in the base category $\one{B}$,
			the mate of the canonical natural isomorphism between reindexing functors gives rise to an isomorphism. 
			\[
				\begin{tikzcd}
					A
					\ar[r,"f"]
					\ar[d,"h"']
					\pullback
						&
						C
						\ar[d,"k"]
					\\
					B
					\ar[r,"g"']
						&
						D
				\end{tikzcd}
				\hspace{3ex}
				\begin{tikzcd}
					\one{E}_A
					\al[r,"f^*"]
					\ar[d,"\Sigma_h"']
					\doublecell[rd]{\cong}
						&
						\one{E}_C
						\ar[d,"\Sigma_k"]
					\\
					\one{E}_B
					\al[r,"g^*"']
						&
						\one{E}_D
				\end{tikzcd}
			\]
	\end{itemize}
	Then, there is an equipment $\dbl{F}\mathrm{r}(F)$ whose accompanying fibration is given by
	the base change
	\[
	\begin{tikzcd}
		\dbl{F}\mathrm{r}(F)_1 
		\ar[r]
		\ar[d]
		\pullback
			&
			\one{E}
			\ar[d, "F"]
		\\
		\one{B}\times\one{B}
		\ar[r, "\times"']
			&
			\one{B}
	\end{tikzcd}
	\hspace{1ex}
	\text{,}	
	\]
	and its horizontal units and compositions are given in the following way, where 
	$\Sigma_\mhyphen$ denotes the opcartesian lift.
	\begin{itemize}
		\item The horizontal unit on $A$ is $\Sigma_{\Delta_A}(!_A)^*\top \in \one{E}_{A\times A}$,
		where $\top \in \one{E}$ is the terminal object in the fibre on the terminal object in $\one{B}$,
		$!_A\colon A\to 1$ is the unique arrow to the terminal object in $\one{B}$,
		and $\Delta_A\colon A\to A\times A$ is the diagonal arrow.
		\item The horizontal composition of $R\colon A\sto B$ and $S\colon B\sto C$ is given by
		$\Sigma_{\pi_{A,C}}(\pi_{\id_A\times\Delta_B\times\id_C})^*(R\times S)$,
		where $R\times S$ is the product of $R$ and $S$ in $\one{E}$ over $A\times B\times B\times C$,
		and $\pi_{A,C}\colon A\times B\times C\to A\times C$ is the projection.
	\end{itemize}

	For the fibration $\one{M}\to\one{C}$ in \cref{def:RelEM},
	an arrow $(f,g)\colon m\to m'$ in $\one{M}$ is opcartesian if and only if $f$ is in the left class,
	and
	such arrows are closed under product since $\zero{E}$ is stable under pullback.
	The strong \textsc{bc} condition
	is satisfied because $\one{M}$ is a reflective subfibration, as we have remarked in \Cref{lem:FrConstruction}.
	Therefore, we can apply this theorem to the fibration $\one{M}\to\one{C}$
	and we obtain the equipment $\dbl{F}\mathrm{r}(\one{M}\to\one{C})=\Rel{\zero{E}, \zero{M}}(\one{C})$.
\end{remark}	

\begin{proposition}
	\label{prop:RelEMisCartBC}
	Let $\one{C}$ be a category with finite limits and
	$(\zero{E}, \zero{M})$ be a stable factorisation system in $\one{C}$.
	Then, $\Rel{\zero{E}, \zero{M}}(\one{C})$ is a cartesian equipment with Beck-Chevalley pullbacks.	
\end{proposition}

\begin{proof}
	Firstly, we show that $\Rel{\zero{E}, \zero{M}}(\one{C})$ is a cartesian double category,
	where finite products are given by the product in $\one{C}$.
	The product functor $\times\colon \Rel{\zero{E}, \zero{M}}(\one{C})\times\Rel{\zero{E}, \zero{M}}(\one{C})\to\Rel{\zero{E}, \zero{M}}(\one{C})$
	and the terminal object functor $1\colon \mathbbm{1}\to\Rel{\zero{E}, \zero{M}}(\one{C})$ are pseudo-functors
	because the products in $\one{C}$ are compatible with the factorisations of arrows in $\one{C}$.

	Secondly, $\Rel{\zero{E}, \zero{M}}(\one{C})$ is an equipment.
	Indeed, it is constructed by the $\dbl{F}\mathrm{r}$-construction,
	which is proven to create an equipment in \cite[Theorem 14.4]{Shu08}.
	The explicit description is given as follows.
	The companion $f_!$ of a vertical arrow $f\colon A\to B$ is given 
	by the $\zero{M}$-image of the graph of $f$, i.e., an arrow
	$\langle \id_A, f\rangle\colon A\to A\times B$,
	and the conjoint $f^*$ is given by the $\zero{M}$-image of the graph $\langle f, \id_A\rangle\colon A\to B\times A$.
	More generally, 
	the restriction of a $\zero{M}$-relation $R\colon A\sto B$ along a vertical arrow $f\colon C\to A$ and $g\colon D\to B$
	is given by the $\zero{M}$-image of the pullback of $R\tto A\times B$
	through $f\times g\colon C\times D\to A\times B$.
	On the other hand, the extension of $R$ along $h\colon A\to X$ and $k\colon B\to Y$ is given by 
	the $\zero{M}$-image of the composite of $R\tto A\times B$ with $h\times k\colon A\times B\to X\times Y$.
	
	Finally, we show that $\Rel{\zero{E}, \zero{M}}(\one{C})$ has Beck-Chevalley pullbacks.
	To show pullback squares are Beck-Chevalley, take a cospan composed of $f\colon A\to D$ and $g\colon B\to D$ in $\one{C}$,
	and consider the following diagram, where every square is a pullback square in $\one{C}$.
	\begin{equation}\label{dgm:ThePullbacksRel}
		\begin{tikzcd}
			C
			\ar[r, "\alpha"]
			\ar[d]
			\ar[rr, bend left, "{\langle h,k\rangle}"]
			\pullback
				&
				R
				\ar[r, "{\langle l_R,r_R\rangle}", tail]
				\ar[d, "\beta"']
				\pullback
					&
					A\times B
					\ar[d,"f\times g"]
			\\
			D
			\ar[r, two heads, "e"]
			\ar[rr, "\Delta"', bend right]
				&
				\cdot
				\ar[r, tail]
					&
					D\times D
		\end{tikzcd}
	\end{equation}
	Here, we factorise $\Delta$ at the bottom of the diagram with respect to $(\zero{E},\zero{M})$.
	Then, $h$ and $k$ become a pair of legs of a pullback square of the cospan in $\one{C}$.
	Moreover, $\alpha$ and $\beta$ define cells of the following form.
	\[
		\begin{tikzcd}[column sep = small]
				&
					C 
					\ar[ld,"h"']
					\ar[rd,"k"]
						&
			\\
			A 
			\sar[rr,"R"']
			\doublecell[rr, shift left=3ex]{\alpha}
			\doublecell[rr, shift right=4ex]{\beta}
			\ar[rd,"f"']
				&
					&
					B
					\ar[ld,"g"]
			\\
				&
				D
					&
		\end{tikzcd}
	\]
	Recall the conditions for cells to be cartesian and opcartesian in $\Rel{\zero{E}, \zero{M}}(\one{C})$ described above.
	Since $\zero{E}$ is stable under pullback, $\alpha$ is in $\zero{E}$
	as an arrow in $\one{C}$ and hence opcartesian as a cell.
	$\beta$ is cartesian because the right-hand square in
	\cref{dgm:ThePullbacksRel} is a pullback square.
	Therefore, the two cells $\alpha$ and $\beta$ constitute a Beck-Chevalley pullback square in $\Rel{\zero{E}, \zero{M}}(\one{C})$.
\end{proof}

\begin{example}
		Let $\one{C}$ be a category with finite limits.
		We have a trivial stable factorisation $(\Iso,\Mor)$ in $\one{C}$ where
		$\Iso$ is the class of all isomorphisms in $\one{C}$ and $\Mor$ is the class of all arrows in $\one{C}$.
		Then, $\Mor$-relations are just spans in $\one{C}$ and $\Rel{\Iso, \Mor}(\one{C})$ is nothing but the double category of spans $\Span(\one{C})$ in $\one{C}$.
\end{example}		
\begin{example}
	\label{ex:RelRegCat}
	Let $\one{C}$ be a regular category.
	We have a stable factorisation $(\Regepi,\Mono)$ in $\one{C}$.
	The double category $\Rel{\Regepi, \Mono}(\one{C})$ is the double category of relations 
	$\Rel{}(\one{C})$ in the sense of \cite{Lam22}.
\end{example}

We now formulate the desirable properties for cartesian equipments
to make them look like $\Rel{\zero{E}, \zero{M}}(\one{C})$.
\begin{definition}
	Let $\dbl{D}$ be a cartesian double category
	and $\zero{M}$ be a class of vertical arrows in $\dbl{D}$.
	Suppose that the following cell exhibits a tabulator $\top R$
	of a horizontal arrow $R\colon A\sto B$:
	\begin{equation}
		\begin{tikzcd}[column sep=small]
				&
				\top R
				\ar[ld, "l"']
				\ar[rd, "r"]
					&
			\\
			A
			\sar[rr, "R"']
			\doublecell[rr, shift left=3ex]{\tau}
				&
					&
					B
		\end{tikzcd}
	\end{equation}
	We say this tabulator is an \emph{$\zero{M}$-tabulator}
	if $\top R\to["\mmbox{\langle l, r\rangle}"]A\times B$ is in $\zero{M}$.
\end{definition}

\begin{lemma}
	\label{lem:rcovtab}
	Let $\dbl{D}$ be a cartesian equipment with tabulators.
	Then, $\dbl{D}$ has $\rclass{\Cov(\dbl{D})}$-tabulators.
\end{lemma}
\begin{proof}
	Let $\langle l,r\rangle \colon  \top R \to A\times B$ exhibit $\top R$ as a tabulator of $R\colon A\sto B$.
	Given an arbitrary commutative diagram in $\oneV{\dbl{D}}$
	\begin{equation}
		\label{eq:lemrdcovtab}
		\begin{tikzcd}
		X 
		\ar[r,"h"]
		\ar[d,"e"']
		&
		\top R 
		\ar[d,"{\langle l,r\rangle}"]
		\\
		Y
		\ar[r,"{\langle f,g\rangle}"']
		&
		A\times B
	\end{tikzcd}		
	\end{equation}
	where $e$ is a cover.
	Then, we gain a cell $\alpha$, as in the diagram below.
	\[
		\begin{tikzcd}[column sep = small,row sep = small]
				&
						X
								\ar[d,"h"]
			\\
				&
						\top R
								\ar[ld,"l"']
								\ar[rd,"r"]
				&
			\\
						A
								\sar[rr,"R"']
								\doublecell[rr, shift left=2ex]{\tab}
				&
				&
						B
		\end{tikzcd}
		\hspace{1ex}
		=
		\hspace{1ex}
		\begin{tikzcd}[column sep = small,row sep = small]
				&
						X
								\ar[d,"e"]
			\\
				&
						Y
								\ar[ld,"f"']
								\ar[rd,"g"]
				&
			\\
						A
								\sar[rr,"R"']
								\doublecell[rr, shift left=2ex]{\alpha}
				&
				&
						B
		\end{tikzcd}
	\]
	The universality of the tabulator $\top R$ implies 
	the existence of
	an arrow $k\colon Y\to\top R$,
	which satisfies $\langle f,g\rangle=k\fatsemi \langle l,r\rangle$,
	and also $e\fatsemi k=h$ again by the uniqueness part of the universality of the tabulator.
	This $k$ is a filler of \cref{eq:lemrdcovtab},
	and one can show by a similar argument that a filler is unique.
\end{proof}

\begin{definition}[{\cite[\S 8]{Lam22}}]
	Let $\dbl{D}$ be a cartesian equipment and
	$\zero{M}$ be a class of vertical arrows in $\dbl{D}$.
	$\dbl{D}$ admits an
	\emph{$\zero{M}$-comprehension scheme}
	if $\dbl{D}$ has strong $\zero{M}$-tabulators
	and, for any morphism of type
	$\langle l,r\rangle\colon X\to A\times B$ in $\zero{M}$,
	the extension
	\begin{equation}
		\begin{tikzcd}[column sep=small]
				&
				X
				\ar[ld, "l"']
				\ar[rd, "r"]
					&
			\\
			A
			\sar[rr, "l^*r_!"']
			\doublecell[rr, shift left=2.5ex]{\opcart}
				&
					&
					B
		\end{tikzcd}
	\end{equation}
	exhibits $X$ as a tabulator of $l^*r_!$.
	An equipment $\dbl{D}$ admits a \emph{left-sided $\zero{M}$-comprehension scheme}
	for a class of arrows $\zero{M}$ in $\oneV{\dbl{D}}$
	if $\dbl{D}$ has left-sided strong $\zero{M}$-tabulators
	and the same condition as above holds for these tabulators;
	i.e., for any morphism $f\colon X\to A$ in $\zero{M}$,
	the opcartesian cell
	\begin{equation}
		\begin{tikzcd}[column sep=small]
				&
				X
				\ar[ld, "f"']
				\ar[rd, "!"]
					&
			\\
			A
			\sar[rr, "f^*!_!"']
			\doublecell[rr, shift left=2.5ex]{\opcart}
				&
					&
					1
		\end{tikzcd}
	\end{equation}
	exhibits $X$ as a tabulator of $f^*!_!$.
\end{definition}

\begin{proposition}
	\label{prop:RelEMhasCompScheme}
	Let $\one{C}$ be a category with finite limits and
	$(\zero{E}, \zero{M})$ be a stable factorisation system in $\one{C}$.
	Then, $\Rel{\zero{E}, \zero{M}}(\one{C})$ admits an $\zero{M}$-comprehension scheme.
\end{proposition}
\begin{proof}
	$\Rel{\zero{E}, \zero{M}}(\one{C})$ has tabulators in an obvious way;
	for a $\zero{M}$-relation $R\colon A\sto B$, we write $\langle l_R,r_R\rangle$
	for the span defining $R$, 
	which is in $\zero{M}$ by definition.
	In $\Rel{\zero{E}, \zero{M}}(\one{C})$, a cell of the following form
	\[
		\begin{tikzcd}[column sep = small]
				&
					C 
					\ar[ld,"f"']
					\ar[rd,"g"]
						&
			\\
			A 
			\sar[rr,"R"']
			\doublecell[rr, shift left=3ex]{\alpha}
				&
					&
						B
		\end{tikzcd}
	\]
	is defined by an arrow $\alpha\colon C\to R$ satisfying
	$\alpha\fatsemi l_R=f$ and $\alpha\fatsemi r_R=g$,
	or equivalently, $\alpha\fatsemi\langle l_R,r_R\rangle=\langle f,g\rangle$.
	Therefore, the span $\langle l_R,r_R\rangle$ gives a tabulator of $R$.
	Furthermore, $\alpha$ is opcartesian if and only if $\alpha$ is in $\zero{E}$
	as an arrow in $\one{C}$,
	since the extension $f^*g_!$ is given by the $\zero{M}$-image of $\langle f,g\rangle$.
	Hence, $\Rel{\zero{E}, \zero{M}}(\one{C})$ has strong tabulators.

	On the other hand, any span $(l,r)$ in $\one{C}$ with $\langle l,r\rangle\colon X\to A\times B$ in $\zero{M}$
	defines an $\zero{M}$-relation $A\sto B$, and the canonical triangle cell
	exhibits $X$ as its tabulator. 
	This concludes that
	$\Rel{\zero{E}, \zero{M}}(\one{C})$ has an $\zero{M}$-comprehension
	scheme.
\end{proof}

\begin{notation}
	For a double category $\dbl{D}$ and a class of vertical arrows $\zero{M}$,
	we write $\one{M}$ for the full subcategory of $\oneV{\dbl{D}}^{\rightarrow}$
	consisting of all vertical arrows in $\zero{M}$,
	similarly in \cref{prop:OFSfromStableSystem}.
	For an object $X\in\dbl{D}$, by $\one{M}\downarrow X$, we mean
	the full subcategory of the slice category $\oneV{\dbl{D}}/X$ consisting of
	vertical arrows in $\zero{M}$.
\end{notation}
\begin{lemma}
	\label{lem:CompScheme}
	Let $\dbl{D}$ be a cartesian equipment and
	$\zero{M}$ be a stable system.
	Suppose that $\dbl{D}$ has strong $\zero{M}$-tabulators.
	Then, $\dbl{D}$ admits an $\zero{M}$-comprehension scheme,
	if and only if
	the functor
	\begin{equation}\label{eq:adjequivforcomp}
		\biH{\dbl{D}}(A, B)\to
		\one{M}\downarrow A\times B
		\text{,}
	\end{equation}
	sending a horizontal arrow to its tabulator is an adjoint equivalence,
	whose pseudo-inverse sends a span to its extension.
\end{lemma}
\begin{proof}
	The adjunction exists when $\dbl{D}$ has $\zero{M}$-tabulators.
	The condition for the unit and the counit to be isomorphism
	is equivalent to the condition of the existence of strong tabulators
	and that every $\zero{M}$-relation $\langle l,r\rangle\colon X\to A\times B$ exhibits $X$ as a tabulator of its extension 
	$l^*r_!$.
\end{proof}

\begin{remark}
	\label{rem:Adjexsist}
Even if $\dbl{D}$ does not admit an $\zero{M}$-comprehension scheme but has $\zero{M}$-tabulators,
we can still define an adjunction between the two categories above, as shown in the proof.
Alternatively, an $\zero{M}$-comprehension scheme is a minimum structure that makes the adjunction an equivalence. 
\end{remark}

\begin{proposition}
	\label{prop:compschemeunbiased}
	The following are equivalent for a discrete cartesian equipment $\dbl{D}$
	and a stable system $\zero{M}$.
	\begin{enumerate}
		\item
			$\dbl{D}$ has strong $\zero{M}$-tabulators.
		\item
			$\dbl{D}$ has left-sided strong $\zero{M}$-tabulators.
	\end{enumerate}
	Also, the following are equivalent.
	\begin{enumerate}
		\item
			$\dbl{D}$ admits an $\zero{M}$-comprehension scheme.
		\item
			$\dbl{D}$ admits a left-sided $\zero{M}$-comprehension scheme.
	\end{enumerate}
\end{proposition}

\begin{proof}
	The statement is a direct consequence of \cref{lem:UnbiasingLemma}.
\end{proof}

\subsection{The characterisation theorem for double categories of relations}\label{subsec:Characterisation}

\begin{lemma}\label{lem:idTab}
	We have the following tabulator in any equipment $\dbl{D}$ with the vertical terminal $1$.
	\begin{equation}
		\label{eq:idTab}
		\begin{tikzcd}[column sep=small]
				&
				A 
				\ar[ld, equal]
				\ar[rd, "!"]
					&
			\\
			A
			\sar[rr, "!_!"']
			\doublecell[rr, shift left=2.5ex]{\tab}
			&
					&
					1
		\end{tikzcd}
	\end{equation}
\end{lemma}
\begin{proof}
	There is a unique cell of the form on the left below
	because it corresponds to a cell of the form on the right,
	which is unique since $1$ is the vertical terminal.
	\[
		\begin{tikzcd}[column sep=small]
				&
				X
				\ar[ld, "f"']
				\ar[rd, "!"]
					&
			\\
			A
			\sar[rr, "!_!"']
			\doublecell[rr, shift left=2.5ex]{\alpha}
				&
					&
					1
		\end{tikzcd}
		\hspace{3ex}
		,
		\hspace{3ex}
		\begin{tikzcd}[column sep=small]
				&
				X
				\ar[ld, "f"']
				\ar[rd, "!"]
					&
			\\
			A
			\doublecell[rr, ]{\Id_!}
			\ar[rd, "!"']
				&
					&
					1
					\ar[ld,equal]
			\\
				&
				1
					&
		\end{tikzcd}
	\]
	Additionally, $f$ is the only vertical arrow that is composed with the cell \Cref{lem:idTab} to give the unique cell $\alpha$. 
	This proves the cell in question to be the tabulator.
\end{proof}

\begin{definition}
	In a double category $\dbl{D}$ with the vertical terminal $1$,
	we say a vertical morphism $f\colon A \to X$ is a \emph{fibration}
	if there exists a horizontal arrow $p\colon X \sto 1$ and a tabulating cell
	\[
		\begin{tikzcd}[column sep=small]
				&
				A
				\ar[ld, "f"']
				\ar[rd, "!"]
					&
			\\
			X
			\sar[rr, "p"']
			\doublecell[rr, shift left=2.5ex]{\tab}
				&
					&
					1
		\end{tikzcd}.
	\]

	We write $\Fib(\dbl{D})$ for the class of fibrations in $\dbl{D}$.
\end{definition}

\begin{remark}
	\label{rem:fibandcompr}
	If $\dbl{D}$ has strong tabulators, then, for a fibration $f\colon A\to X$ in the definition above,
	$p$ is uniquely determined as $f^*!_!$ up to isomorphism.
	Therefore, a cartesian equipment $\dbl{D}$ admits a left-sided $\zero{M}$-comprehension scheme for a stable system $\zero{M}$
	if it has strong $\zero{M}$-tabulators and $\zero{M}\subset\Fib(\dbl{D})$ holds.   

	We see this class as a candidate for the right class of an orthogonal factorisation system,
	but it is not closed under composition in general.
\end{remark}

\begin{example}
	In $\Prof$, fibrations are precisely discrete fibrations,
	while in $\spFib$, fibrations are defined similarly and coincide with split fibrations in the ordinary sense.
\end{example}

\begin{lemma}\label{lem:FibStable}
	In any equipment $\dbl{D}$ with the vertical terminal object,
	$\Fib(\dbl{D})$ is stable under pullbacks.
\end{lemma}
\begin{proof}
	Suppose there is a pullback in $\oneV{\dbl{D}}$ of the following form, and 
	$f$ exhibits $A$ as a tabulator of $p\colon X\sto 1$.
	\[
		\begin{tikzcd}
			B
			\ar[r, "h"]
			\ar[d, "f'"']
			\pullback
				&
				A
				\ar[d, "f"]
			\\
			Y
			\ar[r, "g"']
				&
				X
		\end{tikzcd}
		\hspace{2ex}
		,
		\hspace{2ex}
		\begin{tikzcd}[column sep=small]
				&
				A
				\ar[ld, "f"']
				\ar[rd, "!"]
					&
			\\
			X
			\sar[rr, "p"']
			\doublecell[rr, shift left=2.5ex]{\tab}
				&
					&
					1
		\end{tikzcd}
	\]
	Define a cell $\bar{h}$ by the following equality.
	We show that $\bar{h}$ exhibits $B$ as the tabulator of $g_!p$.
	\[
		\begin{tikzcd}[column sep=small]
				&
				B
				\ar[d, "h"]
					&
			\\
				&
				A
				\ar[ld, "f"']
				\ar[rd, "!"]
					&
			\\
			X
			\sar[rr, "p"']
			\doublecell[rr, shift left=2.5ex]{\tab}
				&
					&
					1
		\end{tikzcd}
		\hspace{2ex}
		=
		\hspace{2ex}
		\begin{tikzcd}[column sep=small]
				&
				B
				\ar[ld, "f'"']
				\ar[rd, "!"]
					&
			\\
			Y
			\sar[rr, "g_!p"']
			\ar[d, "g"']
			\doublecell[rr, shift left=3.0ex]{\bar{h}}
			\doublecell[rrd]{\cart}
				&
					&
					1
					\ar[d, equal]
			\\
			X
			\sar[rr, "p"']
				&
					&
					1
		\end{tikzcd}
	\]
	Fix a vertical arrow $y\colon Z\to Y$.
	Firstly, the cartesian cell defining $g_!p$
	assures that a cell $\alpha_1$ of the form
	on the left below corresponds to a cell $\alpha_2$.
	Secondly, such a cell corresponds to a vertical arrow
	$\alpha_3$ below, making the diagram commute through postcomposing
	the tabulator of $p$.
	Finally, by the pullback assumed at the beginning,
	$\alpha_3$ corresponds to another
	vertical arrow $\alpha_4$ under the equality $\alpha_3=\alpha_4\fatsemi h$.
	\[
		\begin{tikzcd}[column sep=small]
				&
				Z
				\ar[ld, "y"']
				\ar[rd, "!"]
					&
			\\
			Y
			\sar[rr, "g_!p"']
			\doublecell[rr, shift left=2.5ex]{\alpha_1}
				&
					&
					1
		\end{tikzcd}
		\hspace{2ex}
		\vline 
		\,
		\vline
		\hspace{2ex}
		\begin{tikzcd}[column sep=small]
				&
				Z
				\ar[ld, "y"']
				\ar[rd, "!"]
					&
			\\
			Y
			\ar[d, "g"']
				&
					&
					1
					\ar[d, equal]
			\\
			X
			\sar[rr, "p"']
			\doublecell[rr, shift left=7.0ex]{\alpha_2}
				&
					&
					1
		\end{tikzcd}
		\hspace{2ex}
		\vline 
		\,
		\vline
		\hspace{2ex}
		\begin{tikzcd}[column sep=small]
				&
				Z
				\ar[ld, "y"']
				\ar[rd, "!"]
				\ar[d, "\alpha_3"description]
					&
			\\
			Y
			\ar[d, "g"']
				&
				A
				\ar[ld, "f"]
				\ar[rd, "!"']
					&
					1
					\ar[d, equal]
			\\
			X
				&
					&
					1
		\end{tikzcd}
		\hspace{2ex}
		\vline 
		\,
		\vline
		\hspace{2ex}
		\begin{tikzcd}[column sep=small]
				&
				Z
				\ar[ldd, "y"', bend right]
				\ar[rdd, "!", bend left]
				\ar[d, "\alpha_4"description]
					&
			\\
				&
				B
				\ar[ld, "f'"]
				\ar[rd, "!"']
					&
			\\
			Y
				&
					&
					1
		\end{tikzcd}
	\]
	Tracing back the sequence of correspondences,
	we see that the bijective correspondence
	$\alpha_4\mapsto\alpha_1$ is obtained by postcomposing $\bar{h}$.
\end{proof}

Now we remark here that with the above lemma, we can rephrase left-sided $\zero{M}$-comprehension scheme by composability of fibrations:
\begin{lemma}
	\label{lem:leftsidedComp_vs_ComposableFib}
	Let $\dbl{D}$ be a double category such that
	$\oneV{\dbl{D}}$ has finite limits.
	The following are equivalent.
	\begin{enumerate}
		\item %
			There exists a stable system $\zero{M}$ on $\oneV{\dbl{D}}$ 
			and $\dbl{D}$ has left-sided $\zero{M}$-comprehension scheme.
		\item %
			$\dbl{D}$ has left-sided strong tabulators
			and $\Fib(\dbl{D})$ is closed under composition.
	\end{enumerate}
	If the first condition holds, the class $\zero{M}$ is precisely $\Fib(\dbl{D})$.
\end{lemma}

\begin{proof}
	If the first condition holds, then $\zero{M}=\Fib(\dbl{D})$ by
	\Cref{rem:fibandcompr} since a morphism in $\Fib(\dbl{D})$ arises as a $\zero{M}$-tabulator of 
	a horizontal arrow to the terminal object.
	This makes $\Fib(\dbl{D})$ closed under composition.
	Conversely, if the second condition holds, then $\Fib(\dbl{D})$ is a stable system by \Cref{lem:FibStable}.
	Note that an isomorphism arises as a tabulator of the terminal horizontal arrow.
	By definition and the first half of the second condition, 
	$\dbl{D}$ has left-sided $\Fib(\dbl{D})$-comprehension scheme.
\end{proof}

In search of an orthogonal factorisation system whose right class is
$\Fib(\dbl{D})$, we should recall the final-functors/discrete-fibrations 
factorisation in the category of categories, also known as the
comprehensive factorisation.
As we shall see afterwards, final functors are characterised in $\Prof$ 
solely by its double categorical structure.
This observation incites us to consider the following definition.

\begin{definition}
	Let $\dbl{D}$ be an equipment with the vertical terminal object $1$,
	and $f\colon A\to X$ be a vertical arrow.
	We say $f$ is \emph{final} if the identity cell of the form below,
	which is the unique cell of this form because of the universality of the terminal $1$,
	satisfies the Beck-Chevalley condition.
	\[
		\begin{tikzcd}[column sep=small]
				&
				A
				\ar[ld, "f"']
				\ar[rd, "!"]
				\doublecell[dd]{=}
					&
			\\
			X
			\ar[rd, "!"']
				&
					&
					1
					\ar[ld, equal]
			\\
				&
				1
					&
		\end{tikzcd}
	\]
	We write $\Fin(\dbl{D})$ for the class of final morphisms.
\end{definition}
\begin{example}
	In $\Prof$, a functor $f\colon A \to X$ is final in the above sense if and only if
	the colimit $\colim_{a\in A} X(x, f(a))$ is terminal in $\Set$
	for each $x\in X$, which means that the comma category $x\downarrow f$
	is connected.
	This is equivalent to the condition that $f$ is a final functor in the ordinary sense.
	There are other known characterisations of final functors in $\Prof$ as in \cite[Theorem 4.67]{Kel82},
	and some of them can be derived directly from this definition of finality in a well-behaved double category.
	See \cite[Example 4.3]{Kou14} for more details.
\end{example}
\begin{lemma}\label{lem:CovInFin}
	Covers are final in any equipment $\dbl{D}$ with the vertical terminal object $1$.
\end{lemma}
\begin{proof}
	Let $e\colon A\to X$ be a cover. We have the following opcartesian and cartesian
	cells, which proves that $e$ is final.
	\[
		\begin{tikzcd}[column sep=small]
				&
				A
				\ar[ld, "e"']
				\ar[rd, "e"]
					&
			\\
			X
			\ar[d, equal]
			\sar[rr, equal]
			\doublecell[rrd]{\opcart}
			\doublecell[rr, shift left=2ex]{\opcart}
				&
					&
					X
					\ar[d, "!"]
			\\
			X
			\ar[rd, "!"']
			\sar[rr]
			\doublecell[rr, shift right=3ex]{\cart}
				&
					&
					1
					\ar[ld, equal]
			\\
				&
				1
					&
		\end{tikzcd}
	\]
\end{proof}
\begin{lemma}
	\label{lem:FinStable}
	Let $\dbl{D}$ be an equipment with the vertical terminal $1$.
	Suppose that there exists a cell $\alpha$ of the following form
	satisfying the Beck-Chevalley condition.
	\[
		\begin{tikzcd}
				&
				\cdot
				\ar[ld, "e"']
				\ar[rd, "h"]
				\doublecell[dd]{\alpha}
					&
			\\
			\cdot
			\ar[rd, "g"']
				&
					&
					\cdot
					\ar[ld, "f"]
			\\
				&
				\cdot
					&
		\end{tikzcd}
	\]
	Then, if $f$ is final, $e$ is as well.
	In particular, the class of final arrows $\Fin(\dbl{D})$  
	are stable under pullbacks if $\dbl{D}$
	has Beck-Chevalley pullbacks.
\end{lemma}
\begin{proof}
	Utilising \cref{lem:Sandwich},
	we achieve an opcartesian cell exhibiting the extension of
	$e$ and $!$ as follows.
	\[
		\begin{tikzcd}
				&
					&
					\cdot
					\ar[ld,"e"']
					\ar[rd,"h"]
						&
							&
			\\
				&
				\cdot
				\ar[ld,equal]
				\ar[rd,"g"description]
				\doublecell[rr, shift left=3ex]{\opcart}
				\doublecell[rr, shift right=3ex]{\cart}
				\sar[rr]
					&
						&
						\cdot
						\ar[ld,"f"description]
						\ar[rd,"!"]
								&
			\\
			\cdot
			\sar[rr, "g_!"']
			\doublecell[rr, shift left=2.5ex]{\opcart}
				&
					&
					\cdot
					\sar[rr,"!_!"']
					\doublecell[rr, shift left=2.5ex]{\opcart}
						&
							& 
							\hspace{2.0ex}1\hspace{2.0ex}
		\end{tikzcd}
	\]
	Since $g_!!_!\cong!_!$, this shows $e$ is final.
\end{proof}
\begin{lemma}\label{lem:CovEqualFin}
	Let $\dbl{D}$ be a discrete cartesian equipment.
	Then $\Fin(\dbl{D})$ coincides with $\Cov(\dbl{D})$.
\end{lemma}
\begin{proof}
	Since we have already shown $\Cov(\dbl{D})\subseteq\Fin(\dbl{D})$
	in \cref{lem:CovInFin},we show the converse.
	Suppose $f\colon A\to X$ is a final and consider the following cells,
	in which the horizontal arrow $\varepsilon$ is the same as the one in \Cref{prop:discselfdual}.
	Then we obtain the whole opcartesian cell below showing that
	$\varepsilon$ is presented by $\langle f,f\rangle$.
	\[
		\begin{tikzcd}
					&
					A
					\ar[ld,"f"']
					\ar[rd,"!"]
						&
							&
			\\
				\hspace{2.0ex}X\hspace{2.0ex}
				\ar[d,"\Delta"']
				\ar[rrd,"!"description]
				\doublecell[rr, shift left=3ex]{\opcart}
				\doublecell[rr, shift right=2ex,xshift=3ex]{\cart}
				\sar[rr]
					&
						&
						\hspace{2.0ex}1\hspace{2.0ex}
						\ar[d,equal]
								&
			\\
			X\times X
			\sar[rr, "\varepsilon"']
			\doublecell[rr, shift left=2ex,xshift=-3.5ex]{\opcart}
				&
					&
					\hspace{2.0ex}1\hspace{2.0ex}
						&
		\end{tikzcd}
	\]
	Recalling the definition of $\varepsilon_f$,
	we know that this cell is exactly the cell corresponding to
	the following vertical identity cell by \cref{lem:UnbiasingLemma}.
	\[
		\begin{tikzcd}[column sep = small]
				&
				A
				\ar[ld,"f"']
				\ar[rd,"f"]
					&
			\\
			\hspace{2.0ex}X\hspace{2.0ex}
			\doublecell[rr, shift left=3ex]{\Id_f}
			\sar[rr, equal]
				&
					&
					\hspace{2.0ex}X\hspace{2.0ex}
		\end{tikzcd}
	\]
	Therefore, \cref{lem:UnbiasingLemma} shows that this cell is
	opcartesian, which means that $f$ is a cover.
\end{proof}
\begin{remark}
	In general, finals do not coincide with covers.
	In $\Prof$, take a category $\one{C}$ with a terminal object $1$ that is not an initial object.
	Now the functor $\ulcorner1\urcorner\colon\one{1}\to\one{C}$ from the terminal category is not initial but final.
	Such a functor never be absolutely dense, because by applying \cref{lem:CovInFin} to $\Prof$ and $\Prof^\hop$, we can show that any absolutely dense functor is final and initial (see also \cite[Proposition 3.6]{EV02}).
\end{remark}

The following theorem generalises the \textit{comprehensive factorisation} observed in \cite{SW73} and \cite[\S 4.7]{Kel82}.
\begin{theorem}
	\label{thm:CompSchemeToFactorisation}
	Let $\dbl{D}$ be an equipment such that $\oneV{\dbl{D}}$ is finitely complete.
	Suppose that we are given a stable system $\zero{M}$ on $\oneV{\dbl{D}}$ and that
	$\dbl{D}$ admits a left-sided $\zero{M}$-comprehension scheme.
	Then $(\Fin(\dbl{D}),\zero{M})$ is an
	orthogonal factorisation system on $\oneV{\dbl{D}}$.
	Moreover, if $\dbl{D}$ has Beck-Chevalley pullbacks, then
	the factorisation system is stable.
\end{theorem}
\begin{remark}
	Using \cref{lem:leftsidedComp_vs_ComposableFib}, we can state this theorem without $\zero{M}$:
	$\zero{M}$ in the conclusion is replaced by $\Fib(\dbl{D})$, and
	the second condition of the theorem is replaced by the following:
	$\dbl{D}$ has left-sided strong tabulators
	and $\Fib(\dbl{D})$ is closed under composition.
\end{remark}
\begin{proof}[Proof of \cref{thm:CompSchemeToFactorisation}]
	Owing to \cref{prop:OFSfromStableSystem},
	it suffices to show that $\one{M}\hto\cod^{\oneV{\dbl{D}}}$
	has the left adjoint and the class of all vertical arrows left orthogonal to all arrows in $\zero{M}$ is exactly $\Fin(\dbl{D})$.

	Note that
	$\one{M}\hto\cod^{\oneV{\dbl{D}}}$ is a fibred functor since
	$\zero{M}$ is a stable system.
	\cite[Theorem 5.3.7]{Vas14} shows that such an adjunction between
	fibrations is checked fibrewise.
	Therefore, to establish this,
	we need to show that for each object $B\in \dbl{D}$,
	the inclusion $\one{M}\downarrow B\to[hook]\oneV{\dbl{D}}/B$ has a left adjoint.

	In other words, we need to prove that
	for each vertical arrow $f \colon  A \to B$, there
	exists an arrow $m \colon  C \to B$ in $\zero{M}$ equipped with
	a \textit{unit} $e \colon  f \to m$ in $\oneV{\dbl{D}}/B$ that is universal.
	Specifically, we show that, for each $n\in \one{M}\downarrow B$, the function
	$\left(\one{M}\downarrow B\right)(m, n)\to\left(\oneV{\dbl{D}}/B\right)(f, n)$
	obtained by precomposing $e$ is a bijection.

	Let $f\colon A\to B$ be an arrow in $\oneV{\dbl{D}}$.
	We take the extension of the span $(f,!)\colon B\to[leftarrow]A\to 1$
	and then obtain the tabulator $C$ of the extension, as shown below.
	\begin{equation}\label{dgm:TheFactorisationFin}
		\begin{tikzcd}[column sep=small]
				&
				A
				\ar[ldd, "f"', bend right=20]
				\ar[rdd, "!", bend left=20]
				\ar[d,"e"]
					&
			\\
				&
				C 
				\ar[ld, "m"']
				\ar[rd, "!"]
					&
			\\
			B
			\sar[rr, "f^*!_!"']
			\doublecell[rr, shift left=2.5ex]{\tab}
				&
					&
					1
		\end{tikzcd}
	\end{equation}
	Since $\dbl{D}$ has left-sided strong $\zero{M}$-tabulators,
	$m$ is in $\zero{M}$.

	Let $n \colon  D \to B$ be a vertical arrow in $\zero{M}$.
	By the left-sided $\zero{M}$-comprehension scheme,
	we obtain the following tabulator.
	\begin{equation}\label{dgm:TabFor_nFin}
		\begin{tikzcd}[column sep=small]
				&
				D 
				\ar[ld, "n"']
				\ar[rd, "!"]
					&
			\\
			B
			\sar[rr, "n^*!_!"']
			\doublecell[rr, shift left=2.5ex]{\tab}
			&
					&
					1
		\end{tikzcd}
	\end{equation}
	Therefore,
	an arrow $\bar{u}\colon m \to n$ in $\one{M}\downarrow B$
	bijectively corresponds to a cell of
	the form on the left below.
	Furthermore, such a cell corresponds to a horizontal cell of the form in
	the middle below by considering the strong tabulator
	in \cref{dgm:TheFactorisationFin}.
	The opcartesian cell defining the extension $f^*!_!$
	ensures that this bijectively corresponds to a cell
	of the form on the right.
	\[
		\begin{tikzcd}[column sep=small]
				&
				C
				\ar[ld, "m"']
				\ar[rd, "!"]
					&
			\\
			B
			\sar[rr, "n^*!_!"']
			\doublecell[rr, shift left=2.5ex]{\gamma_1}
			&
					&
					1
		\end{tikzcd}
		\hspace{2ex}
		\vline 
		\,
		\vline
		\hspace{2ex}
		\begin{tikzcd}
			B
			\sar[r, "f^*!_!"]
			\ar[d, equal]
			\doublecell{\gamma_2}
				&
				1
				\ar[d, equal]
			\\
			B
			\sar[r, "n^*!_!"']
				&
				1
		\end{tikzcd}
		\hspace{2ex}
		\vline 
		\,
		\vline
		\hspace{2ex}
		\begin{tikzcd}[column sep=small]
				&
				A
				\ar[ld, "f"']
				\ar[rd, "!"]
					&
			\\
			B
			\sar[rr, "n^*!_!"']
			\doublecell[rr, shift left=2.5ex]{\gamma_3}
			&
					&
					1
		\end{tikzcd}
	\]
	Finally, by considering the tabulator \cref{dgm:TabFor_nFin} again,
	this corresponds
	to a vertical arrow $u\colon A\to D$.
	It remains to show the bijection $\bar{u}\mapsto u$ is obtained
	by precomposing $e$ to $\bar{u}$. 
	By tracing the above correspondence, we obtain
	the following equations.
	\[
		\begin{tikzcd}[column sep=small]
				&
				C
				\ar[d, "\bar{u}"]
					&
			\\
				&
				D
				\ar[ld, "n"']
				\ar[rd, "!"]
					&
			\\
			B
			\sar[rr, "n^*!_!"']
			\doublecell[rr, shift left=2.5ex]{\tab}
				&
						&
						1
		\end{tikzcd}
		\hspace{1ex}
		=
		\hspace{1ex}
		\begin{tikzcd}[column sep=small]
				&
				C
				\ar[ld, "m"']
				\ar[rd, "!"]
					&
			\\
			B
			\sar[rr, "n^*!_!"']
			\doublecell[rr, shift left=2.5ex]{\gamma_1}
			&
					&
					1
		\end{tikzcd}
		\hspace{1ex}
		=
		\hspace{1ex}
		\begin{tikzcd}[column sep=small]
				&
				C
				\ar[ld, "m"']
				\ar[rd, "!"]
					&
			\\
			B
			\sar[rr, "f^*!_!"']
			\doublecell[rr, shift left=2.5ex]{\tab}
			\ar[d, equal]
			\doublecell[rrd]{\gamma_2}
				&
						&
						1
					\ar[d, equal]
			\\
			B
			\sar[rr, "n^*!_!"']
				&
						&
						1
		\end{tikzcd}
	\]
	By precomposing $e$ to these diagrams, we obtain the following equations,
	hence verifying $e\fatsemi \bar{u}=u$.
	\[
		\begin{tikzcd}[column sep=small]
				&
				A
				\ar[d, "e\fatsemi \bar{u}"]
					&
			\\
				&
				D
				\ar[ld, "n"']
				\ar[rd, "!"]
					&
			\\
			B
			\sar[rr, "n^*!_!"']
			\doublecell[rr, shift left=2.5ex]{\tab}
				&
						&
						1
		\end{tikzcd}
		\hspace{1ex}
		=
		\hspace{1ex}
		\begin{tikzcd}[column sep=small, row sep=small]
				&
				A
				\ar[d, "e"]
					&
			\\
				&
				C
				\ar[ld, "m"']
				\ar[rd, "!"]
					&
			\\
			B
			\sar[rr, ]
			\doublecell[rr, shift left=2.5ex]{\tab}
			\ar[d, equal]
			\doublecell[rrd]{\gamma_2}
				&
						&
						1
					\ar[d, equal]
			\\
			B
			\sar[rr, "n^*!_!"']
				&
						&
						1
		\end{tikzcd}
		\hspace{1ex}
		=
		\hspace{1ex}
		\begin{tikzcd}[column sep=small]
				&
				A
				\ar[ld, "f"']
				\ar[rd, "!"]
					&
			\\
			B
			\sar[rr, "n^*!_!"']
			\doublecell[rr, shift left=2.5ex]{\gamma_3}
			&
					&
					1
		\end{tikzcd}
	\]
	Here, the last equation follows from the definition of the tabulator in
	\cref{dgm:TheFactorisationFin}.

	Now, we show the left class precisely consists of final arrows.
	The vertical arrows $e$ and $m$ in \cref{dgm:TheFactorisationFin} give the factorisation
	for each vertical arrow $f \colon A\to B$
	as seen in the proof of \cref{prop:OFSfromStableSystem}.
	Vertical arrows in the left class of the resulting orthogonal factorisation system
	are characterised as those arrows mapped to an isomorphism by the reflection 
	of the inclusion $\one{M}\hto\cod^{\oneV{\dbl{D}}}$.
	Therefore, it suffices to show that, for each $f$,
	$f$ is final if and only if
	$m$ in \cref{dgm:TheFactorisationFin} is an isomorphism.

	If $f$ is final, then $m$ is an isomorphism due to $f^*!_!$ being
	isomorphic to $!_!$, and its tabulator must be isomorphic to the identity
	by \cref{lem:idTab}.
	For the converse, suppose $m$ is an isomorphism. Then, $f^*!_!\cong m^*!_!$ is isomorphic to $!_!$, and $f$ is final.
\end{proof}

\begin{example}
	The double category $\Prof$ has strong tabulators, and discrete fibrations are composable. \cref{thm:CompSchemeToFactorisation} gives the classical comprehensive factorisation system in \cite{SW73,Kel82}.

	More generally, applying this theorem to the double category of 
	internal categories $\Prof(\one{E})$ in a category $\one{E}$ with pullbacks and pullback-stable reflexive coequalisers in the sense of \cite{Kou14},
	we obtain the comprehensive factorisation system internal to $\one{E}$, 
	which is a special case of \cite[Theorem 6.3]{TY21}.
	The existence of strong tabulators is shown in \cite[Theorem 5.15]{Kou14},
	which originally goes back to \cite{Bet96}.

	Since any tabulator of a profunctor between groupoids is a groupoid,
	the sub-double category $\Prof(\one{E})_{\mathrm{gpd}}$ of $\Prof(\one{E})$ spanned by 
	the internal groupoids also has strong tabulators.
	This leads to the comprehensive factorisation for internal groupoids, observed in \cite{Bou87}.
\end{example}

\begin{corollary}
	Let $\dbl{D}$ be an equipment with a vertical terminal object and
	left-sided strong tabulators
	such that $\oneV{\dbl{D}}$ has pullbacks.
	Suppose, moreover, that $\Fib(\dbl{D})$ is closed under composition.
	Then, $(\Fin(\dbl{D}),\Fib(\dbl{D}))$ is an
	orthogonal factorisation system on $\oneV{\dbl{D}}$.
\end{corollary}
\begin{proof}
	$\Fib(\dbl{D})$ is a stable system by \cref{lem:FibStable},
	and it follows from \cref{thm:CompSchemeToFactorisation} where
	$\zero{M}$ is defined as $\Fib(\dbl{D})$.
\end{proof}

\begin{lemma}\label{lem:FactorisationOfSpan}
	Let $\dbl{D}$ be a cartesian equipment with Beck-Chevalley pullbacks
	and assume that $\dbl{D}$ admits a left-sided $\Fib(\dbl{D})$-comprehension scheme.
	If we have two opcartesian triangles
	\begin{equation}
			\begin{tikzcd}[column sep = small]
					&
						C 
						\ar[ld,"f"']
						\ar[rd,"g"]
							&
				\\
				A 
				\sar[rr,"p"']
				\doublecell[rr, shift left=2ex]{\opcart}
					&
						&
							B
			\end{tikzcd}
			\hspace{1ex}
			,
			\hspace{1ex}
			\begin{tikzcd}[column sep = small]
					&
					\top p
					\ar[ld,"l"']
					\ar[rd,"r"]
						&
				\\
					A 
					\sar[rr,"p"']
					\doublecell[rr, shift left=3ex]{\tab}
						&
							&
							B
			\end{tikzcd}
			\hspace{1ex}
			,
	\end{equation}
	then the unique canonical vertical arrow $e\colon C\to\top p$ that satisfies 
	$e\fatsemi l=f,\ e\fatsemi r=g$ is final.
\end{lemma}
\begin{proof}
	Let $\bar{p}\colon A\times B\sto 1$ be the horizontal arrow
	that corresponds to $p$ in the way described in
	\cref{rem:DagcpctCorrespondence}.
	By \cref{lem:UnbiasingLemma},
	we obtain the following cells whose composite is opcartesian,
	in which, by the proof of \cref{thm:CompSchemeToFactorisation},
	we observe that $e\fatsemi \langle l,r\rangle$ gives the factorisation of
	$\langle f, g\rangle$, and hence $e$ is final.
	\[
		\begin{tikzcd}[column sep=small]
				&
				C
				\ar[ldd, "{\langle f,g\rangle}"', bend right=20]
				\ar[rdd, "!", bend left=20]
				\ar[d,"e"]
					&
			\\
				&
				\top p
				\ar[ld, "{\langle l,r\rangle}"description]
				\ar[rd, "!"]
					&
			\\
			A\times B
			\sar[rr, "\bar{p}"']
			\doublecell[rr, shift left=2.5ex]{\tab}
				&
					&
					\hspace{2ex}1\hspace{2ex}
		\end{tikzcd}
	\]
\end{proof}
\begin{lemma}\label{lem:Characterisation}
	Let $\dbl{D}$ be a cartesian equipment with Beck-Chevalley pullbacks and
	$\zero{M}$ be a stable system on $\oneV{\dbl{D}}$.
	Suppose that $\dbl{D}$ admits a left-sided $\zero{M}$-comprehension scheme.
	Then,
	there is an equivalence
	$\dbl{D}\simeq\Rel{\Fin(\dbl{D}), \zero{M}}(\oneV{\dbl{D}})$.
\end{lemma}
\begin{proof}
	By \cref{thm:CompSchemeToFactorisation}, $\zero{M}$ is identical to $\Fib(\dbl{D})$.
	We write $\dbl{R}$ for
	$\Rel{\Fin(\dbl{D}),\Fib(\dbl{D})}(\oneV{\dbl{D}})$
	for brevity.
	For any horizontal arrow $R\colon A\sto B$ in $\dbl{R}$,
	we write $\langle l_R, r_R\rangle\colon \abs{R}\to A \times B$
	for the span defining $R$, and we mean by $F(R)$
	the extension $l_R^*{r_R}_!$ of the span.

	By \cref{prop:compschemeunbiased},
	we see that $\dbl{D}$ admits a $\Fib(\dbl{D})$-comprehension scheme,
	not only the left-sided one.
	Therefore, we have the equivalence
	considered in \cref{lem:CompScheme} below.
	\[
		G \colon  \biH{\dbl{D}}(A, B) \to \biH{\dbl{R}}(A, B)
	\]
	Then, $F$ is the pseudo-inverse
	of $G$. Note that they define
	a fibrewise adjoint equivalence between the bifibrations
	defining the equipments $\dbl{D}$ and $\dbl{R}$.

	Recall again from \cite[Theorem 5.3.7]{Vas14} that a fibrewise adjoint
	extends to an adjoint in the slice 2-category
	$\CATbi/\oneV{\dbl{D}}\times\oneV{\dbl{D}}$ if 
	it is stable under reindexing functors. In this case, to show that
	$F$ extends to an adjoint equivalence in $\CATbi/\oneV{\dbl{D}}\times\oneV{\dbl{D}}$,
	it suffices to check that $F$ is stable under extensions;
	that is,
	for each pair of vertical arrows $u\colon  A \to X$ and $v\colon  B\to Y$
	and any $\zero{M}$-relation $\langle l_R, r_R\rangle\colon \abs{R}\to A\times B$,
	the extension of $F(R)$ along $(u, v)$ is naturally isomorphic to
	$F(\langle l', r'\rangle)$,
	where $\langle l', r'\rangle$ is the image of
	$\langle l_R\fatsemi u, r_R\fatsemi v\rangle$
	with respect to
	$(\Fin(\dbl{D}),\Fib(\dbl{D}))$.
	\cref{lem:FactorisationOfSpan} shows that such a factorisation
	$\langle l_R\fatsemi u, r_R\fatsemi v\rangle = e\fatsemi \langle l',r'\rangle$ is obtained 
	by considering the following tabulator.
	\[
		\begin{tikzcd}[column sep=small]
				&
				\abs{R}
				\ar[ld, "l_R"']
				\ar[rd, "r_R"]
				\ar[d, "e"]
					&
			\\
			A
			\ar[d, "u"']
				&
				\cdot
				\ar[ld,"l'"']
				\ar[rd,"r'"]
						&
						B
					\ar[d,"v"]
			\\
			X
			\sar[rr, "u^*F(R)v_!"']
			\doublecell[rr, shift left=2.5ex]{\tab}
				&
						&
						Y
		\end{tikzcd}
		\hspace{2ex}
		=
		\hspace{2ex}
		\begin{tikzcd}[column sep=small]
				&
				\abs{R}
				\ar[ld, "l_R"']
				\ar[rd, "r_R"]
					&
			\\
			A
			\ar[d, "u"']
			\sar[rr, "F(R)"']
			\doublecell[rrd,yshift=-1ex]{\opcart}
			\doublecell[rr, shift left=3ex]{\opcart}
				&
						&
						B
					\ar[d,"v"]
			\\
			X
			\sar[rr, "u^*F(R)v_!"']
				&
						&
						Y
		\end{tikzcd}
	\]
	Since tabulators are strong and $e$ is a cover by \cref{lem:CovEqualFin}, we have the canonical isomorphism
	$u^*F(R)v_!\cong F(\langle l',r'\rangle)$.

	Now that we obtain an adjoint equivalence $G\colon \dbl{D}_1\adjointleft\dbl{R}_1\lon F$
	in $\CATbi/\oneV{\dbl{D}}\times\oneV{\dbl{D}}$,
	it suffices to show that they are compatible with horizontal compositions.

	Recall that, for two horizontal arrows
	$R\colon  A \sto B$ and $S\colon  B \sto C$ in $\dbl{R}$,
	the composite $RS$ is the $\Fib(\dbl{D})$-image of
	the legs of the composite of the spans defining $R$ and $S$.
	Since $\dbl{D}$ has Beck-Chevalley pullbacks,
	we obtain the following opcartesian and cartesian cells in $\dbl{D}$,
	where the square in the middle is a pullback square.
	\begin{equation}
		\begin{tikzcd}
				&
					&
					\cdot
					\ar[ld,""']
					\ar[lldd,"f"', bend right = 40]
					\ar[rd,""]
					\ar[rrdd,"g", bend left = 40]
						&
							&
			\\
				&
				\abs{R}
				\ar[ld,"l_R"']
				\ar[rd,"r_R"description]
				\doublecell[rr, shift left=3ex]{\opcart}
				\doublecell[rr, shift right=3ex]{\cart}
				\sar[rr]
					&
						&
						\abs{S}
						\ar[ld,"l_S"description]
						\ar[rd,"r_S"]
								&
			\\
			A 
			\sar[rr, "F(R)"']
			\doublecell[rr, shift left=2.5ex]{\opcart}
				&
					&
					B 
					\sar[rr,"F(S)"']
					\doublecell[rr, shift left=2.5ex]{\opcart}
						&
							& 
							C
		\end{tikzcd}
	\end{equation}
	Thanks to \cref{lem:Sandwich}, the whole cell is opcartesian.
	The composite $RS$ is
	the $\Fib(\dbl{D})$-image of $\langle f, g\rangle$,
	and by \cref{lem:FactorisationOfSpan},
	it is realised by 
	taking the tabulator of the composite $F(R)F(S)$ in $\dbl{D}$.
	This means that $F(R)F(S)$ and $F(RS)$ are isomorphic to each other.

	On the other hand, for two horizontal arrows
	$p\colon A \sto B$ and $q\colon B\sto C$ 
	in $\dbl{D}$, consider the cells in $\dbl{D}$ presented below, with $l_p$ and $r_p$ denoting $l_{G(p)}$ and $r_{G(p)}$ 
	and similarly for $l_q$ and $r_q$.
	\begin{equation}
		\begin{tikzcd}[column sep = small]
				&
					&
					\cdot
					\ar[ld,""']
					\ar[lldd,"f"', bend right = 40]
					\ar[rd,""]
					\ar[rrdd,"g", bend left = 40]
						&
							&
			\\
				&
				\abs{G(p)}
				\ar[ld,"l_p"']
				\ar[rd,"r_p"description]
				\doublecell[rr, shift left=3ex]{\opcart}
				\doublecell[rr, shift right=3ex]{\cart}
				\sar[rr]
					&
						&
						\abs{G(q)}
						\ar[ld,"l_q"description]
						\ar[rd,"r_q"]
								&
			\\
			A 
			\sar[rr, "p"']
			\doublecell[rr, shift left=2.5ex]{\tab}
				&
					&
					B 
					\sar[rr,"q"']
					\doublecell[rr, shift left=2.5ex]{\tab}
						&
							& 
							C
		\end{tikzcd}
	\end{equation}
	Here, the square in the middle is a pullback square.
	Since tabulators are strong, a similar discussion concludes that
	the composite $G(p)G(q)$ in $\dbl{R}$ is obtained through the factorisation
	of $\langle f, g\rangle$, which is achieved by taking the tabulator
	of the composite $pq$. Therefore, we conclude $G(pq)$ and $G(p)G(q)$
	are naturally isomorphic.
\end{proof}

\begin{remark}
	Even if we drop the assumption that $\dbl{D}$ admits an $\zero{M}$-comprehension scheme,
	then we can still show that there is the oplax/lax adjunction by the same construction.
	This is because we still have the adjunction mentioned above between 
	$\biH{\dbl{D}}(A,B)$ and $\biH{\dbl{R}}(A,B)$ for each pair of objects $A$ and $B$,
	as in \cref{rem:Adjexsist}. 
	This method to construct the oplax/lax adjunction dates back to \cite[\S 5]{Nie12}
	and has further been developed in \cite{Ale18,Lam22}.
\end{remark}

Finally, we have:

\begin{theorem}\label{thm:MainThm}
	The following are equivalent for a double category $\dbl{D}$.
	\begin{enumerate}
		\item
			$\dbl{D}$ is equivalent to $\Rel{\zero{E}, \zero{M}}(\one{C})$
			for some category $\one{C}$ with finite limits and a stable factorisation system
			$(\zero{E}, \zero{M})$ on $\one{C}$.
		\item
			$\dbl{D}$ is a cartesian equipment with
			Beck-Chevalley pullbacks and
			an $\zero{M}$-comprehension scheme for some stable system
			$\zero{M}$ on $\oneV{\dbl{D}}$.
		\item
			$\dbl{D}$ is a cartesian equipment with
			Beck-Chevalley pullbacks and a
			left-sided $\zero{M}$-comprehension scheme for some stable system
			$\zero{M}$ on $\oneV{\dbl{D}}$.
		\item
			$\Fib(\dbl{D})$ is closed under composition, and
			$\dbl{D}$ is a cartesian equipment with
			Beck-Chevalley pullbacks and
			strong tabulators.
		\item
			$\Fib(\dbl{D})$ is closed under composition, and
			$\dbl{D}$ is a cartesian equipment with
			Beck-Chevalley pullbacks and
			left-sided strong tabulators.
	\end{enumerate}
\end{theorem}
By \emph{\acf{DCR}}, we mean 
the double categories with the conditions proven to be equivalent in this theorem from now on.
\begin{proof} 
	$i)\Rightarrow ii)$ is proved in \cref{prop:RelEMisCartBC,prop:RelEMhasCompScheme}.
	$ii)\Rightarrow iii)$ is trivial and $iii)\Rightarrow i)$ follows from \cref{lem:Characterisation}.
	From the same lemma, $iii)\Rightarrow v)$ follows.
	To show $v)\Rightarrow iii)$, we can take $\Fib(\dbl{D})$ as $\zero{M}$.
	$iv)\Rightarrow v)$ is trivial, and $ii)$ and $v)$ imply $iv)$.
\end{proof}

\section{Several classes of double categories of relations}\label{sec:severalclasses}

\subsection{Local properties --- proper factorisation systems}\label{sec:properness}

The first local property we consider is called unit-pureness.
The unit-pure condition for double categories is introduced in \cite{Ale18},
and Lambert develops the theory of double categories of relations based on this condition in \cite{Lam22}.
We shall see how this condition works in the context of double categories of relations in our framework.

\begin{definition}[{\cite[Definition 4.3.7]{Ale18}}]\label{defn:unitpure}
	A double category $\dbl{D}$ is called \emph{unit-pure}
	if $\Id\colon \dbl{D}_0\to\dbl{D}_1$ is fully faithful.
	In more concrete terms,
	$\dbl{D}$ is unit-pure if every cell of the form
	\begin{equation*}\label{eq:unit-pure}
		\begin{tikzcd}[column sep=small]
			X
			\sar[r,"\Id_X"]
			\ar[d, "f"']
			\doublecell[rd]{\cdot}
				&
				X
				\ar[d,"g"]
			\\
			Y
			\sar[r,"\Id_Y"']
				&
				Y
		\end{tikzcd}
	\end{equation*}
	is an identity cell,
	implying $f$ and $g$ are equal. 
\end{definition}

Another way to express the unit-pureness is the following: for every object $X$,
the identity cell 
\begin{equation}\label{eq:unit-pure2}
	\begin{tikzcd}[column sep=small]
			&
			X 
			\ar[ld, "\id_X"']
			\ar[rd, "\id_X"]
				&
		\\
			X
			\sar[rr, "\Id_X"']
				&
				&
				X
	\end{tikzcd}
\end{equation}
exhibits $X$ as a tabulator of $\Id_X$.
What is more, the unit-pureness of $\dbl{D}$ is equivalent to the vertical 2-category $\biV{\dbl{D}}$ being locally discrete.

\begin{lemma}[{\cite[Proposition 5.9]{Kou14}}]\label{lem:comma} 
	Let $\dbl{D}$ be a double category.
	Suppose that
	we are given
	a cartesian cell $\alpha$ and a tabulating cell $\tau$ 
	of the form below.
	\begin{equation}
		\begin{tikzcd}[column sep=small]
				&
				\mmbox{\top C(f, g)}
				\ar[ld, "l"']
				\ar[rd, "r"]
					&
			\\
			A
			\ar[rd, "f"']
			\sar[rr, "\mmbox{C(f, g)}"]
			\doublecell[rr, shift left=3.5ex]{\tau}
			\doublecell[rr, shift right=3ex]{\alpha}
				&
					&
					B
					\ar[ld, "g"]
			\\
				&
				C
					&
		\end{tikzcd}
	\end{equation}
	Then,
	the composite of $\tau$ and $\alpha$ exhibits
	$\top C(f, g)$ as
	a comma object $f\downarrow g$
	in $\biV{\dbl{D}}$.
\end{lemma}
\begin{proof}
	Let $s\colon Z\to A$ and $t\colon Z\to B$ be vertical arrows 
	and $\beta$ be a cell of the form 
	\begin{equation*}
		\begin{tikzcd}[column sep=small, row sep=small]
				&
				Z
				\ar[ld, "s"']
				\ar[rd, "t"]
					&
			\\
			A
			\ar[rd, "f"']
					&
					\beta
					&
					B
					\ar[ld, "g"]
			\\
				&
				C
					&
		\end{tikzcd}.
	\end{equation*}	
	Since $\alpha$ is cartesian, $\beta$ factors uniquely through $\alpha$,
	and $\tau$, in turn, gives a unique vertical arrow $u\colon Z\to\top C(f, g)$.
	It is the only vertical arrow
	that is composed with $\tau$ and $\alpha$ to give $\beta$.
	This is verified by consecutively applying the universal property of $\alpha$ and $\tau$.
	Thus, the composition of $\tau$ and $\alpha$ exhibits $\top C(f, g)$ as a comma object.
\end{proof}

\begin{corollary}\label{cor:comma}
	Let $\dbl{D}$ be an equipment with tabulators.
	Then, its vertical 2-category $\biV{\dbl{D}}$ has comma objects.
	In particular, if $\dbl{D}$ is also unit-pure, then $\oneV{\dbl{D}}$
	has pullbacks.
	Moreover, if $\dbl{D}$ is unit-pure and has strong tabulators,
	$\dbl{D}$ has Beck-Chevalley pullbacks.
\end{corollary}

\begin{lemma}\label{lem:covepi}
	For a unit-pure double category $\dbl{D}$,
	any cover is epic in $\oneV{\dbl{D}}$, and any inclusion is monic in $\oneV{\dbl{D}}$.
\end{lemma}
\begin{proof}
	We show covers are epic, and the other is the vertical dual.
	Let $e\colon A\to B$ be a cover in $\dbl{D}$. To show $e$ is epic, take
	parallel vertical arrows $f, g\colon B\to C$ in $\dbl{D}$ with
	$e\fatsemi f=e\fatsemi g$.
	The identity vertical cell $e\fatsemi f\to[Rightarrow]e\fatsemi g$ factors through the cover
	$e$, whence we obtain a vertical cell
	$\alpha\colon f\to[Rightarrow]g\colon B\to C$ and hence $f=g$ holds because
	$\dbl{D}$ is unit-pure.
\end{proof}

\begin{corollary}[{\cite[Lemma 4.6]{Lam22}}] 
	\label{lem:monoinc}
	In a unit-pure double category $\dbl{D}$ with strong tabulators,
	the class of inclusions coincides with the class of monomorphisms.
\end{corollary}
\begin{proof}
	It follows from \cref{lem:monoincl,lem:covepi} and \cref{cor:comma}.
\end{proof}

Thus, we have the characterisation of unit-pure double categories of relations
described more concisely than \cref{thm:MainThm}.

\begin{theorem}\label{thm:unitpuredbl}
	For a double category $\dbl{D}$,
	the following are equivalent:
	\begin{enumerate}
		\item%
			$\dbl{D}$ is a unit-pure double category of relations.
		\item%
			$\dbl{D}$ is equivalent to $\Rel{\zero{E}, \zero{M}}(\one{C})$
			for some category $\one{C}$ with finite limits and a left-proper stable factorisation system
			$(\zero{E}, \zero{M})$ on $\one{C}$.
		\item%
			$\dbl{D}$ is a unit-pure cartesian equipment with an
			$\zero{M}$-comprehension scheme for 
			some stable system $\zero{M}$ on $\oneV{\dbl{D}}$.
		\item%
			$\Fib(\dbl{D})$ is closed under composition,
			and $\dbl{D}$ is a unit-pure cartesian equipment with strong tabulators. 
		\item%
			$\dbl{D}$ is a unit-pure discrete cartesian equipment with a
			left-sided $\zero{M}$-comprehension scheme for 
			some stable system $\zero{M}$ on $\oneV{\dbl{D}}$.
	\end{enumerate}
	In particular, if $ii)$, $iii)$, or $v)$ holds,
	then $\zero{M}$ is the same class as $\Fib(\dbl{D})$.
\end{theorem}

\begin{proof}
	For the equivalence of $i)$ and $ii)$, it suffices to show that,
	$\Rel{\zero{E}, \zero{M}}(\one{C})$ is unit-pure if and only if 
	$(\zero{E}, \zero{M})$ is left-proper.
	The left class of the double category of relations is the class of final morphisms,
	hence the class of covers by \cref{lem:CovEqualFin}.
	If $\Rel{\zero{E}, \zero{M}}(\one{C})$ is unit-pure,
	then by \cref{lem:covepi}, every arrow in $\zero{E}$ is epic.
	Thus, $(\zero{E}, \zero{M})$ is left-proper.
	Conversely, if $(\zero{E}, \zero{M})$ is left-proper,
	then every diagonal arrow $\Delta_X\colon X\to X\times X$ is in $\zero{M}$
	since an arrow is epimorphic precisely when it is left orthogonal to every diagonal arrow. 
	By the $\zero{M}$-comprehension scheme,
	every cell of the form \cref{eq:unit-pure2} exhibits $X$ as a tabulator of $\Delta_X$,
	which shows that $\Rel{\zero{E}, \zero{M}}(\one{C})$ is unit-pure.

	By \cref{cor:comma},
	the existence of strong tabulators and the unit-pureness of $\dbl{D}$ leads to 
	the existence of Beck-Chevalley pullbacks.
	Therefore, the conditions $i)$, $iii)$ and $iv)$ are equivalent by \cref{thm:MainThm}.
	In particular, this implies that $\dbl{D}$ is discrete if $\dbl{D}$ is unit-pure and has strong tabulators, and hence $iii)$ implies $v)$.
	Conversely, $v)$ implies $iii)$ by \cref{prop:compschemeunbiased}.
\end{proof}

We proceed to the characterisation of locally posetal double categories of relations.

\begin{definition}
	\label{defn:locallyPreordered}
	Let $\dbl{D}$ be a double category.
	We say that $\dbl{D}$ is \emph{locally preordered}
	if there must be at most one cell framed by a pair of vertical arrows and a pair of horizontal arrows.
	We say that $\dbl{D}$ is \emph{locally posetal}
	if it is locally preordered and
	the vertical 2-category $\biV{\dbl{D}}$ is locally posetal.
\end{definition}

In some papers \cite{GP99,Ste24},
a locally preordered double category is called a \emph{flat double category}.

\begin{remark}
	There are two things to note here.

	Firstly, observe that
	given a cell in $\dbl{D}$, we can take the restriction of the bottom horizontal arrow and consider it as a 2-cell in $\biH{\dbl{D}}$,
	and hence
	an equipment $\dbl{D}$ is locally preordered if and only if
	the horizontal bicategory $\biH{\dbl{D}}$ is locally preordered,

	Secondly, the definition of local posetality might seem inappropriate for its name at first,
	but
	this is grounded on the fact
	that an equipment $\dbl{D}$ is locally posetal if and only if 
	it is equivalent in $\DblCat$ to one with stricter conditions:
	not only the vertical 2-category $\biV{\dbl{D}}$ is locally posetal, but 
	the horizontal bicategory $\biH{\dbl{D}}$ is locally posetal.
	This is because equivalences in $\DblCat$ loosens how skeletal the horizontal bicategory of a double category is,
	whereas it leaves 2-cells of its vertical 2-category essentially unchanged.
\end{remark}

\begin{lemma}
	\label{lem:DagLocPos}
	A discrete cartesian equipment $\dbl{D}$
	is locally posetal
	if and only if it is locally preordered and unit-pure.
\end{lemma}

\begin{proof}
	Suppose $\dbl{D}$ is locally preordered.
	Unit-pureness means that the vertical 2-category $\biV{\dbl{D}}$ 
	is locally discrete.
	Hence, $\dbl{D}$ is locally posetal if it is unit-pure.

	From \cref{prop:pseudoFunctorialityDagger},
	it follows that $\biV{\dbl{D}}$ is 2-equivalent to its 2-cell dual $\biV{\dbl{D}}^\co$ through a 2-functor that is
	the identity on underlying categories, and hence 2-cells in $\biV{\dbl{D}}$ are all invertible.
	Therefore, $\dbl{D}$ is locally posetal if and only if $\biV{\dbl{D}}$ is locally discrete,
	or equivalently, $\dbl{D}$ is unit-pure.
\end{proof}

\begin{theorem}\label{thm:locorddbl}
	A double category of relations is locally preordered 
	if and only if the accompanying factorisation system is right-proper.
\end{theorem}

\begin{proof}
	Let $\dbl{D}$ be a double category of relations
	and $(\zero{E}, \zero{M})$ be the accompanying factorisation system.
	We have $\biH{\dbl{D}}(X,Y)\simeq\one{M}\downarrow X\times Y$
	for any objects $X$ and $Y$ in $\dbl{D}$.
	If the factorisation system is right-proper,
	$\one{M}\downarrow X\times Y$ is a preordered set, making $\biH{\dbl{D}}$ locally preordered.
	Conversely, if $\biH{\dbl{D}}$ is locally preordered,
	then $\one{M}\downarrow X$ is a preorder for every object $X$ in $\dbl{D}$.
	In general, binary products are idempotent in a preordered set, which, for this preordered set, means that any object is monic as an arrow in $\oneV{\dbl{D}}$,
	and this concludes $\zero{M}\subset\Mono$.
\end{proof}

\begin{remark}
	We briefly mention the work of \cite{JW00,HNST22} in relation to locally preordered double categories of relations.
	Right-proper factorisation systems are the central subject of \cite{JW00},
	in which the authors deal with functional relations, or \textit{maps} that we are now to discuss in \cref{subsec:cauchyreg}.

	In \cite{HNST22}, it is proved that the bicategories of $\zero{M}$-relations are allegories if and only if the factorisation system is right-proper.
	This result also follows from \cref{thm:locorddbl,rem:involution}.
	They further show that such bicategories are unital tabular allegories.
	From our perspective, this is a consequence of the fact that the locally-preordered double categories of relations have a terminal object,
	have strong tabulators, and admit a $\zero{M}$-comprehension scheme where $\zero{M}$ is included in the class of inclusions.
\end{remark}

\begin{corollary}\label{cor:locpos}
	A double category of relations is locally posetal 
	if and only if the accompanying factorisation system is proper.
\end{corollary}
\begin{proof} 
	The statement follows from \Cref{thm:unitpuredbl,lem:DagLocPos,thm:locorddbl}.
\end{proof}

\subsection{Cauchy double categories of relations}
\label{subsec:cauchyreg}

The relations between two objects in the usual sense, or $\Mono$-relations in our terminology,
are typical examples of the notion of relations.
However, we need to require the category to be regular if we want to have a double category of $\Mono$-relations
as mentioned in \cref{ex:RelRegCat}.
The aim of this subsection is to determine the conditions crucial for these particular double categories of relations.
To do so, we use the notion of Cauchy double categories 
which was introduced by Par\'e .

\begin{definition}[{\cite[Definition 19]{Par21}}]
	A double category $\dbl{D}$ is \emph{Cauchy} if
	any adjoint $p\colon X\adjointleft Y\lon q$ in the bicategory $\bi{H}(\dbl{D})$
	is representable, 
	namely, of the form $f_!\colon X\adjointleft Y\lon f^*$ for some $f\colon X\to Y$.
\end{definition}

The name comes from an observation in \cite[Theorem 2]{BD86} that
for a small category $\one{C}$, 
it is Cauchy complete if and only if every profunctor with a right adjoint profunctor from every small category to $\one{C}$ 
is a representable profunctor.
In terms of double categories,
this means that Cauchy complete categories are precisely the ones such that every horizontal arrow to it is representable in 
the double category $\Prof$.
If we define Cauchy completeness for an object in a double category by the same condition,
then Cauchy double categories are double categories in which every object is Cauchy complete.

The critical property of Cauchy double categories emerges in the following lemma.

\begin{lemma}
	\label{lem:covincortho}
	Let $\dbl{D}$ be a unit-pure Cauchy equipment.
	Then, every cover is left orthogonal to every inclusion in $\dbl{D}$.
\end{lemma}

\begin{proof}
	Let $m\colon A\to B$ be an inclusion and $e\colon C\to D$ be a cover and suppose we are given
	the following diagram in $\oneV{\dbl{D}}$.
	\[
		\begin{tikzcd}
				C
					\ar[d,"e"']
					\ar[r,"g"]
			&
				A
					\ar[d,"m"]
			\\
				D
					\ar[r,"f"']
			&
				B
		\end{tikzcd}
	\]
	The first step is to show that  $e^*g_!$ is left adjoint to $m_!f^*$.
	We employ string diagrams for equipments which were introduced in \cite{Mye18}.
	The two horizontal identity cells corresponding to this commutative square are displayed as follows:
	
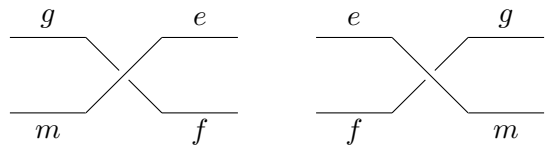
\[
		\begin{tikzpicture}
	\begin{pgfonlayer}{nodelayer}
		\node [style=none] (0) at (-3, 1) {};
		\node [style=none] (1) at (-1, 1) {};
		\node [style=none] (2) at (-3, -1) {};
		\node [style=none] (3) at (-1, -1) {};
		\node [style=none] (4) at (1, 1) {};
		\node [style=none] (5) at (3, 1) {};
		\node [style=none] (6) at (1, -1) {};
		\node [style=none] (7) at (3, -1) {};
		\node [style=none] (8) at (2, 1.5) {$e$};
		\node [style=none] (9) at (2, -1.5) {$f$};
		\node [style=none] (10) at (-2, 1.5) {$g$};
		\node [style=none] (11) at (-2, -1.5) {$m$};
	\end{pgfonlayer}
	\begin{pgfonlayer}{edgelayer}
		\draw (1.center) to (0.center);
		\draw (2.center) to (3.center);
		\draw (5.center) to (4.center);
		\draw (7.center) to (6.center);
		\draw (1.center) to (6.center);
		\draw [style=white] (4.center) to (3.center);
		\draw (4.center) to (3.center);
	\end{pgfonlayer}
\end{tikzpicture}
		\hspace{6ex}
		\begin{tikzpicture}
	\begin{pgfonlayer}{nodelayer}
		\node [style=none] (0) at (-3, 1) {};
		\node [style=none] (1) at (-1, 1) {};
		\node [style=none] (2) at (-3, -1) {};
		\node [style=none] (3) at (-1, -1) {};
		\node [style=none] (4) at (1, 1) {};
		\node [style=none] (5) at (3, 1) {};
		\node [style=none] (6) at (1, -1) {};
		\node [style=none] (7) at (3, -1) {};
		\node [style=none] (8) at (2, 1.5) {$g$};
		\node [style=none] (9) at (2, -1.5) {$m$};
		\node [style=none] (10) at (-2, 1.5) {$e$};
		\node [style=none] (11) at (-2, -1.5) {$f$};
	\end{pgfonlayer}
	\begin{pgfonlayer}{edgelayer}
		\draw (1.center) to (0.center);
		\draw (2.center) to (3.center);
		\draw (5.center) to (4.center);
		\draw (7.center) to (6.center);
		\draw (4.center) to (3.center);
		\draw [style=white] (1.center) to (6.center);
		\draw (1.center) to (6.center);
	\end{pgfonlayer}
\end{tikzpicture}
	\]
	Here we use the convention of crossing one string over another to indicate the direction 
	of equality of vertical arrows.
	The unit and the counit of $e^*g_!\dashv m_!f^*$
	are given as below:
	\begin{equation}
		\begin{tikzpicture}
	\begin{pgfonlayer}{nodelayer}
		\node [style=none] (1) at (-2, 0.5) {};
		\node [style=none] (2) at (-2, -2) {};
		\node [style=none] (3) at (2, 0.5) {};
		\node [style=none] (4) at (2, 0) {};
		\node [style=none] (5) at (0.5, 0) {};
		\node [style=none] (6) at (1.5, -1) {};
		\node [style=none] (7) at (-0.75, 0) {};
		\node [style=none] (8) at (-0.75, -2) {};
		\node [style=none] (9) at (1.5, 0) {};
		\node [style=none] (10) at (0.5, -1) {};
		\node [style=none] (13) at (2, -2) {};
		\node [style=none] (14) at (2, -1) {};
		\node [style=none] (15) at (0, -1) {};
		\node [style=none] (16) at (0, -2) {};
		\node [style=none] (17) at (-2.5, -0.75) {$e^*$};
		\node [style=none] (18) at (0.5, -1.75) {$m_!$};
		\node [style=huge cell] (19) at (0, 1) {$\varepsilon^{-1}$};
		\node [style=none] (20) at (-1.25, -1.25) {$g_!$};
		\node [style=none] (21) at (2.5, -1.5) {$f^*$};
	\end{pgfonlayer}
	\begin{pgfonlayer}{edgelayer}
		\draw (1.center) to (2.center);
		\draw (3.center) to (4.center);
		\draw (5.center) to (7.center);
		\draw (7.center) to (8.center);
		\draw (4.center) to (9.center);
		\draw (14.center) to (6.center);
		\draw (10.center) to (15.center);
		\draw (15.center) to (16.center);
		\draw (14.center) to (13.center);
		\draw (5.center) to (6.center);
		\draw [style=white] (9.center) to (10.center);
		\draw (9.center) to (10.center);
	\end{pgfonlayer}
\end{tikzpicture}
		\hspace{6ex}
		\begin{tikzpicture}
	\begin{pgfonlayer}{nodelayer}
		\node [style=none] (1) at (-2, -1) {};
		\node [style=none] (2) at (-2, 1.5) {};
		\node [style=none] (3) at (2, -1) {};
		\node [style=none] (4) at (2, -0.5) {};
		\node [style=none] (5) at (0.5, -0.5) {};
		\node [style=none] (6) at (1.5, 0.5) {};
		\node [style=none] (7) at (-0.75, -0.5) {};
		\node [style=none] (8) at (-0.75, 1.5) {};
		\node [style=none] (9) at (1.5, -0.5) {};
		\node [style=none] (10) at (0.5, 0.5) {};
		\node [style=none] (11) at (2, 1.5) {};
		\node [style=none] (12) at (2, 0.5) {};
		\node [style=none] (13) at (0, 0.5) {};
		\node [style=none] (14) at (0, 1.5) {};
		\node [style=huge cell] (15) at (0, -1.5) {$\eta^{-1}$};
		\node [style=none] (24) at (0.5, 1) {$e^*$};
		\node [style=none] (25) at (-2.5, 0) {$m_!$};
		\node [style=none] (26) at (2.5, 1) {$g_!$};
		\node [style=none] (27) at (-1.25, 0.5) {$f^*$};
	\end{pgfonlayer}
	\begin{pgfonlayer}{edgelayer}
		\draw [in=270, out=90] (1.center) to (2.center);
		\draw (3.center) to (4.center);
		\draw (5.center) to (7.center);
		\draw (7.center) to (8.center);
		\draw (4.center) to (9.center);
		\draw (12.center) to (6.center);
		\draw (10.center) to (13.center);
		\draw (13.center) to (14.center);
		\draw (12.center) to (11.center);
		\draw (5.center) to (6.center);
		\draw [style=white] (9.center) to (10.center);
		\draw (9.center) to (10.center);
	\end{pgfonlayer}
\end{tikzpicture}
	\end{equation}
	Here, $\eta^{-1}$ and $\varepsilon^{-1}$ are the inverses of the unit and the counit
	for the representable adjunctions for $m$ and $e$, respectively.
	They exist because $m$ is an inclusion and $e$ is a cover, respectively.
	By definition, 
	the following diagrams all amount to vacancy or just two strings
	as string diagrams:
	\[
	\begin{tikzpicture}
	\begin{pgfonlayer}{nodelayer}
		\node [style=big cell] (0) at (0, -1) {$\eta^{-1}$};
		\node [style=none] (1) at (-1, -0.5) {};
		\node [style=none] (2) at (1, -0.5) {};
		\node [style=none] (3) at (-1, 1) {};
		\node [style=none] (4) at (1, 1) {};
	\end{pgfonlayer}
	\begin{pgfonlayer}{edgelayer}
		\draw (3.center) to (1.center);
		\draw (4.center) to (2.center);
		\draw (3.center) to (4.center);
	\end{pgfonlayer}
\end{tikzpicture}
	\hspace{1ex}
	,
	\hspace{1ex}
	\begin{tikzpicture}
	\begin{pgfonlayer}{nodelayer}
		\node [style=big cell] (0) at (0, 0) {$\eta^{-1}$};
		\node [style=none] (1) at (-1, -1.75) {};
		\node [style=none] (2) at (1, -1.75) {};
		\node [style=none] (3) at (-1, -1) {};
		\node [style=none] (4) at (1, -1) {};
		\node [style=none] (5) at (-1, 0.5) {};
		\node [style=none] (6) at (-1, 1.75) {};
		\node [style=none] (7) at (1, 1.75) {};
		\node [style=none] (8) at (1, 0.5) {};
	\end{pgfonlayer}
	\begin{pgfonlayer}{edgelayer}
		\draw (3.center) to (1.center);
		\draw (4.center) to (2.center);
		\draw (3.center) to (4.center);
		\draw (6.center) to (5.center);
		\draw (8.center) to (7.center);
	\end{pgfonlayer}
\end{tikzpicture}
	\hspace{1ex}
	,
	\hspace{1ex}
	\begin{tikzpicture}
	\begin{pgfonlayer}{nodelayer}
		\node [style=big cell] (0) at (0, 0.5) {$\varepsilon^{-1}$};
		\node [style=none] (1) at (-1, 0) {};
		\node [style=none] (2) at (1, 0) {};
		\node [style=none] (3) at (-1, -1.5) {};
		\node [style=none] (4) at (1, -1.5) {};
	\end{pgfonlayer}
	\begin{pgfonlayer}{edgelayer}
		\draw (3.center) to (1.center);
		\draw (4.center) to (2.center);
		\draw (3.center) to (4.center);
	\end{pgfonlayer}
\end{tikzpicture}
	\hspace{1ex}
	,
	\hspace{1ex}
	\begin{tikzpicture}
	\begin{pgfonlayer}{nodelayer}
		\node [style=big cell] (0) at (0, -0.25) {$\varepsilon^{-1}$};
		\node [style=none] (1) at (-1, 1.75) {};
		\node [style=none] (2) at (1, 1.75) {};
		\node [style=none] (3) at (-1, 0.75) {};
		\node [style=none] (4) at (1, 0.75) {};
		\node [style=none] (5) at (-1, -0.75) {};
		\node [style=none] (6) at (-1, -1.75) {};
		\node [style=none] (7) at (1, -0.75) {};
		\node [style=none] (8) at (1, -1.75) {};
	\end{pgfonlayer}
	\begin{pgfonlayer}{edgelayer}
		\draw (3.center) to (1.center);
		\draw (4.center) to (2.center);
		\draw (3.center) to (4.center);
		\draw (5.center) to (6.center);
		\draw (7.center) to (8.center);
	\end{pgfonlayer}
\end{tikzpicture}
	\]

	The diagrammatic argument allows us to verify the triangle identity as follows:
	\begin{equation}
		\begin{tikzpicture}
	\begin{pgfonlayer}{nodelayer}
		\node [style=big cell] (0) at (-1.5, 2.25) {$\varepsilon^{-1}$};
		\node [style=none] (1) at (-3, 1.75) {};
		\node [style=none] (2) at (-3, -2.75) {};
		\node [style=none] (3) at (0, 1.75) {};
		\node [style=none] (4) at (0, 1.25) {};
		\node [style=none] (5) at (-1.5, 1.25) {};
		\node [style=none] (6) at (-0.5, 0.25) {};
		\node [style=none] (7) at (-2.5, 1.25) {};
		\node [style=none] (8) at (-2.5, -2.75) {};
		\node [style=none] (9) at (-0.5, 1.25) {};
		\node [style=none] (10) at (-1.5, 0.25) {};
		\node [style=none] (12) at (0, 0.25) {};
		\node [style=none] (13) at (-2, 0.25) {};
		\node [style=none] (14) at (-2, -1.25) {};
		\node [style=none] (18) at (3, -1.25) {};
		\node [style=none] (19) at (3, -0.75) {};
		\node [style=none] (20) at (1.5, -0.75) {};
		\node [style=none] (21) at (2.5, 0.25) {};
		\node [style=none] (22) at (0, -0.75) {};
		\node [style=none] (24) at (2.5, -0.75) {};
		\node [style=none] (25) at (1.5, 0.25) {};
		\node [style=none] (26) at (3, 3.25) {};
		\node [style=none] (27) at (3, 0.25) {};
		\node [style=none] (28) at (1, 0.25) {};
		\node [style=none] (29) at (1, 3.25) {};
		\node [style=huge cell] (30) at (0.75, -1.75) {$\eta^{-1}$};
	\end{pgfonlayer}
	\begin{pgfonlayer}{edgelayer}
		\draw (1.center) to (2.center);
		\draw (3.center) to (4.center);
		\draw (5.center) to (7.center);
		\draw (7.center) to (8.center);
		\draw (4.center) to (9.center);
		\draw (12.center) to (6.center);
		\draw (10.center) to (13.center);
		\draw (13.center) to (14.center);
		\draw (5.center) to (6.center);
		\draw [style=white] (9.center) to (10.center);
		\draw (9.center) to (10.center);
		\draw (18.center) to (19.center);
		\draw (20.center) to (22.center);
		\draw (19.center) to (24.center);
		\draw (27.center) to (21.center);
		\draw (25.center) to (28.center);
		\draw (28.center) to (29.center);
		\draw (27.center) to (26.center);
		\draw (20.center) to (21.center);
		\draw [style=white] (24.center) to (25.center);
		\draw (24.center) to (25.center);
		\draw (12.center) to (22.center);
	\end{pgfonlayer}
\end{tikzpicture}
		\hspace{1ex}
		=   
		\hspace{1ex}    
		\begin{tikzpicture}
	\begin{pgfonlayer}{nodelayer}
		\node [style=big cell] (0) at (-1.5, 1) {$\varepsilon^{-1}$};
		\node [style=none] (1) at (-3, 0.5) {};
		\node [style=none] (2) at (-3, -3) {};
		\node [style=none] (3) at (0, 0.5) {};
		\node [style=none] (4) at (0, 0) {};
		\node [style=none] (5) at (-1.5, 0) {};
		\node [style=none] (6) at (-0.5, -1) {};
		\node [style=none] (7) at (-2.5, 0) {};
		\node [style=none] (8) at (-2.5, -3) {};
		\node [style=none] (9) at (-0.5, 0) {};
		\node [style=none] (10) at (-1.5, -1) {};
		\node [style=none] (12) at (-2, -1) {};
		\node [style=none] (13) at (-2, -1.5) {};
		\node [style=none] (14) at (3, -1.5) {};
		\node [style=none] (15) at (3, -1) {};
		\node [style=none] (16) at (1.5, -1) {};
		\node [style=none] (17) at (2.5, 0) {};
		\node [style=none] (19) at (2.5, -1) {};
		\node [style=none] (20) at (1.5, 0) {};
		\node [style=none] (21) at (3, 4) {};
		\node [style=none] (22) at (3, 0) {};
		\node [style=none] (23) at (1, 0) {};
		\node [style=none] (24) at (1, 3) {};
		\node [style=huge cell] (25) at (0.75, -2) {$\eta^{-1}$};
		\node [style=none] (26) at (0, 3) {};
		\node [style=none] (27) at (0, 2.25) {};
		\node [style=none] (28) at (-3, 2.25) {};
		\node [style=none] (29) at (-3, 4) {};
	\end{pgfonlayer}
	\begin{pgfonlayer}{edgelayer}
		\draw (1.center) to (2.center);
		\draw (3.center) to (4.center);
		\draw (5.center) to (7.center);
		\draw (7.center) to (8.center);
		\draw (4.center) to (9.center);
		\draw (10.center) to (12.center);
		\draw (12.center) to (13.center);
		\draw (5.center) to (6.center);
		\draw [style=white] (9.center) to (10.center);
		\draw (9.center) to (10.center);
		\draw (14.center) to (15.center);
		\draw (15.center) to (19.center);
		\draw (22.center) to (17.center);
		\draw (20.center) to (23.center);
		\draw (23.center) to (24.center);
		\draw (22.center) to (21.center);
		\draw (16.center) to (17.center);
		\draw [style=white] (19.center) to (20.center);
		\draw (19.center) to (20.center);
		\draw (6.center) to (16.center);
		\draw (24.center) to (26.center);
		\draw (26.center) to (27.center);
		\draw (28.center) to (27.center);
		\draw (29.center) to (28.center);
	\end{pgfonlayer}
\end{tikzpicture}
		\hspace{1ex}
		=
		\hspace{1ex}
		\begin{tikzpicture}
	\begin{pgfonlayer}{nodelayer}
		\node [style=none] (2) at (-3, -3) {};
		\node [style=none] (4) at (0, 0) {};
		\node [style=none] (5) at (-1.5, 0) {};
		\node [style=none] (6) at (-0.5, -1) {};
		\node [style=none] (7) at (-2.5, 0) {};
		\node [style=none] (8) at (-2.5, -3) {};
		\node [style=none] (9) at (-0.5, 0) {};
		\node [style=none] (10) at (-1.5, -1) {};
		\node [style=none] (11) at (-2, -1) {};
		\node [style=none] (12) at (-2, -1.5) {};
		\node [style=none] (13) at (3, -1.5) {};
		\node [style=none] (14) at (3, -1) {};
		\node [style=none] (15) at (1.5, -1) {};
		\node [style=none] (16) at (2.5, 0) {};
		\node [style=none] (17) at (2.5, -1) {};
		\node [style=none] (18) at (1.5, 0) {};
		\node [style=none] (19) at (3, 4) {};
		\node [style=none] (20) at (3, 0) {};
		\node [style=none] (21) at (1, 0) {};
		\node [style=none] (22) at (1, 3) {};
		\node [style=huge cell] (23) at (0.75, -2) {$\eta^{-1}$};
		\node [style=none] (24) at (0, 3) {};
		\node [style=none] (25) at (0, 0) {};
		\node [style=none] (26) at (-3, -3) {};
		\node [style=none] (27) at (-3, 4) {};
	\end{pgfonlayer}
	\begin{pgfonlayer}{edgelayer}
		\draw (5.center) to (7.center);
		\draw (7.center) to (8.center);
		\draw (4.center) to (9.center);
		\draw (10.center) to (11.center);
		\draw (11.center) to (12.center);
		\draw (5.center) to (6.center);
		\draw [style=white] (9.center) to (10.center);
		\draw (9.center) to (10.center);
		\draw (13.center) to (14.center);
		\draw (14.center) to (17.center);
		\draw (20.center) to (16.center);
		\draw (18.center) to (21.center);
		\draw (21.center) to (22.center);
		\draw (20.center) to (19.center);
		\draw (15.center) to (16.center);
		\draw [style=white] (17.center) to (18.center);
		\draw (17.center) to (18.center);
		\draw (6.center) to (15.center);
		\draw (22.center) to (24.center);
		\draw (24.center) to (25.center);
		\draw (27.center) to (26.center);
	\end{pgfonlayer}
\end{tikzpicture}
		\hspace{1ex}
		=
		\hspace{1ex}
		\begin{tikzpicture}
	\begin{pgfonlayer}{nodelayer}
		\node [style=none] (1) at (-0.5, 3.5) {};
		\node [style=none] (2) at (-0.5, -3.5) {};
		\node [style=none] (7) at (0.75, 3.5) {};
		\node [style=none] (8) at (0.75, -3.5) {};
		\node [style=none] (26) at (-0.5, 1.25) {};
		\node [style=none] (27) at (-0.5, 3.5) {};
	\end{pgfonlayer}
	\begin{pgfonlayer}{edgelayer}
		\draw (1.center) to (2.center);
		\draw (7.center) to (8.center);
	\end{pgfonlayer}
\end{tikzpicture}
	\end{equation}
	The other triangle identity is verified dually.
	Therefore, we have $e^*g_!\dashv m_!f^*$.

	Since $\dbl{D}$ is Cauchy, $e^*g_!$ is representable, i.e., of the form 
	$h_!\colon D\sto A$ for some $h\colon D\to A$.
	Then, $h$ is a filler of the square by the following argument.
	The isomorphism $h_!\cong e^*g_!$ gives a cell
	whose vertical arrows are $e\fatsemi h$ and $g$
	and whose horizontal arrows are both identities,
	which leads to the equality $e\fatsemi h=g$ by the unit-pure property.
	The same argument shows that $h\fatsemi m=f$.
	The uniqueness follows since $m$ is monic and $e$ is epic by \Cref{lem:covepi}.
\end{proof}

The Cauchyness of double categories of relations gains significance in the context of double categories of relations in its capacity to express the unique choice principle,
where a unit and a counit of an adjunction $p \colon  A \adjointleft B \lon q$ are the double-categorical counterparts of the existence and the uniqueness 
of an element of $B$ relating to each element of $A$.
This point of view was taken by Rosolini in \cite{Ros99}, for instance.
The following theorem measures the extent to which a double category of relations is capable of `unique choice'.

\begin{theorem}
	\label{thm:mapunitpure}
	Let $\dbl{D}$ be a unit-pure double category of relations.
	If $x_1\colon X\to A$ and $x_2\colon X\to B$ give the tabulator of a horizontal arrow $P\colon A \sto B$ with a right adjoint,
	then $x_1\colon X\to A$ is a cover and a monomorphism simultaneously.
\end{theorem}

In the case of double categories of spans,
$x_1$ becomes an isomorphism. The result for this case is proven in \cite[Proposition 2.2]{CKS84}
by constructing the inverse of $x_1$ from the triangle identity of this adjunction.
In the case of double categories of relations for proper factorisation systems,
the same result as ours is proven in \cite[Theorem 3.3]{Kel91}. 
In the proof, Kelly shows that $x_1$ is precomposed by another arrow and becomes a cover
because of the existence of the unit of adjunction.
Then, he uses the fact that every cover is an epimorphism to conclude that $x_1$ is a cover.
Neither of these arguments works in our case,
but we amalgamate the two arguments to show that $x_1$ is a cover,
and then, we show that $x_1$ is a mono in the same way as Kelly's argument.

\begin{proof}[Proof of \Cref{thm:mapunitpure}]
	Let $P\colon A \adjointleft B \lon  Q $ be an adjoint pair in $\bi{H}(\dbl{D})$.
	Take the tabulators of $P$, $Q$, $PQ$, $QP$ and $PQP$ as below. 
	\begin{equation*}
		\begin{tikzcd}[column sep = small]
			& 
			X 
				\ar[ld,"x_1"']
				\ar[rd,"x_2"]	
			&
			\\
			A 
				\sar[rr,"P"']
				\doublecell[rr, shift left=3ex]{\tab}	
			&
			&
			B 
		\end{tikzcd},
		\hspace{0.5ex}
		\begin{tikzcd}[column sep = small]
			& 
			Y 
				\ar[ld,"y_1"']
				\ar[rd,"y_2"]	
			&
			\\
			B 
				\sar[rr,"Q"']
				\doublecell[rr, shift left=3ex]{\tab}	
			&
			&
			A 
		\end{tikzcd},
		\hspace{1ex}
		\begin{tikzcd}[column sep = small]
			& 
			C
				\ar[ld,"c_1"']
				\ar[rd,"c_2"]	
			&
			\\
			A 
				\sar[r,"P"']
				\doublecell[rr, shift left=3ex]{\tab}
			&
			B 
				\sar[r,"Q"']
			&
			A 
		\end{tikzcd},
		\hspace{0.5ex}
		\begin{tikzcd}[column sep = small]
			& 
			D
				\ar[ld,"d_1"']
				\ar[rd,"d_2"]	
			&
			\\
			B
				\sar[r,"Q"']
				\doublecell[rr, shift left=3ex]{\tab}
			&
			A 
				\sar[r,"P"']
			&	
			B	
		\end{tikzcd},
		\hspace{1ex}
		\begin{tikzcd}[column sep = tiny]
			&&&
			G
				\ar[llld,"g_1"']
				\ar[rrrd,"g_2"]
			&&&
			\\
			A 
				\sar[rr,"P"']
				\doublecell[rrrrrr, shift left=3ex]{\tab}
			&&
			B 
				\sar[rr,"Q"']
			&&
			A 
				\sar[rr,"P"']
			&&
			B	
		\end{tikzcd}
	\end{equation*}
	The aim for now is to show that $x_1$ is a cover.
	The unit of the adjunction corresponds to a vertical arrow $\eta\colon  A \to C$ such that 
	$\eta\fatsemi c_1=\eta\fatsemi c_2=\id_A$,
	and the counit of the adjunction ensures that $d_1$ and $d_2$ are the same, so 
	we let $\varepsilon\colon D\to B $ denote this arrow.

	We form the following diagrams using pullbacks and the universality of tabulators.
	\begin{equation*}
		\begin{tikzcd}[column sep = small]
			&&&
			Z 
			\ar[ld,"z_1"'] 
			\ar[rd,"z_2"]
			\ar[lllddd,"z_3"',bend right=35]
			\ar[rrrddd,"z_4", bend left=35]
			\ar[dd,phantom,"\mbox{\rotatebox[origin=c]{45}{$\lrcorner$}}"{marking,near start}]
			&&&
			\\
			&& 
			U
			\ar[ld,"u_1"']
			\ar[rd,"u_2"]
			\ar[dd,phantom,"\mbox{\rotatebox[origin=c]{45}{$\lrcorner$}}"{marking,near start}]
			&&
			V
			\ar[ld,"v_1"']
			\ar[rd,"v_2"]
			\ar[dd,phantom,"\mbox{\rotatebox[origin=c]{45}{$\lrcorner$}}"{marking,near start}]
			&&
			\\ 
			& 
			X 
			\ar[ld,"x_1"']
			\ar[rd,"x_2"]
			&&
			Y
			\ar[ld,"y_1"']
			\ar[rd,"y_2"]
			&&
			X 
			\ar[ld,"x_1"']
			\ar[rd,"x_2"]
			\\
			A
			\sar[rr,"P"']
			\doublecell[rr, shift left=3ex]{\tab}
			&&
			B
			\sar[rr,"Q"']
			\doublecell[rr, shift left=3ex]{\tab}
			&&
			A
			\sar[rr,"P"']
			\doublecell[rr, shift left=3ex]{\tab}			
			&&
			B
		\end{tikzcd},
		\hspace{1ex}
		\begin{tikzcd}
			& 
			U
			\ar[ld, "u_1"']
			\ar[d,"\varphi", two heads]
			\ar[rd,"u_2"]
			\\
			X
			\ar[d,"x_1"']
			\ar[r,phantom,"="]
			& 
			C 
			\ar[ld,"c_1"']
			\ar[rd,"c_2"]
			\ar[r,phantom,"="]
			&
			Y
			\ar[d,"y_2"]
			\\
			A
			\sar[r,"P"']
			\doublecell[rr, shift left=3ex]{\tab}
			&
			B
			\sar[r,"Q"']
			& 
			A
		\end{tikzcd},
		\hspace{1ex}
		\begin{tikzcd}
			& 
			V
			\ar[ld, "v_1"']
			\ar[d,"\theta", two heads]
			\ar[rd,"v_2"]
			\\
			Y
			\ar[d,"y_1"']
			\ar[r,phantom,"="]
			& 
			D
			\ar[ld,"\varepsilon"']
			\ar[rd,"\varepsilon"]
			\ar[r,phantom,"="]
			&
			X
			\ar[d,"x_2"]
			\\
			B
			\sar[r,"Q"']
			\doublecell[rr, shift left=3ex]{\tab}
			&
			A
			\sar[r,"P"']
			& 
			B
		\end{tikzcd}
	\end{equation*}
	\begin{equation*}
		\begin{tikzcd}[column sep = small]
			&&&
			Z 
			\ar[ld,"z_1"'] 
			\ar[rd,"\chi", two heads]
			\ar[dd,phantom,"\mbox{\rotatebox[origin=c]{45}{$\lrcorner$}}"{marking,near start}]
			&&&
			\\
			&&
			U
			\ar[rd,"\varphi"', two heads]
			\ar[ld,"u_1"']
			&&
			W 
			\ar[ld,"w_1"']
			\ar[rd, "w_2"]
			\ar[dd,phantom,"\mbox{\rotatebox[origin=c]{45}{$\lrcorner$}}"{marking,near start}]
			&&
			\\ 
			&
			X
			\ar[ld,"x_1"']
			\ar[rr,phantom,"="]
			&&
			C 
			\ar[llld,"c_1"]
			\ar[rd,"c_2"]
			&&
			X 
			\ar[ld,"x_1"']
			\ar[rd,"x_2"]
			&
			\\
			A
			\sar[rr,"P"']
			&&
			B
			\sar[rr,"Q"']
			\doublecell[rr, shift left=3ex]{\tab}
			&&
			A
			\sar[rr,"P"']
			\doublecell[rr, shift left=3ex]{\tab}			
			&&
			B
		\end{tikzcd},
		\hspace{2ex}
		\begin{tikzcd}[column sep = tiny]
			&&&
			W 
				\ar[d,"\psi",two heads]
				\ar[llld,"w_1"']
				\ar[rrrd,"w_2"]
			\\
			C
				\ar[d,"c_1"']
			&&&
			G
				\ar[llld,"g_1"']
				\ar[rrrd,"g_2"]
			&&&
			X
				\ar[d,"x_2"]
			\\
			A 
				\sar[rr,"P"']
				\doublecell[rrrrrr, shift left=3ex]{\tab}
			&&
			B 
				\sar[rr,"Q"']
			&&
			A 
				\sar[rr,"P"']
			&&
			B	
		\end{tikzcd}
	\end{equation*}
	\begin{equation*}
		\begin{tikzcd}[column sep = small]
			&&&
			Z 
			\ar[rd,"z_2"] 
			\ar[ld,"\zeta"', two heads]
			\ar[dd,phantom,"\mbox{\rotatebox[origin=c]{45}{$\lrcorner$}}"{marking,near start}]
			&&&
			\\
			&&
			T
			\ar[ld,"t_1"']
			\ar[rd,"t_2"]
			\ar[dd,phantom,"\mbox{\rotatebox[origin=c]{45}{$\lrcorner$}}"{marking,near start}]
			&&
			V
			\ar[ld,"\theta"', two heads]
			\ar[rd, "v_2"]
			&&
			\\ 
			&
			X
			\ar[ld,"x_1"']
			\ar[rd,"x_2"]
			&&
			D
			\ar[ld,"\varepsilon"']
			\ar[rrrd,"\varepsilon"]
			\ar[rr,phantom,"="]
			&&
			X 
			\ar[rd,"x_2"]
			&
			\\
			A
			\sar[rr,"P"']
			\doublecell[rr, shift left=3ex]{\tab}
			&&
			B
			\sar[rr,"Q"']
			\doublecell[rr, shift left=3ex]{\tab}			
			&&
			A
			\sar[rr,"P"']
			&&
			B
		\end{tikzcd}
		,
		\hspace{2ex}
		\begin{tikzcd}[column sep = tiny]
			&&&
			T
				\ar[d,"\lambda",two heads]
				\ar[llld,"t_1"']
				\ar[rrrd,"t_2"]
			\\
			X
				\ar[d,"x_1"']
			&&&
			G
				\ar[llld,"g_1"']
				\ar[rrrd,"g_2"]
			&&&
			D
				\ar[d,"\varepsilon"]
			\\
			A 
				\sar[rr,"P"']
				\doublecell[rrrrrr, shift left=3ex]{\tab}
			&&
			B 
				\sar[rr,"Q"']
			&&
			A 
				\sar[rr,"P"']
			&&
			B	
		\end{tikzcd}
	\end{equation*}
	By \cref{lem:FactorisationOfSpan} and \cref{lem:FinStable},
	$\varphi$, $\psi$, $\chi$, $\theta$, $\zeta$, and $\lambda$ are all covers. 
	We would like to see 
	how the triangle identities of this adjoint behave.
	In the left diagram below,
	the unit cell is composited with $P$ horizontally, and the whole cell yields the corresponding arrow $\iota\fatsemi\psi$ from $X$ to $G$
	since $X$ and $G$ are exhibited as the tabulators of $P$ and $PQP$, respectively.
	Similarly, from the counit, we get the corresponding arrow from $G$ to $X$ as shown on the right below.
	
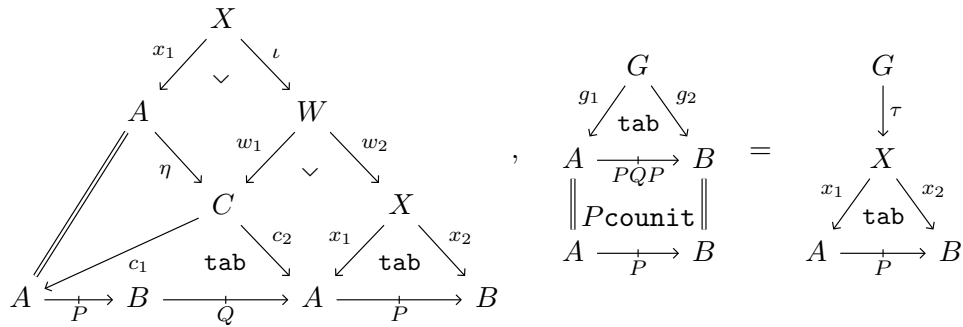
\begin{equation*}
			\begin{tikzcd}[column sep = small]
				&&&
				X
				\ar[ld,"x_1"'] 
				\ar[rd,"\iota"]
				\ar[dd,phantom,"\mbox{\rotatebox[origin=c]{45}{$\lrcorner$}}"{marking,near start}]
				&&&
				\\
				&&
				A
				\ar[rd,"\eta"']
				\ar[lldd,equal]
				&&
				W 
				\ar[ld,"w_1"']
				\ar[rd, "w_2"]
				\ar[dd,phantom,"\mbox{\rotatebox[origin=c]{45}{$\lrcorner$}}"{marking,near start}]
				&&
				\\ 
				&&&
				C
				\ar[llld,"c_1"]
				\ar[rd,"c_2"]
				&&
				X 
				\ar[ld,"x_1"']
				\ar[rd,"x_2"]
				&
				\\
				A
				\sar[rr,"P"']
				&&
				B
				\sar[rr,"Q"']
				\doublecell[rr, shift left=3ex]{\tab}
				&&
				A
				\sar[rr,"P"']
				\doublecell[rr, shift left=3ex]{\tab}			
				&&
				B
			\end{tikzcd},
			\hspace{2ex}
			\begin{tikzcd}[column sep = tiny]
				&
				G 
				\ar[ld,"g_1"']
				\ar[rd,"g_2"]
				&
				\\
				A
				\ar[d,equal] 
				\sar[rr,"PQP"']
				\doublecell[rr, shift left=3ex]{\tab}
				&&
				B
				\ar[d,equal]
				\\
				A
				\sar[rr,"P"']
				\doublecell[rr, shift left=3ex]{P\mathtt{counit}}
				&&
				B
			\end{tikzcd}
			\hspace{1ex}
			=
			\hspace{1ex}
			\begin{tikzcd}[column sep = tiny]
				&
				G 
				\ar[d,"\tau"]
				&
				\\
				& 
				X 
				\ar[ld,"x_1"']
				\ar[rd,"x_2"]
				&
				\\
				A
				\sar[rr,"P"']
				\doublecell[rr, shift left=3ex]{\tab}
				&&
				B
			\end{tikzcd}
	\end{equation*} 
	One of the triangle identities amounts to the equality of $\iota\fatsemi \psi\fatsemi \tau$ and $\id_X$.

	One observes that $z_1\fatsemi u_1$ and $\chi\fatsemi \psi\fatsemi \tau$ are the same arrows from $Z$ to $X$.
	To verify this equality, we compose the tabulator of $P$ to have the following equality.
	\begin{equation*}
		\begin{tikzcd}[column sep = tiny]
				&
				Z
				\ar[d,"\chi\fatsemi \psi\fatsemi \tau"]
					&
			\\
				& 
				X 
				\ar[ld,"x_1"']
				\ar[rd,"x_2"]
					&
			\\
			A
			\sar[rr,"P"']
			\doublecell[rr, shift left=3ex]{\tab}
				&
					&
					B
		\end{tikzcd}
		\hspace{0.2ex}
		=
		\hspace{0.2ex}
		\begin{tikzcd}[column sep = tiny]
			&
			Z 
			\ar[d,"\chi\fatsemi \psi"]
			&
			\\
			&
			G 
			\ar[ld,"g_1"']
			\ar[rd,"g_2"]
			&
			\\
			A
			\ar[d,equal] 
			\sar[rr,"PQP"']
			\doublecell[rr, shift left=3ex]{\tab}
			&&
			B
			\ar[d,equal]
			\\
			A
			\sar[rr,"P"']
			\doublecell[rr, shift left=3ex]{P\mathtt{counit}}
			&&
			B
		\end{tikzcd}
		\hspace{0.2ex}
		=
		\hspace{0.2ex}
		\begin{tikzcd}[column sep = tiny]
			&
			Z
			\ar[ld,"z_3"']
			\ar[rd,"z_4"]
			&
			\\
			A
			\ar[d,equal] 
			\sar[rr,"PQP"']
			\doublecell[rr, shift left=3ex]{\tab}
			&&
			B
			\ar[d,equal]
			\\
			A
			\sar[rr,"P"']
			\doublecell[rr, shift left=3ex]{P\mathtt{counit}}
			&&
			B
		\end{tikzcd}
		\hspace{0.2ex}
		=
		\hspace{0.2ex}
		\begin{tikzcd}[row sep = small, column sep = small]
			&&&
			Z 
			\ar[rd,"z_2"] 
			\ar[ld,"\zeta"', two heads]
			\ar[dd,phantom,"\mbox{\rotatebox[origin=c]{45}{$\lrcorner$}}"{marking,near start}]
			&
			\\
			&&
			T
			\ar[ld,"t_1"']
			\ar[rd,"t_2"]
			\ar[dd,phantom,"\mbox{\rotatebox[origin=c]{45}{$\lrcorner$}}"{marking,near start}]
			&&
			V
			\ar[ld,"\theta"', two heads]
			\ar[d, "v_2"]
			\\ 
			&
			X
			\ar[ld,"x_1"']
			\ar[rd,"x_2"]
			&&
			D
			\ar[ld,"\varepsilon"']
			\ar[rd,"\varepsilon"]
			\ar[r,phantom,"="]
			&
			X 
			\ar[d,"x_2"]
			\\
			A
			\ar[d, equal]
			\sar[rr,"P"']
			\doublecell[rr, shift left=2ex]{\tab}
			&&
			B
			\ar[d, equal]
			\sar[rr,"QP"']
			\doublecell[rr, shift left=2ex]{\tab}
			&&
			B
			\ar[d,equal]
			\\
			A 
			\sar[rr,"P"']
			\doublecell[rr, shift left=2ex]{\verteq}
			&&
			B 
			\sar[rr, equal]
			\doublecell[rr, shift left=2ex]{\mathtt{counit}}
			&&
			B
		\end{tikzcd}
		\hspace{0.2ex}
		=
		\hspace{0.2ex}
		\begin{tikzcd}[column sep = tiny]
			&
			Z
			\ar[d,"z_1\fatsemi u_1"]
			&
			\\
			& 
			X 
			\ar[ld,"x_1"']
			\ar[rd,"x_2"]
			&
			\\
			A
			\sar[rr,"P"']
			\doublecell[rr, shift left=3ex]{\tab}
			&&
			B
		\end{tikzcd}
	\end{equation*}
	We present one more diagram in which all faces of the cube on the top are pullback squares.
	\begin{equation*}
		\begin{tikzcd}[row sep = small, column sep = small]
			E
			\ar[rr,"e_1",two heads]
			\ar[dd,"e_2"']
			\ar[rd,"e_3"]
			&
			&
			X 
			\ar[dd,"\iota",near end]
			\ar[rd,"x_1"]
			& 
			\\
			&
			F 
			\ar[rr,"f_1",crossing over, near start,two heads]
			&
			&
			A 
			\ar[dd,"\eta"']
			\ar[dddd, bend left=30, equal]
			\\
			Z 
			\ar[rr,"\chi",near end, two heads]
			\ar[rd,"z_1"']
			&
			&
			W 
			\ar[rd,"w_1"']
			&
			\\
			&
			U 
			\ar[dd,"u_1"']
			\ar[rr,"\varphi",two heads,near start]
			\ar[from=uu,"f_2"',crossing over, near start]
			&
			&
			C 
			\ar[dd,"c_1"']
			\\
			\\
			& 
			X
			\ar[rr,"x_1"']
			&&
			A
		\end{tikzcd}
	\end{equation*}
	We have the following equalities.
	\begin{align*}
		E \to["e_3"] F \to["f_2"] U \to["u_1"] X &= E \to["e_2"] Z \to["z_1"] U \to["u_1"] X \\
		&= E \to["e_2"] Z \thto["\chi"] W \to["\psi"] G \to["\tau"] X & (\text{From the argument above})\\
		&= E \thto["e_1"] X \to["\iota"] W \to["\psi"] G \to["\tau"] X \\
		&= E \thto["e_1"] X & (\text{The consequence of the triangle identity})
	\end{align*}
	Since a cover is an epimorphism, $e_1$ is an epimorphism, and so is $F\tonormal["f_2\fatsemi u_1"] X$.
	We know that $F\thto["f_1"] A= F \tonormal["f_2\fatsemi u_1"] X \to["x_1"] A$ is a cover,
	and hence, we conclude that $x_1$ is a cover as well
	since, by the general theory of orthogonal factorisation systems, the property of being a cover is left-cancellable by epimorphisms.

	By the same argument using the other part of the triangle identity, 
	we can show that $y_2$ is a cover.
	The stability of covers leads $v_1$ and $v_2$ to be covers.
	We have the following two opcartesian cells.
	\begin{equation*}
			\begin{tikzcd}[column sep = small,row sep = small]
				&
					V
						\ar[d,"v_1"]
				\\
				&
					Y
						\ar[ld,"y_1"']
						\ar[rd,"y_2"]
				&
				\\
					B
						\sar[rr,"Q"']
						\doublecell[rr, shift left=2ex]{\tab}
				&
				&
					A 
			\end{tikzcd}
			\hspace{1ex}
			,
			\hspace{1ex}
			\begin{tikzcd}[column sep = small,row sep = small]
				&
					V
						\ar[d,"v_2"]
				\\
				&
					X
						\ar[ld,"x_1"']
						\ar[rd,"x_2"]
				&
				\\
					A
						\sar[rr,"P"']
						\doublecell[rr, shift left=2ex]{\tab}
				&
				&
					B
			\end{tikzcd}
	\end{equation*}
	Here, $v_1\fatsemi y_2=v_2\fatsemi x_1$, by definition, and $v_1\fatsemi y_1=\theta\fatsemi \varepsilon=v_2\fatsemi x_2$.
	This means that $Q=P^\dag$ by \cref{prop:pseudoFunctorialityDagger,prop:ConjointDagger}.
	By \cref{lem:UnbiasingLemma}, we suppose without loss of generality that $X=Y$, 
	$x_1=y_2$ and $x_2=y_1$ and also $v_1=v_2$ by the universal property of the tabulator.
	However, the pair of $v_1$ and $v_2$ can be seen as a kernel pair of $x_1$.
	Therefore, $x_1$ is monic.
\end{proof}

\begin{remark} 
	As the reviewer pointed out, this theorem has already been proven in Milius's master thesis \cite[Corollary 4.21]{Mil00}.
	However, the proof in this paper is different from the one in \cite{Mil00}
	in that we did not use the notions of \textit{single-valuedness} and \textit{totality} of a relation
	and instead used the condition of being an adjoint pair as it is.
\end{remark}

\begin{lemma}
	\label{lem:Cauchyantimono}
	Let $\dbl{D}$ be a unit-pure double category of relations.
	Then $\dbl{D}$ is Cauchy if and only if the accompanying factorisation system is anti-right-proper.
\end{lemma}

\begin{proof}
	Suppose $\dbl{D}= \Rel{\zero{E},\zero{M}}(\one{C})$ is Cauchy.
	Observe that $\zero{M}$ is the class of all vertical arrows right orthogonal to all arrows in $\zero{E}$, and $\zero{E}$ is equal to $\Cov(\dbl{D})$. 
	So $\zero{M}$ contains all monomorphisms by \cref{lem:covincortho,lem:monoinc}.

	Conversely, assume the right class of the accompanying factorisation system $\zero{M}$ contains all monomorphisms.
	Then, monic covers are all isomorphisms.
	Therefore, every map is representable by \cref{thm:mapunitpure}.
\end{proof}

\begin{theorem}
	\label{thm:Cauchydouble}
	The following are equivalent for a double category $\dbl{D}$.
	\begin{enumerate}
		\item $\dbl{D}$ is equivalent 
		to $\Rel{\zero{E},\zero{M}}(\one{C})$ for some category $\one{C}$
		with finite limits and an anti-right-proper stable factorisation system $(\zero{E},\zero{M})$.
		\item $\dbl{D}$ is a unit-pure, Cauchy double category of relations.
	\end{enumerate}
\end{theorem}

\begin{proof}
	The direction $ii)\Rightarrow i)$ is a consequence of \cref{thm:MainThm} and \cref{lem:Cauchyantimono}.
	We prove the other direction $i)\Rightarrow ii)$.
	If $\zero{M}$ contains every monomorphism, then every vertical arrow in $\zero{E}$ is a strong epimorphism, hence an epimorphism.
	Thus, the factorisation system is left-proper, and the double category is unit-pure by \cref{thm:unitpuredbl}.
	By \cref{lem:Cauchyantimono}, the double category is Cauchy.
\end{proof}

\begin{remark}
    \label{rem:RegepiMon}
    A similar result to \Cref{thm:Cauchydouble} already exists \cite[Theorem 5.1]{Pav95} in the context of bicategories.
    If we rephrase it in the language of double categories,
    the conditions of \Cref{thm:Cauchydouble} are also equivalent to the following:
    $\dbl{D}$ is equivalent 
    to $\Rel{\zero{E},\zero{M}}(\one{C})$ for some category $\one{C}$
    with finite limits and a stable factorisation system $(\zero{E},\zero{M})$
    where $\zero{E}$ is induced by the class of all the regular epimorphisms $\Regepi$ in $\one{C}$.
	A more direct and simplified proof of this fact is given in \cite[Theorem 4.23]{Mil00}.
    Since the factorisation system $(\zero{E},\zero{M})$ is uniquely determined by the structure of double categories,
    we can conclude that a stable factorisation system $(\zero{E},\zero{M})$ on a category with finite limits is anti-right-proper if and only if
    $\zero{E}$ is induced by $\Regepi$.
    This fact was conjectured to be false in general in \cite{Pav95}.
    What we have shown is that it is indeed true, but we provide a direct proof of this fact as follows.
\end{remark}

\begin{proposition}
    \label{prop:RegepiMon}
    Let $\one{C}$ be a category with finite limits.
    Then a stable factorisation system $(\zero{E},\zero{M})$ on $\one{C}$ is anti-right-proper 
    if and only if $\zero{E}$ is included in the class of regular epimorphisms $\Regepi$.
\end{proposition}

\begin{proof}
    If $\zero{E}$ is included in $\Regepi$, then all the arrows in $\zero{E}$ are orthogonal to all the monomorphisms,
    which implies that the factorisation system $(\zero{E},\zero{M})$ is anti-right-proper.

    Conversely, suppose that $(\zero{E},\zero{M})$ is anti-right-proper,
    and let $f\colon A\to B$ be an arrow in $\zero{E}$.
    Since $\one{C}$ has finite limits and any arrow orthogonal to all the monomorphisms is an epimorphism,
    $f$ is an epimorphism.
    We take the pullback of $f$ along itself, in other words, the kernel pair of $f$,
    and name it as $k_1,k_2\colon K\to A$.
    We prove that $f$ is a coequaliser of $k_1$ and $k_2$.
    Let $g\colon A\to C$ be an arrow such that $g\circ k_1=g\circ k_2$.
    Since $f$ is an epimorphism, it suffices to show that $g$ factors through $f$.

    The pairing $\langle f,g\rangle\colon A\to B\times C$ has the same kernel pair $k_1,k_2$ as $f$
    because $g$ coequalises $k_1$ and $k_2$.
    Take the $(\zero{E},\zero{M})$-factorisation of $\langle f,g\rangle$ as 
    the composite of $e\colon A\to D$ in $\zero{E}$ and $\langle m,n\rangle\colon D\to B\times C$ in $\zero{M}$.
    Then, take the kernel pairs of $\langle m,n\rangle$ and $m$ and we have the following diagrams:
    \begin{equation}
        \label{eq:RegepiMon}
        \begin{tikzcd}
            K
            \ar[r, two heads]
            \ar[d, two heads]
            \ar[rr, bend left, two heads, "k_1"]
            \ar[dd, bend right, two heads, "k_2"']
            \pullback
            &
            \cdot 
            \ar[r]
            \ar[d, two heads]
            \pullback
            &
            A
            \ar[d, "e", two heads]
            \\
            \cdot
            \pullback
            \ar[r, two heads]
            \ar[d]
            &
            X 
            \pullback
            \ar[r, tail, "a"]
            \ar[d, tail, "b"']
            &
            D 
            \ar[d, tail, "{\langle m,n\rangle}"']
            \\
            A
            \ar[r, "e"']
            &
            D 
            \ar[r, tail,  "{\langle m,n\rangle}"']
            &
            B\times C
        \end{tikzcd}
        \quad
        ,   
        \quad
        \begin{tikzcd}
            K
            \ar[r, two heads]
            \ar[d, two heads]
            \ar[rr, bend left, two heads, "k_1"]
            \ar[dd, bend right, two heads, "k_2"']
            \pullback
            &
            \cdot 
            \ar[r]
            \ar[d, two heads]
            \pullback
            &
            A
            \ar[d, "e", two heads]
            \\
            \cdot
            \pullback
            \ar[r, two heads]
            \ar[d]
            &
            Y 
            \pullback
            \ar[r, two heads, "x"]
            \ar[d, two heads, "y"']
            &
            D 
            \ar[d, two heads, "{m}"]
            \\
            A
            \ar[r, "e"']
            &
            D 
            \ar[r, two heads,  "{m}"']
            &
            B
        \end{tikzcd}
        .
    \end{equation}
    By the universal property of the pullbacks, there exist canonical arrows $p\colon K \to X$ and $u\colon X\to Y$ 
    which make the following diagram commute:
    \[
        \begin{tikzcd}[row sep=small, column sep=small]
            K 
            \ar[rrrr,"k_1"]
            \ar[dddd,"k_2"']
            \ar[rd,"p"]
            &&&&
            A
            \ar[dd,"e"]
            \\
            &
            X 
            \ar[rd,"u "]
            \ar[rrrd, bend left, "a",tail]
            \ar[dddr, bend right, "b"',tail]
            \\
            &&
            Y
            \ar[rr,"x",two heads]
            \ar[dd,"y"',two heads]
            &&
            D
            \ar[dd,"m"]
            \\
            \\
            A
            \ar[rr,"e"']
            &&
            D 
            \ar[rr,"m"']
            &&
            B 
        \end{tikzcd}
        .
    \]
    Note that $x,y$ are in $\zero{E}$ by its stability under pullbacks,
    and $a,b$ are in $\zero{M}$.
    Then, $p$ and $p\fatsemi u$ are in $\zero{E}$ since these arrows fit into the diagrams \eqref{eq:RegepiMon} 
    and so they are in $\zero{E}$.
    Again by the cancellation property of $\zero{E}$, $u$ belongs to $\zero{E}$.
    Thus, $a=u\fatsemi x$ and $b=u\fatsemi y$ are in $\zero{E}$,
    but since they are also in $\zero{M}$, they are isomorphisms.
    Therefore, $u$ is a split monomorphism, hence an isomorphism, and so are $x$ and $y$.
    This implies that $m$ is monomorphic, and because of $e\fatsemi m=f$ and the cancellation property of $\zero{E}$,
    it also belongs to $\zero{E}$, which means that $m$ is an isomorphism.
    Therefore, $g=f\fatsemi(m\inv\fatsemi n)$, which shows that $f$ is a coequaliser of $k_1$ and $k_2$.
\end{proof}

Consequently, we have the following characterisation of double categories of relations on regular categories.
Here, we recapture the result in \cite{Lam22} as the equivalence of $i)$ and $iv)$.

\begin{theorem}[{\cite[Theorem 10.2]{Lam22}}]
	\label{thm:Regcat}
	The following are equivalent for a double category $\dbl{D}$.
	\begin{enumerate}
		\item $\dbl{D}$ is equivalent to $\Rel{\Regepi,\Mono}(\one{C})$ for some regular category $\one{C}$.
		\item $\dbl{D}$ is a locally posetal, Cauchy double category of relations.
		\item $\dbl{D}$ is a unit-pure, locally preordered, Cauchy double category of relations.
		\item $\dbl{D}$ is a locally posetal discrete cartesian equipment with a $\Mono$-comprehension scheme.
		\item $\dbl{D}$ is a locally posetal discrete cartesian equipment with a left-sided $\Mono$-comprehension scheme.
	\end{enumerate}
\end{theorem}

\begin{proof}
	The equivalence between $ii)$ and $iii)$ is a direct consequence of \cref{lem:DagLocPos}.
	One can deduce the equivalence of $i)$ and $ii)$ by 
	combining \cref{cor:locpos} and \cref{thm:Cauchydouble}
	since they show the equality of $\Mono$ and $\zero{M}$.
	Note that regular categories are precisely those categories with stable factorisation systems whose right class is the class of monomorphisms.
	The implication from $i)$ to $iv)$ and $v)$ is a consequence of \cref{thm:MainThm},
	and the opposite directions follow from \cref{lem:DagLocPos} and \cref{thm:unitpuredbl}. 
\end{proof}

In light of this theorem, 
unit-pure Cauchy double categories of relations can be considered as a persuasive generalisation
of the double categories of relations on regular categories.
As observed in \cite{CKS84},
the double category of spans in a finitely complete category is also an example of a Cauchy double
category of relations. 

In the paper \cite{Sch15}, the concept of \textit{regular double categories} is defined
differently from the definition of Cauchy double categories of relations.
A regular double category is defined in the paper
using what the author called \textit{normal collapses}.
The difference between the two double categorical generalisations for regular categories is rooted in 
the difference between the two definitions of regular categories.
One definition regards a regular category as a category with finite limits and coequalisers of kernel pairs in which
regular epimorphisms are stable under pullbacks,
and the other definition regards it 
as a category with finite limits and stable images.
As a classical result, these two definitions for a category with finite limits are known to be equivalent.
The concept of normal collapses is a translation of the concept of coequalisers of kernel pairs 
into the language of (virtual) equipments, and hence
the work in \cite{Sch15} can be seen as a generalisation of the former definition in terms of regular double categories.
On the other hand, the latter definition is closely related to our generalisation,
since we can readily rephrase the definition in terms of factorisation systems:
a regular category is a category with a stable factorisation system whose morphisms in the right class are monomorphisms.
How these two ways of generalisation are related is a natural question, but is not discussed in this paper.

In the rest of this subsection, we would like to discuss a process to obtain a Cauchy double category of relations from a general double category of relations. 
To that end, we restrict the double categories in question to be unit-pure. 
The work by Kelly \cite{Kel91} can be understood as this process in the case of proper factorisation systems,
although the approach was taken in the ordinary category theory. 
The following discussion extends it to unit-pure double categories of relations.

Note that when an equipment $\dbl{D}$ is Cauchy and unit-pure, 
vertical arrows are in one-to-one correspondence with isomorphism classes of adjoint pairs in the horizontal bicategory $\bi{H}(\dbl{D})$.
This is because if two vertical arrows define the isomorphic adjoint pairs, then they are isomorphic as 1-cells in $\bi{V}(\dbl{D})$ and unit-pureness implies that they are equal.
This observation leads us to the following definitions.

\begin{definition}\label{def:mapdoublecategory}
	Let $\bi{B}$ be a bicategory, and $\dbl{D}$ be a double category.
	We say $\dbl{D}$ is a \emph{map double category} of $\bi{B}$
	if $\biH{\dbl{D}}=\bi{B}$ and $\dbl{D}$ is a Cauchy equipment.
\end{definition}

\begin{remark}
	In the literature, the map bicategory 
	$\biMap(\bi{B})$ of a bicategory $\bi{B}$ is defined as the locally full sub-bicategory of $\bi{B}$ whose objects are the objects of $\bi{B}$ and whose 1-cells are the 1-cells with right adjoints in $\bi{B}$.
	See, for example, \cite{Ale18} for the relation between these two concepts.
\end{remark}

Since vertical arrows in Cauchy equipments are `almost the same' as
left adjoints in the horizontal bicategory,
map double categories are almost identical for fixed bicategory $\bi{B}$.
However, the equivalences in $\DblCat$ do not conceive this `sameness'
since the 2-functoriality of $\bi{V}$ implies that
equivalences in $\DblCat$ restrict to equivalences in $2\Catbi$; i.e.,
2-equivalences rather than \textit{biequivalences}.
But if we restrict map double categories to sufficiently simple ones,
the 2-category $\DblCat$ gives a sufficient framework to understand
\textit{the} map double category of a bicategory.

\begin{proposition} 
	Let $\bi{B}$ be a bicategory.
	If $\dbl{D}$ and $\dbl{D}'$ are map double categories of $\bi{B}$ that are unit-pure,
	then $\dbl{D}$ and $\dbl{D}'$ are isomorphic in $\DblCat$.
\end{proposition}
\begin{proof}
	In general, observe that for each left adjoint $p\colon A \sto B$,
	up to an invertible vertical cell,
	there is at most one vertical arrow $f\colon A\to B$ such that
	$p$ is a companion of $f$.
	In particular, if $\dbl{D}$ is unit-pure, such a vertical arrow is actually
	unique if exists.
	Moreover, when we fix such a pair $(p,f)$,
	the canonical cartesian and opcartesian cells $\alpha$ and $\beta$ as shown in \cref{prop:repradj}, exhibiting $p$ as a companion of $f$,
	are also unique. To check this, suppose two such pairs $(\alpha,\beta)$ and $(\alpha',\beta')$ are given.
	Since $\dbl{D}$ is unit-pure,
	two vertical cells obtained by postcomposing $\alpha$ to $\beta$ and $\beta'$ are the same. Therefore, we have $\beta=\beta'$ for $\alpha$ is cartesian.
	Similar arguments hold for conjoints considering the horizontal opposite. Therefore, utilising \cref{prop:repradj},
	we observe that the unit and counit of the horizontal adjoint are also unique if we fix a companion $f_!$ and a conjoint $f^*$.

	Suppose we have two unit-pure map double categories $\dbl{D}$ and $\dbl{D}'$ of $\bi{B}$.
	Choose adjoints $f_!\dashv f^*$ in $\biH{\dbl{D}}$ for each $f$, and
	define a pseudo-functor $F\colon\dbl{D}\to\dbl{D}'$ as follows.
	\begin{itemize}
		\item%
			$F$ is the identity on the horizontal bicategory.
		\item%
			For each vertical arrow $f\colon A\to B$ in $\dbl{D}$,
			$F(f)\colon A\to B$ is defined as the unique vertical arrow in $\dbl{D}'$ whose companion is $f_!$.
			Note that $F(f)$ is also the unique vertical arrow whose conjoint is $f^*$.
		\item%
			Observe that any cell in $\dbl{D}$ is uniquely decomposed as the following composite.
			\[
				\begin{tikzcd}
					W
					\sar[rrr, "q"]
					\ar[d, equal]
					\doublecell[rrrd]{v}
						&
							&
								&
								X
								\ar[d, equal]
					\\
					W
					\ar[rd, "f"']
					\sar[r, "f_!"]
					\doublecell[rd, shift left=2.0ex]{\alpha}
						&
						Y
						\doublecell[rd]{\verteq}
						\sar[r,"p"]
						\ar[d, equal]
							&
							Z
							\ar[d, equal]
							\sar[r, "g^*"]
								&
								X
								\doublecell[ld, shift right=2.0ex]{\gamma}
								\ar[ld, "g"]
					\\
						&
						Y
						\sar[r, "p"']
							&
							Z
								&
				\end{tikzcd}
			\]
			The discussion above asserts that $F(\alpha)$ and $F(\gamma)$ are uniquely determined, and this defines how $F$ assigns cells.
		\item%
			A tedious but straightforward discussion shows that this indeed defines a pseudo-functor.
			For example, to show the horizontal functoriality, it suffices to check $F$ sends counits of the adjoint induced from
			selected pairs of companions and conjoints to those induced from $F(\alpha)$ and $F(\gamma)$, which follows from what we discussed above.
	\end{itemize}
	$F$ is invertible because of the uniqueness appearing in the construction,
	and this completes the proof.
\end{proof}

\begin{definition}
	Let $\dbl{D}$ be an equipment.
	An equipment $\dbl{E}$ is a \emph{Cauchisation} of $\dbl{D}$
	if $\bi{H}(\dbl{E})=\bi{H}(\dbl{D})$ holds and $\dbl{E}$ is Cauchy and unit-pure,
	and denoted by $\Cau(\dbl{D})$.
	In terms of \cref{def:mapdoublecategory},
	$\Cau(\dbl{D})$ is a unit-pure map double category of $\bi{H}(\dbl{D})$.
\end{definition}

We sometimes use the word `map' to refer to an isomorphism class of adjoint pairs in $\bi{H}(\dbl{D})$,
or equivalently, a vertical arrow in $\Cau(\dbl{D})$.
We have named this `Cauchisation' because of the following proposition.

\begin{proposition}
	\label{prop:univpropcauchy}
	Let $\dbl{D}$ be an equipment and $\dbl{E}$ be a Cauchisation of $\dbl{D}$.
	Then, we have a canonical pseudo-functor $I\colon\dbl{D}\to\dbl{E}$.
	Moreover, for any Cauchy unit-pure equipment $\dbl{F}$ and a pseudo-functor $F\colon\dbl{D}\to\dbl{F}$,
	there exists a unique pseudo-functor $G\colon\dbl{E}\to\dbl{F}$ such that $F=G\circ I$.
\end{proposition}

\begin{proof}
	We define $I$ as the identity on the horizontal bicategory and send each vertical arrow $f$ to $(\id,f)$.
	Since a cell in $\dbl{D}$ of the form below is in one-to-one correspondence with a horizontal cells 
	$\widetilde{\alpha}\colon Rg_!\Rightarrow f_!S$ in $\bi{H}(\dbl{D})$,
	so $I$ sends such a cell to the cell on the right below.
	Note that a pseudo-functor preserves companions.
	\[
	\begin{tikzcd}
		A 
			\ar[d,"f"']
			\sar[r,"R"]
			\ar[dr, phantom, "\alpha"]
		&
			B
			\ar[d,"g"]
		\\
			C
			\sar[r,"S"']	
		&
			D
	\end{tikzcd}
	\hspace{2ex},
	\hspace{2ex}
	\begin{tikzcd}[column sep = large, row sep = small]
		A 
			\ar[d,equal]
			\sar[r,"R"]
			\doublecell[rd]{\verteq}
		&
		B 
			\ar[rd,"Ig",bend left=15]
			\ar[d,equal]
		\\
		A
			\sar[r,"R"']
		    \ar[d,equal]
			\ar[drr, phantom, "\widetilde{\alpha}"]
		&
		B 
			\sar[r,"g_!"']
			\ar[r, phantom, shift left=1.5ex,xshift=-1ex, "\opcart"]
		&
		D
			\ar[d,equal]
		\\
		A
			\sar[r,"f_!"]
			\ar[rd, "If"', bend right=15]
			\ar[r, phantom, shift right=1.5ex, xshift=1ex,"\cart"]
			&
		C
			\sar[r,"S"]
			\ar[d,equal]
			\doublecell[rd]{\verteq}
		&
		D
			\ar[d,equal]
		\\
		&
		C
			\sar[r,"S"']
		&
		D
	\end{tikzcd}
	\]
	Thus defined $I$ is easily shown to be a pseudo-functor.

	In a unit-pure Cauchy equipment, a vertical arrow is uniquely determined by its representative adjoint pair.
	Therefore, for any pseudo-functor $F\colon\dbl{D}\to\dbl{F}$,
	assignment of the image of a vertical arrow in $\dbl{E}$ is uniquely determined by the image of its representative adjoint pair.
	Cells are also in one-to-one correspondence with horizontal cells by composing the companion cells,
	which implies that $G$ is uniquely determined and also well-defined.
\end{proof}

We move on to the construction of $\Cau(\dbl{D})$ for 
a unit-pure double category of relations $\dbl{D}$.
For a monic cover $m$, we see that $(m_!,m^*)$ is an adjoint equivalence in $\bi{H}(\dbl{D})$.
Therefore, $m^*$ gives a vertical isomorphism in $\Cau(\dbl{D})$.
Recall from \cref{thm:mapunitpure} that every map in $\bi{H}(\dbl{D})$ is of the form $m^*f_!$ for some monic cover $m$ and some vertical arrow $f$.
Before we construct $\Cau(\dbl{D})$, we probe into the properties of monic covers.

\begin{remark}
	As shown in \Cref{lem:CovEqualFin},
	a class of covers $\Cov(\dbl{D})$ and a class of final arrows $\Fin(\dbl{D})$ are the same in a discrete cartesian equipment $\dbl{D}$,
	and hence $\rclass{\Cov(\dbl{D})}$ is the same as $\Fib(\dbl{D})$ in a double category of relations $\dbl{D}$.
	When we address the left class of the accompanying factorisation system of a double category of relations,
	we focus on the aspect as the class of covers, so we use the notation $\Cov(\dbl{D})$.
	Although we write $\Fib(\dbl{D})$ for the right class of the accompanying factorisation system,
	an important property of the class of fibrations in a double category of relations
	for the sequel is that it is the class of vertical arrows right orthogonal to all covers.
\end{remark}

\begin{lemma}\label{lem:monocovcancel} 
	Let $\dbl{D}$ be a unit-pure double category of relations.
	The class of monic covers is closed under composition and 
	right-cancellable, that is, if $m\fatsemi n$ and $n$ are monic covers, then $m$ is also a monic cover.
\end{lemma}
\begin{proof}
	It is straightforward to check that the class of monic covers is
	closed under composition.
	Take composable arrows $m\colon A\to B$ and $n\colon B\to C$;
	suppose $m\fatsemi n$ and $n$ are monic covers.
	Then, $m$ is a monic cover, being a pullback of $m\fatsemi n$ through $n$.
\end{proof}

\begin{corollary}\label{cor:existcauchy}
	A unit-pure double category of relations $\dbl{D}$ has a Cauchisation $\Cau(\dbl{D})$.
	It comprises the following data.
	\begin{itemize}
		\item An object is an object of $\dbl{D}$.
		\item A vertical arrow from $A$ to $B$ is an isomorphism
		class of $\Fib(\dbl{D})$-relations 
		$\langle m,f\rangle$ where $m\colon X\to A$ is a monic cover and $f\colon X\to B$ is a vertical arrow.
		We write it just as $(m,f)$.
		\item A horizontal arrow from $A$ to $B$ is a horizontal arrow from $A$ to $B$ in $\dbl{D}$.
		\item A cell of the form on the left below is a cell in $\dbl{D}$ of the form on the right below:
		\[
		\begin{tikzcd}
			A 
				\ar[d,"{(m,f)}"']
				\sar[r,"R"]
			&
				B
				\ar[d,"{(n,g)}"]
			\\
				C
				\sar[r,"S"']	
			&
				D
		\end{tikzcd},
		\hspace{2ex}
		\begin{tikzcd}
			A 
				\ar[d,equal]
				\sar[r,"R"]
			&
			B 
				\sar[r,"n^*g_!"]
			&
			D
				\ar[d,equal]
			\\
			A
				\sar[r,"m^*f_!"']
			&
			C
				\sar[r,"S"']
			&
			D
		\end{tikzcd}
		\]
	\end{itemize}
\end{corollary} 
\begin{proof} 
Restricting the adjoint equivalence \cref{eq:adjequivforcomp} in the case of $\zero{M}=\Fib(\dbl{D})$,
we obtain the one to one correspondence between vertical arrows defined above and isomorphism classes of adjoint pairs in $\bi{H}(\dbl{D})$.
The composition of vertical arrows is derived from the composition of 1-cells in $\bi{H}(\dbl{D})$ through the adjoint equivalence.
Cells are defined above in a unique way to make $m^*f_!$ a companion of $(m,f)$ for every vertical arrow $(m,f)$.
Composition of cells is defined similarly to the composition of cells in the double category of quintets \cite{GP04},
but since $\bi{H}(\dbl{D})$ is a bicategory, we need to take isomorphism classes of adjoint pairs into account.
The composition of cells then does not depend on the choice of representatives.
The essential point of the proof for this to be a double category is that the vertical composition is strictly associative and unital.
The composition of vertical arrows is strictly associative since we take the isomorphism classes of adjoint pairs.
Since the composition of 2-cells in $\bi{H}(\dbl{D})$ is strictly associative,
the vertical composition of cells is also strictly associative.
Since $m^*f_!$ has a right adjoint $f^*m_!$, it is an equipment.
By \cref{thm:mapunitpure}, the double category $\Cau(\dbl{D})$ is Cauchy.
Therefore, the data above constitute a Cauchy equipment.

To see that it is unit-pure, suppose there is a vertical cell of the following form on the left;
then, we have the cell in $\dbl{D}$ on the right.
\[
		\begin{tikzcd}
			A 
				\ar[d,"{(m,f)}"']
				\sar[r,equal]
				\doublecell[rd]{k}
			&
				A
				\ar[d,"{(n,g)}"]
			\\
				B
				\sar[r,equal]	
			&
				B
		\end{tikzcd},
		\hspace{2ex}
		\begin{tikzcd}
			A 
				\ar[d,equal]
				\sar[r,"n^*g_!"]
				\doublecell[rd]{k}
			&
			B 
				\ar[d,equal]
			\\
			A
				\sar[r,"m^*f_!"']
			&	
			B 
		\end{tikzcd}
		\]
		where $(m\colon X\to A,f\colon X\to B),\ (n\colon Y\to A,g\colon Y \to B)$ are $\Fib(\dbl{D})\mhyphen$ relations, and $k\colon X\to Y$ is 
		a vertical arrow corresponding to the cell.
		By the equality $k\fatsemi n=m$, we know that $k$ is a monic cover by \cref{lem:monocovcancel}.
		Furthermore, it follows from $k\fatsemi \langle n,g\rangle=\langle m,f\rangle$ 
		that $k$ belongs to $\Fib(\dbl{D})$ by its cancellability.
		Therefore, $k\colon X\to Y$ is an isomorphism,
		meaning that $(m,f)$ and $(n,g)$ are the same as a vertical arrow in $\Cau(\dbl{D})$.
\end{proof}

According to the proof,
the composition of maps can be manifestly presented as follows.
Given two vertical arrows $(m,f)\colon A\to B$ and $(n,g)\colon B\to C$ in $\Cau(\dbl{D})$,
and take the pullback of $f$ and $n$ as below.
\begin{equation*}
	\begin{tikzcd}[column sep = tiny,row sep =small]
			&
				&
				Z
				\ar[ld,"k"']
				\ar[rd,"t"]
				\ar[dd,phantom,"\mbox{\rotatebox[origin=c]{45}{$\lrcorner$}}"{marking,near start}]
					&
		\\
			&
			X
			\ar[ld,"m"']
			\ar[rd,"f"']
				&
					&
						Y
							\ar[ld,"n"]
							\ar[rd,"g"]
		\\
		A
			&
				&
				B
					&
						&
							C
	\end{tikzcd}
\end{equation*}
Since $m$ and $n$ are monic covers, so is $k\fatsemi m$.
However, we do not know whether $\langle k\fatsemi m, t\fatsemi g\rangle$ 
belongs to $\Fib(\dbl{D})$.
All we can say is that the composite of $(m,f)$ and $(n,g)$ corresponds to $(k\fatsemi m)^*(t\fatsemi g)_!$.
In this way, it is sometimes useful to write vertical arrows as the form of a map $m^*f_!$ for some monic cover $m$ and some vertical arrow $f$,
even when $\langle m,f\rangle$ is not in $\Fib(\dbl{D})$.
Note that, as a vertical arrow in $\Cau(\dbl{D})$, $m^*f_!$ is the composite of $m^*$ and $f_!$;
this can be checked by the previous observation on the composition of maps.

\begin{remark}
A canonical pseudo-functor $\dbl{D}\to\Cau(\dbl{D})$ exists
and is the identity on the horizontal part and sends each vertical arrow to its companion.
This pseudo-functor is faithful as a functor $\dbl{D}_0\to\Cau(\dbl{D})_0$
since the original double category $\dbl{D}$ is unit-pure.
The functor $\dbl{D}_1\to\Cau(\dbl{D})_1$ obtained by the above pseudo-functor is faithful as well.
This is verified by using the cells defining the companions of vertical arrows.
\end{remark}

\begin{remark}
	\label{rem:equaldomain}
	In $\Cau(\dbl{D})$, every vertical arrow is of the form $m^*f_!$ for some monic cover $m$ and some vertical arrow $f$,
	which means that it is in the image of the canonical pseudo-functor $\dbl{D}\to\Cau(\dbl{D})$ up to a precomposition of isomorphism
	since $m^*$ has an inverse $m_!$ in $\oneV{\Cau(\dbl{D})}$.
	Furthermore, if we have two vertical arrows $p=(m,t)\colon A\to B$ and $q=(n,s)\colon A \to C$ in $\Cau(\dbl{D})$ with the same domain, 
	we can take a monic cover $m\colon X\to A$ in $\dbl{D}$ such that $R$ and $S$ are presented as $k^*f_!$ and $k^*g_!$ respectively for some vertical arrows $f\colon X\to B$ and $g\colon X\to C$,
	even though $\langle k,f\rangle$ and $\langle k,g\rangle$ are not necessarily in $\Fib(\dbl{D})$.
	This is achieved by taking the pullback of $m$ and $n$.
\end{remark}

\begin{lemma}
	Let $\dbl{D}$ be a unit-pure double category of relations.
	$\Cau(\dbl{D})$ is a cartesian equipment
	with strong tabulators.
\end{lemma}

\begin{proof}
	Firstly,
	we prove that the canonical pseudo-functor $\dbl{D}\to\Cau(\dbl{D})$ preserves strong tabulators.
	Let $A \sto["R"] B$ be a horizontal arrow and $\langle l, r \rangle\colon \top R \to A\times B$ exhibit $X$ as a tabulator of $R$ in $\dbl{D}$.
	Take an arbitrary cell of the following form in $\Cau(\dbl{D})$.
	\[
		\begin{tikzcd}[column sep=small]
				&
				X
				\ar[ld, "m^*f_!"']
				\ar[rd, "n^*g_!"]
					&
			\\
			A
			\doublecell[rr, shift left=3ex]{\alpha}
			\sar[rr, "R"']
				&
					&
					B	
		\end{tikzcd}
	\]
	Then, we can assume $m$ and $n$ are the same vertical arrows from some object $Y$ to $X$ by \Cref{rem:equaldomain}.
	Moreover, by precomposing the isomorphism $m_!$ on the top of the triangle cell,
	we obtain the cell in the image of the canonical faithful pseudo-functor $\dbl{D}\to\Cau(\dbl{D})$.
	Let $k\colon Y\to X$ be the unique vertical arrow in $\dbl{D}$ such that $k\fatsemi l=f$ and $k\fatsemi r=g$.
	Then, we have the following decomposition of the cell in $\Cau(\dbl{D})$.
	\[ 
	\begin{tikzcd}[column sep=small, row sep=small]
			&
			Y
			\ar[ldd, "f_!"']
			\ar[rdd, "g_!"]
				&
		\\
		\\
		A
		\doublecell[rr, shift left=3ex]{\alpha}
		\sar[rr, "R"']
			&
				&
				B	
	\end{tikzcd}
	\hspace{2ex}
	=
	\hspace{2ex}
	\begin{tikzcd}[column sep=small, row sep=small]
		&
		Y
		\ar[d,"k_!"']
		&
		\\
		&
		\top R 
		\ar[ld,"l_!"']
		\ar[rd,"r_!"]
		&
		\\
		A
		\doublecell[rr, shift left=2ex]{\tau}
		\sar[rr, "R"']
			&
				&
				B
	\end{tikzcd}
	\]
	On the right-hand side, $\tau$ is the cell corresponding to the tabulating cell.
	Suppose there is another map $s^*t_!\colon Y\to \top R$ such that $s\colon Z\to Y$ is a monic cover in $\dbl{D}$ and $s^*t_!\fatsemi l_!=f_!$ and $s^*t_!\fatsemi r_!=g_!$.
	Then, we have the following equality of cells in $\dbl{D}$.
	\[
		\begin{tikzcd}[column sep=small, row sep=small]
			&
			Z
			\ar[ld,"s"']
			\ar[d,"t"']
			\ar[rd,"s"]
			&
			\\
			Y
			\ar[d,"f"']
			&
			\top R 
			\ar[ld,"l"']
			\ar[rd,"r"]
			&
			Y
			\ar[d,"g"]
			\\
			A
			\doublecell[rr, shift left=2ex]{\tab}
			\sar[rr, "R"']
			&
			&
			B
		\end{tikzcd}
		\hspace{2ex}	
		=
		\hspace{2ex}
		\begin{tikzcd}[column sep=small, row sep=small]
			&
			Z
			\ar[d,"s"']
			&
			\\
			&
			Y
			\ar[ld,"f"']
			\ar[rd,"g"]
			&
			\\
			A
			\doublecell[rr, shift left=2ex]{\alpha}
			\sar[rr, "R"']
				&
					&
					B
		\end{tikzcd}
		\hspace{2ex}	
		=
		\hspace{2ex}
		\begin{tikzcd}[column sep=small, row sep=small]
			&
			Z
			\ar[d,"s\fatsemi k"']
			&
			\\
			&
			\top R
			\ar[ld,"l"']
			\ar[rd,"r"]
			&
			\\
			A
			\doublecell[rr, shift left=2ex]{\tab}
			\sar[rr, "R"']
				&
					&
					B
		\end{tikzcd}
	\]
	Therefore, we have $t=s\fatsemi k$, and hence $s^*t_!=k_!$ as a vertical arrow in $\Cau(\dbl{D})$.
	This leads the cell $\tau$ in $\Cau(\dbl{D})$ to be a tabulating cell of $R$.
	
	Secondly, we prove that $\Cau(\dbl{D})$ is cartesian.
	Recall from \cite[Proposition 4.2.3]{Ale18} that 
	an equipment $\dbl{E}$ is cartesian 
	if and only if its vertical category has finite products, 
	its horizontal bicategory has finite products locally,
	and the lax-functors $1\colon \mathbbm{1}\to\dbl{E}$ and $\times\colon \dbl{E}\times\dbl{E}\to\dbl{E}$
	induced by the two kinds of the finite products are pseudo-functors.
	In this case, however, it suffices to show only the first condition.
	This is because
	$\dbl{D}$ is cartesian and
	nothing changes on the horizontal bicategory when we take the Cauchisation.
	Note that the definition for a lax functor to be a pseudo-functor only requires properties describable in the horizontal bicategory. 

	In general, given an equipment whose horizontal bicategory has finite products locally,
	the tabulator of a local terminal $\top\colon A \sto B$ gives the product of $A$ and $B$ in the vertical category.
	In the case of $\Cau(\dbl{D})$, the local terminals and tabulators are the same as those for $\dbl{D}$, and hence we conclude that
	the vertical category $\oneV{\Cau(\dbl{D})}$ has binary products which are the same as those for $\dbl{D}$.
	It remains to show that it has $1$ as its terminal object, but for each object $A$,
	all the maps from $A$ to the terminal object $1$ are displayed as $m^*!_!$ for some monic cover $m$,
	and they are all isomorphic to $!_!$ since $m\fatsemi !=!$.
	This means that $!_!$ is the unique isomorphic class of maps from $A$ to $1$, and $1$ is the terminal object in the vertical category.
\end{proof}

\begin{lemma}
	\label{lem:monoincincau'}
	Let $\dbl{D}$ be a unit-pure double category of relations 
	and $\Cau(\dbl{D})$ be its Cauchisation.
	\begin{enumerate}
		\item Let $m$ be a monic cover and $f$ be a vertical arrow in $\dbl{D}$.
		Then $m^*f_!$ is an inclusion (resp. a cover) in $\Cau(\dbl{D})$
		if and only if $f$ is an inclusion (resp. a cover) in $\dbl{D}$.
		\item A vertical arrow $p\colon A\to B$ in $\Cau(\dbl{D})$ is right orthogonal to all covers
		if and only if there are a vertical arrow $f'\colon C\to B$ in $\dbl{D}$ that is right orthogonal to all covers in $\dbl{D}$
		and an isomorphism $s\colon A\to C$ in $\Cau(\dbl{D})$ such that $p=s\fatsemi f'_!$.
	\end{enumerate}
\end{lemma}

The reader should be aware that even if $m^*f_!$ is right orthogonal to all covers in $\Cau(\dbl{D})$, 
$f$ is not necessarily right orthogonal to all covers in $\dbl{D}$.
The key to the proof of this lemma is that 
every vertical arrow in $\Cau(\dbl{D})$ is written as $f_!$ for some vertical arrow $f$ in $\dbl{D}$ up
to isomorphism.

\begin{proof}
	Above all, observe that for each monic cover $m\colon A\to B$,
	a map $m_!$ is an isomorphism, $m^*$ being the inverse.
	Therefore, to see if $B \sto["m^*"]A \sto["f_!"] Y$ satisfies one of the conditions in $\Cau(\dbl{D})$,
	we only need to check whether $f_!$ does so.
	
	Firstly, we give a proof for the case of inclusions and covers.
	Note that the canonical fibred functor 
	$\left(
		\begin{tikzcd}[row sep=small, ampersand replacement=\&]
			\dbl{D}_1
			\ar[d, "{\langle \src,\tgt\rangle}"']
			\\	
			\dbl{D}_0\times\dbl{D}_0
		\end{tikzcd}
	 \right)
	 \to
	 \left(
		\begin{tikzcd}[row sep=small, ampersand replacement=\&]
			\Cau(\dbl{D})_1
			\ar[d, "{\langle \src,\tgt\rangle}"']
			\\	
			\Cau(\dbl{D})_0\times\Cau(\dbl{D})_0	
		\end{tikzcd}		
				 \right)
	$ 
	given by the pseudo-functor $\dbl{D}\to\Cau(\dbl{D})$ is a base change.
	Since the condition to be an inclusion or a cover is described merely by the canonical fibred category and 
	the horizontal identity arrow, the statement follows.

	Secondly, we move on to the case of arrows right orthogonal to all covers.
	Observe that in a unit-pure equipment,
	two vertical arrows, $f,g \colon A \to B$, are identical if and only if $f_!$ and $g_!$ are identical.
	On one hand, we will show that if $f$ is right orthogonal to all covers, then $f_!$ is right orthogonal to all covers in $\Cau(\dbl{D})$.
	Take a commutative square in $\oneV{\Cau(\dbl{D})}$ as below, where $r$ is a cover in $\Cau(\dbl{D})$.
	\[
		\begin{tikzcd}
			X
				\sar[r,"p"]
				\sar[d,"r"']
			&
				A
				\sar[d,"f_!"]
			\\
				Y
				\sar[r,"q"']
			&
				B
		\end{tikzcd} 
	\]		
	Since orthogonality is invariant under isomorphisms, we assume without loss of generality that $r$ is of form $e_!$ where $e$ is a cover.
	Using composition and pullback, we further assume that $p$ and $q$ are also representable, as below.
	\[
		\begin{tikzcd}
			X
				\sar[r,"h_!"]
				\sar[d,"e_!"']
			&
				A
				\sar[d,"f_!"]
			\\
				Y
				\sar[r,"g_!"']
			&
				B
		\end{tikzcd} 
	\]
	Then, the orthogonality in $\dbl{D}$ of $f$ and $e$ leads to a vertical arrow $k$ 
	with $e\fatsemi k=h$ and $k\fatsemi f=g$ to exist.
	So we have the filler $k_!$.
	Suppose
	we have another filler $n^*k'_!\colon Y\sto A$ where $n$ is a monic cover.
	Then, we have the following commutative square in $\dbl{D}$ by the unit-pure property of $\dbl{D}$.
	\[
		\begin{tikzcd}
			X
				\ar[r,"k'"]
				\ar[d,"n"']
			&
				A
				\ar[d,"f"]
			\\
				Y
				\ar[r,"g"']
				\ar[ru, dotted, "d"'description]
			&
				B
		\end{tikzcd}
	\]
	Then, we have a vertical arrow $d$ as in the diagram.
	$d_!$ is the same as $n^*k'_!$ as a vertical arrow in $\Cau(\dbl{D})$,
	and $d$ is equal to the unique map $k$.
	Thus, the filler in $\oneV{\Cau(\dbl{D})}$ is unique.

	On the other hand, we would like to prove that if $f_!$ is right orthogonal to all covers in $\Cau(\dbl{D})$,
	then $f_!$ has a presentation $n^*f'_!$ where $n$ is a monic cover, and $f'$ is right orthogonal to all covers.
	Take the $(\Cov(\dbl{D}),\Fib(\dbl{D}))$-factorisation of $f$ as $e\fatsemi m$.
	By the assumption, we have the unique filler $n^*g_!$ of the following square.
	\[
		\begin{tikzcd}
			A
				\sar[r,equal]
				\sar[d,"e_!"']
			&
				A
				\sar[d,"f_!"]
			\\
				Y
				\sar[r,"m_!"']
				\sar[ur,dotted,"n^*g_!" description]
			&
				B
		\end{tikzcd}
	\]
	From the bottom right triangle, we have $n\fatsemi m=g\fatsemi f$ since $n^*$ is the inverse of $n_!$.
	Since $m$ is monic, it follows that $n=g\fatsemi e$.
	Combining this with the commutativity of the top left triangle, 
	we deduce that $g_!$ is the inverse of $e_!n^*$, $f_!$ is equal to
	$g^*n_!m_!$, and $g^*n_!$ is an isomorphism.
\end{proof}

\begin{theorem}
	If $\dbl{D}$ is a unit-pure double category of relations,
	so is $\Cau(\dbl{D})$.
\end{theorem}
\begin{proof}
	We have already done with the most part, including \cref{lem:rcovtab}.
	What remains to be proved is that $\Cau(\dbl{D})$ admits a
	$\rclass{\Cov(\Cau(\dbl{D}))}$-comprehension scheme.
	By \cref{lem:monoincincau'}, every map from $C$ to $A\times B$ in $\rclass{\Cov(\Cau(\dbl{D}))}$
	is identified with a map $\langle f, g \rangle_!$ up to precomposition of isomorphisms where $\langle f, g \rangle\colon C'\to A\times B$ is an
	arrow in $\Fib(\dbl{D})$ .
	The pair $f,g$ exhibits $C'$ as the tabulator of $f^*g_!$ by the comprehension scheme of $\dbl{D}$.
	Then, the pair of $f_!$ and $g_!$ exhibits $C'$ as the tabulator of $f^*g_!$ in $\Cau(\dbl{D})$
	by a similar argument as in the proof of \cref{lem:monoincincau'}.
\end{proof}

\begin{example}
	Let $\one{E}$ be a quasi-topos and $\Rel{\Epi,\Regmono}(\one{E})$ 
	be the double category of relations defined by the stable factorisation system
	$(\Epi,\Regmono)$ on $\one{E}$.
	Horizontal arrows in $\Rel{\Epi,\Regmono}(\one{E})$ are strong relations,
	meaning they are subobjects of the product of two objects, which
	are characterised by the weak subobject classifier $\Omega$.
	The Cauchisation of this equipment is equivalent to the double category of relations
	defined by the topos $\Cs(\one{E})$ of coarse objects.
	See \cite[A 2.6]{Joh02} for the definition of coarse objects and the discussion on them.
\end{example}

\subsection{Contrasting the literature --- relations on regular categories and spans}
\label{sec:Contrast}
In this section, we contrast our work with some characterisation theorems in the literature;
one of the bicategories of relations on regular categories \cite{CW87},
one of the bicategories of spans \cite{LWW10},
and one of the double categories of spans \cite{Ale18}.

Before we begin the discussion, we recall the notion of \textit{cartesian bicategories}.
In short, a cartesian bicategory is a bicategory $\bi{B}$ satisfying the following. See \cite{CKWW07} for more detail.
\begin{itemize}
	\item%
		The subbicategory $\bi{M}$ of maps has finite products.
	\item%
		Each its hom-category has finite products.
	\item%
		A certain derived tensor product on $\bi{B}$, extending the product structure of $\bi{M}$, is functorial.
\end{itemize}
The definition of cartesian double categories \cref{def:cartesian} 
is a double-categorical analogue of this notion introduced in \cite{Ale18}.
Note that, when we construct a cartesian double category from a cartesian bicategory through $\Map$, 
the resulting double category is inevitably Cauchy.

First, we reconstruct the classical characterisation theorem \cite[Theorem 3.5]{CW87} for relations on regular categories.
Let us introduce the concept central to their characterisation, rephrased in terms of double categories.
\begin{definition}[{\cite[Definition 3.1]{CW87},}]
	A cartesian equipment $\dbl{D}$ is \emph{functionally complete} if, for each horizontal arrow of the form $r\colon X\sto 1$,
	there exist an object $X_r$, an inclusion $i\colon X_r\to X$, and an opcartesian cell of the following form.
	\[
		\begin{tikzcd}[column sep=small]
				&
				X_r
				\ar[ld, "i"']
				\ar[rd, "!"]
					&
			\\
			X
			\doublecell[rr, shift left=2ex]{\opcart}
			\sar[rr, "r"']
				&
					&
					1
		\end{tikzcd}
	\]
\end{definition}
\begin{remark} 
	This definition differs considerably from the one defined in \cite[Definition 4.8]{Lam22}.
	Therein, the term `functionally complete' is used as having strong $\Mono$-tabulators in our terminology.
\end{remark}
In \cite{CW87}, it is shown that the bicategory of relations arising from a regular category is characterised as a
locally posetal cartesian bicategory which is functionally complete and such that objects are discrete. 
Therefore, we obtain the following characterisation of $\Rel{\Regepi, \Mono}(\one{C})$.
Note that a locally posetal cartesian bicategory is essentially the same as a locally posetal and Cauchy cartesian equipment
\begin{theorem}[{\cite[Theorem 3.5]{CW87}}]
	\label{thm:RegcatCW}
	The following are equivalent for a double category $\dbl{D}$.
	\begin{enumerate}
		\item $\dbl{D}$ is equivalent 
		to $\Rel{\Regepi,\Mono}(\one{C})$ for some regular category $\one{C}$.
		\item $\dbl{D}$ is a locally posetal, discrete Cauchy cartesian equipment that is functionally complete.
	\end{enumerate}
\end{theorem}

However, now that we have more general characterisation theorems for double categories of relations,
the essential part of the argument in \cite{CW87} can be extracted as the following lemma.
\begin{lemma}
	Let $\dbl{D}$ be a locally posetal, discrete, Cauchy cartesian equipment and $\alpha$ be a cell of the following form.
	Suppose, moreover, that $l$ is an inclusion.
	\[
		\begin{tikzcd}[column sep=small]
				&
				A
				\ar[ld, "l"']
				\ar[rd, "!"]
					&
			\\
			X
			\doublecell[rr, shift left=3ex]{\alpha}
			\sar[rr, "p"']
				&
					&
					1
		\end{tikzcd}
	\]
	Then $\alpha$ is tabulating if it is opcartesian.
\end{lemma}
\begin{proof}
	Note that by \cref{lem:DagLocPos}, a locally posetal discrete cartesian equipment is unit-pure.
	Take another triangle $\beta$ of the following form, and
	we show there exists a unique $k\colon B\to A$ that, composited with $\alpha$, gives $\beta$.
	Since $l$ is monic by \cref{lem:covepi}, the uniqueness part is trivial,
	and all we have to check is $k\fatsemi l=f$ since $\dbl{D}$ is locally posetal.

	Define $u\colon B\sto A$ and $v\colon A\sto B$ as the following restrictions.
	\[
		\begin{tikzcd}[column sep=small]
				&
				B
				\ar[ld, "f"']
				\ar[rd, "!"]
					&
			\\
			X
			\doublecell[rr, shift left=3ex]{\beta}
			\sar[rr, "p"']
				&
					&
					1
		\end{tikzcd}
		\hspace{4ex}
		,
		\hspace{4ex}
		\begin{tikzcd}[column sep=small]
			B
			\sar[rr, "u"]
			\doublecell[rr, shift right=2.5ex]{\cart}
			\ar[rd, "f"']
				&
					&
					A
					\ar[ld, "l"]
			\\
				&
				X
					&
		\end{tikzcd}
		\begin{tikzcd}[column sep=small]
			A
			\sar[rr, "v"]
			\doublecell[rr, shift right=2.5ex]{\cart}
			\ar[rd, "l"']
				&
					&
					B
					\ar[ld, "f"]
			\\
				&
				X
					&
		\end{tikzcd}
	\]
	Take an opcartesian cell $\omega$ of the form described below and define $\bar\alpha'$ through the correspondence
	obtained in \cref{lem:UnbiasingLemma}.
	\[
		\begin{tikzcd}[column sep=small]
				&
				A
				\ar[ld, "l"']
				\ar[rd, "!"]
					&
			\\
			X
			\doublecell[rr, shift left=3ex]{\alpha}
			\sar[rr, "p"']
			\ar[d , "\Delta"']
				&
					&
					1
					\ar[d, equal]
			\\
			X\times X
			\doublecell[rr, shift left=3.5ex]{\omega}
			\sar[rr, "p'"']
				&
					&
					1
		\end{tikzcd}
		\hspace{2ex}
		\vline
		\,
		\vline
		\hspace{2ex}
		\begin{tikzcd}[column sep=small]
				&
				A
				\ar[dl, "l"']
				\ar[dr, "l"]
			\\
			X
			\doublecell[rr, shift left=3ex]{\bar\alpha'}
			\sar[rr, "\bar{p}'"']
				&
					&
					X
		\end{tikzcd}
	\]
	Since $\alpha$ and $\omega$ are opcartesian, \cref{lem:UnbiasingLemma} shows $\bar\alpha'$ is also opcartesian.

	Define, moreover, the following two horizontal cells $\psi$ and $\varphi$ as follows.
	\[
		\begin{tikzcd}[row sep=small]
				&
				B
				\sar[dr, "u"]
				\ar[ddd, "f"description]
				\doublecell[ddd, yshift=1ex,shift right=4.0ex]{\cart}
				\doublecell[ddd, yshift=1ex,shift left=4.0ex]{\cart}
					&
			\\
			A
			\sar[ur, "v"]
			\ar[rdd, "l"']
				&
					&
					A
					\ar[ldd, "l"]
			\\
			\\
				&
				X
					&
		\end{tikzcd}
		=
		\begin{tikzcd}[row sep=small]
				&
				B
				\sar[dr, "u"]
					&
			\\
			A
			\sar[ur, "v"]
			\sar[rr,equal]
			\ar[rdd, "l"']
			\doublecell[rr,  shift left=2ex]{\psi}
			\doublecell[rr, shift right=3ex]{\cart}
				&
					&
					A
					\ar[ldd, "l"]
			\\
			\\
				&
				X
					&
		\end{tikzcd}
	\]
	\[
		\begin{tikzcd}[column sep=small]
				&
				B
				\ar[ld, "f"']
				\ar[rd, "!"]
					&
			\\
			X
			\doublecell[rr, shift left=3ex]{\beta}
			\sar[rr, "p"']
			\ar[d , "\Delta"']
				&
					&
					1
					\ar[d, equal]
			\\
			X\times X
			\doublecell[rr, shift left=3ex]{\omega}
			\sar[rr, "p'"']
				&
					&
					1
		\end{tikzcd}
		\hspace{2ex}
		\vline
		\,
		\vline
		\hspace{2ex}
		\begin{tikzcd}[column sep=small]
				&
					&
					B
					\ar[lld, equal]
					\ar[rrd, equal]
						&
							&
			\\
			B
			\sar[rr, "u"{description, near end}]
			\ar[dr, "f"']
			\doublecell[rrrr, shift left=3ex]{\varphi}
			\doublecell[rr, shift right=3ex]{\cart}
				&
					&
					A
					\sar[rr, "v"{description, near start}]
					\ar[dl, "l"description]
					\ar[dr, "l"description]
					\doublecell[rr, shift right=3ex]{\cart}
						&
							&
							B
							\ar[dl, "f"]
			\\
				&
				X
				\doublecell[rr, shift left=3ex]{\bar\alpha'}
				\sar[rr, "\bar{p}'"']
					&
						&
						X
							&
		\end{tikzcd}
	\]
	Here, the last bijective correspondence is that observed in \cref{lem:UnbiasingLemma}, and
	the existence of $\varphi$ follows from the fact that $\bar\alpha'$ is opcartesian,
	considering \cref{lem:Sandwich}. This shows $u$ is left adjoint to $v$; hence, $u$ is represented by a vertical arrow $k$, which satisfies
	$k\fatsemi l=f$ by the unit-pureness of $\dbl{D}$.
\end{proof}

\begin{proof}[Alternative proof of \cref{thm:RegcatCW}] 
	It is clear that the double category of the form $\Rel{\Regepi,\Mono}(\one{C})$ satisfies the second condition of \cref{thm:RegcatCW}.
	We show the converse.
	For an equipment $\dbl{D}$ satisfying the condition of \Cref{lem:forRegcatCW},
	inclusions in $\dbl{D}$ and monomorphisms in $\oneV{\dbl{D}}$ coincide by \Cref{lem:monoinc,lem:DagLocPos}.
	Using \cref{lem:monoinc} and this lemma,
	$\dbl{D}$ admitting a left-sided $\Mono$-comprehension scheme follows from discreteness and functionally completeness,
	which leads to the characterisation of $\dbl{D}$ as $\Rel{\Regepi,\Mono}(\one{C})$ for some regular category $\one{C}$ by \Cref{thm:unitpuredbl}.
\end{proof}

Second, we adapt our characterisation theorem to the double category of spans in accordance with the previous work \cite{LWW10,Ale18}.
To begin with, we introduce the double categorical analogue of the notions appeared in \cite{LWW10}.
Our terminology is slightly different from that of \cite{Ale18}.

\begin{definition}
	\label{defn:comonoid}
	Let $\dbl{D}$ be a double category.
	A \emph{comonoid} in $\dbl{D}$ is a comonad in $\biH{\dbl{D}}$. In other words,
	a comonoid is a tuple $G=(A,G,\delta_G,\alpha_G)$ such that $A$ is an object in $\dbl{D}$, $G$ is a horizontal endoarrow on $A$, and
	$\delta_G$ and $\alpha_G$ are cells in $\dbl{D}$ satisfying the following equalities.
	\[
		\begin{tikzcd}[column sep=small]
				&
				A
				\sar[rr, "G"]
				\ar[ld, equal]
				\doublecell[rr,shift right=3ex]{\delta_G}
					&
						&
						A
						\ar[rd, equal]
							&
			\\
			A
			\sar[rr, "G"]
			\ar[rd, equal]
			\doublecell[rr, shift right=3ex]{\alpha_G}
				&
					&
					A
					\ar[ld, equal]
					\sar[rr, "G"]
					\doublecell[rd]{\verteq}
						&
							&
							A
							\ar[ld, equal]
			\\
				&
				A
				\sar[rr, "G"']
					&
						&
						A
							&
		\end{tikzcd}
		\hspace{2ex}
		=
		\hspace{2ex}
		\begin{tikzcd}[column sep=small]
			A
			\sar[rr, "G"]
			\doublecell[rr, shift right=3.5ex]{\verteq}
			\ar[d, equal]
				&
					&
					A
					\ar[d, equal]
			\\
			A
			\sar[rr, "G"']
				&
					&
					A
		\end{tikzcd}
		\hspace{2ex}
		=
		\hspace{2ex}
		\begin{tikzcd}[column sep=small]
				&
				A
				\sar[rr, "G"]
				\ar[ld, equal]
				\doublecell[rr, shift right=3ex]{\delta_G}
					&
						&
						A
						\ar[rd, equal]
							&
			\\
			A
			\sar[rr, "G"]
			\ar[rd, equal]
				&
					&
					A
					\doublecell[ld]{\verteq}
					\sar[rr, "G"]
					\ar[rd, equal]
					\doublecell[rr, shift right=3ex]{\alpha_G}
						&
							&
							A
							\ar[ld, equal]
			\\
				&
				A
				\sar[rr, "G"']
					&
						&
						A
							&
		\end{tikzcd}
	\]
	\[
		\begin{tikzcd}
			A
			\sar[rrr,"G"]
			\ar[d,equal]
			\doublecell[rrrd]{\delta_G}
				&
					&
						&
						A
						\ar[d,equal]
			\\
			A
			\sar[rr,"G"]
			\ar[d,equal]
			\doublecell[rrd]{\delta_G}
				&
					&
					A
					\sar[r]
					\ar[d,equal]
					\doublecell[rd]{\verteq}
						&
						A
						\ar[d,equal]
			\\
			A
			\sar[r,"G"']
				&
				A
				\sar[r,"G"']
					&
					A
					\sar[r,"G"']
						&
						A
		\end{tikzcd}
		\hspace{2ex}
		=
		\hspace{2ex}
		\begin{tikzcd}
			A
			\sar[rrr,"G"]
			\ar[d,equal]
			\doublecell[rrrd]{\delta_G}
				&
					&
						&
						A
						\ar[d,equal]
			\\
			A
			\sar[r,"G"]
			\ar[d,equal]
			\doublecell[rd]{\verteq}
				&
				A
				\sar[rr]
				\ar[d,equal]
				\doublecell[rrd]{\delta_G}
					&
						&
						A
						\ar[d,equal]
			\\
			A
			\sar[r,"G"']
				&
				A
				\sar[r,"G"']
					&
					A
					\sar[r,"G"']
						&
						A
		\end{tikzcd}
	\]
	More formally, a comonoid in $\dbl{D}$ is a strictly normal pseudo double functor into $\dbl{D}$ from the `walking' comonoid $\mathbbm{c}\mathrm{omnd}$,
	which is a double category constructed from the monoidal category $\Delta^\op_a$ consisting of finite ordinals and monotone maps.

	Note that if we have an adjoint $P\colon A\adjointleft B\lon Q$ in $\biH{\dbl{D}}$, we obtain a comonoid $QP$ on $B$.
	In particular, if $\dbl{D}$ is an equipment, we obtain a comonoid $f^*f_!$ for each vertical arrow $f\colon A \to B$.
\end{definition}

\begin{definition}\label{defn:comodule}
	Let $G=(A, G,\delta_G,\alpha_G)$ be a comonoid.
	A \emph{comodule}%
	\footnote{%
		If we follow the terminology of \cite{CS10},
		a comodule of $G$ can be seen as an object $X$ equipped with a \textit{comonoid homomorphism} from $\Id_X$ to $G$, or vertical arrow in $\Mod(\dbl{D}^\vop)$
		from $G$ to $\Id_X$,
		where $\Id_X$ is seen as the comonoid on $X$ induced from the coherence cells for this horizontal identity.
	}
	of $G$ is a pair $(f,\nu)$ of a vertical arrow $f$ and a cell $\nu$ satisfying the following equalities.
	\begin{equation}\label{eq:comodule}
		\begin{tikzcd}[column sep=small]
				&
					&
					B
					\ar[d, "f"]
					\ar[lld, "f"', bend right]
					\ar[rrd, "f", bend left]
						&
							&
			\\
			A
			\sar[rr, "G"']
				&
					&
					A
					\doublecell[rr, shift left=3ex]{\nu}
					\doublecell[ll, shift right=3ex]{\nu}
					\sar[rr, "G"']
						&
							&
							A
		\end{tikzcd}
		=
		\begin{tikzcd}[column sep=small]
				&
					&
					B
					\ar[ld, "f"']
					\ar[rd, "f"]
						&
							&
			\\
				&
				A
				\sar[rr, "G"]
				\ar[ld, equal]
				\doublecell[rr, shift right=3ex]{\delta_G}
				\doublecell[rr, shift left=3ex]{\nu}
					&
						&
						A
						\ar[rd, equal]
							&
			\\
			A
			\sar[rr, "G"']
				&
					&
					A
					\sar[rr, "G"']
						&
							&
							A
		\end{tikzcd}\text{,}
		\hspace{5ex}
		\begin{tikzcd}[column sep=small]
				&
				B
				\ar[ld, "f"']
				\ar[rd, "f"]
					&
			\\
			A
			\sar[rr, "G"']
			\ar[rd, equal]
			\doublecell[rr, shift right=3ex]{\alpha_G}
			\doublecell[rr, shift left=3ex]{\nu}
				&
					&
					A
					\ar[ld, equal]
			\\
				&
				A
					&
		\end{tikzcd}
		=
		\begin{tikzcd}[column sep=small]
				&
				B
				\ar[d, bend right=40, "f"']
				\doublecell[d]{\Id}
				\ar[d, bend left=40, "f"]
					&
			\\
				&
				A
					&
		\end{tikzcd}
	\end{equation}
\end{definition}
\begin{definition}\label{def:coEMcomonoid}
	Let $\dbl{D}$ be a double category.
	For a comonoid $G=(G,\delta_G,\alpha_G)$,
	a \emph{co-Eilenberg-Moore object} $\coEM{G}=(\coEM{G}, u,\upsilon)$
	is the one-dimensional universal comodule of $G$; i.e., 
	there is a bijective correspondence between the following data,
	obtained through postcomposing $(u,\upsilon)$ to $h$.
	\[
		\text{a comodule of $G$}
		\begin{tikzcd}[column sep=small]
				&
				X
				\ar[ld, "g"']
				\ar[rd, "g"]
					&
			\\
			A
			\sar[rr, "G"']
			\doublecell[rr, shift left=3ex]{\mu}
				&
					&
					A
		\end{tikzcd}
		\vline
		\,
		\vline
		\text{ a vertical arrow}
		\begin{tikzcd}
			X
			\ar[d, "h"]
			\\
			\coEM{G}
		\end{tikzcd}
	\]

	In other words, a co-Eilenberg-Moore object for a comonoid $(G,\delta_G,\alpha_G)$
	is the one-dimensional double limit
	in the sense of \cite{GP99}
	of the diagram (= functor from $\mathbbm{m}\mathrm{nd}$) corresponding to $G$.

	A double category $\dbl{D}$ has \emph{co-Eilenberg-Moore objects for comonoids} if
	for each comonoid $G$ in $\dbl{D}$, there exists a co-Eilenberg-Moore object for $G$.
	A co-Eilenberg-Moore object $(\coEM{G},u,\upsilon)$ of a comonoid $G$
	is \emph{strong} if $\upsilon$ is an opcartesian cell, and it is an \emph{$\zero{M}$-co-Eilenberg-Moore object} for 
	a class $\zero{M}$ of vertical arrows of $\dbl{D}$
	if $u$ belongs to $\zero{M}$.
	We say $\dbl{D}$ has \emph{(strong/$\zero{M}$-) co-Eilenberg-Moore objects for comonoids} if 
	every comonoid has a (strong/$\zero{M}$-) co-Eilenberg-Moore object.
\end{definition}
	The notion of a co-Eilenberg-Moore object coincides with the notion of a \textit{collapse} in the vertical opposite $\dbl{D}^\vop$ seen as a virtual double category
	in the sense of \cite{Sch15}.
\begin{definition}\label{def:comonoidic}
	Let $\dbl{D}$ be an equipment and $f\colon A\to B$ be a vertical arrow in $\dbl{D}$.
	As mentioned in \cref{defn:comonoid}, $f^*f_!$ is a comonoid, and one can easily check that the canonical opcartesian cell
	\begin{equation}\label{eq:comonoidic}
		\begin{tikzcd}[column sep=small]
				&
				A
				\ar[ld, "f"']
				\ar[rd, "f"]
					&
			\\
			B
			\sar[rr, "f^*f_!"']
			\doublecell[rr, shift left=2.5ex]{\opcart}
				&
					&
					B
		\end{tikzcd}
	\end{equation}
	exhibits $A$ as a comodule of this comonoid.
	We say $f$ is \emph{comonoidic} if
	this comodule
	is a co-Eilenberg-Moore object of the comonoid $f^*f_!$.
\end{definition}

The paper \cite{Ale18} puts more emphasis on \textit{copointed endomorphisms} than comonoids (comonads in the paper)
in the characterisation of double categories of spans on a finitely
complete category.
They are a loosened version of comonoids, lacking the comultiplications.
The counterparts of \cref{defn:comonoid,defn:comodule,def:coEMcomonoid,def:comonoidic}
for this notion are defined as follows.

\begin{definition}[cf. {\cite[\S 3]{Ale18}}]
	Let $\dbl{D}$ be a double category.
	\begin{enumerate} 
		\item%
			A \emph{copointal} of $\dbl{D}$ 
			is a copointed endomorphism.
			In other words, a copointal is a triple $p=(X, p, \alpha_p)$ consisting of
			an object $X$ of $\dbl{D}$, a horizontal arrow $p\colon X\sto X$ of $\dbl{D}$, and
			a horizontal 2-cell $\alpha_p\colon p\Rightarrow\Id_X$ of $\dbl{D}$.
		\item%
			A \emph{comodule (for copointal)} of a copointal $p$ is a pair $(f,\nu)$ of a vertical arrow $f$ and a cell $\nu$ satisfying 
			the second equality of \Cref{eq:comodule} for $\alpha_p$ instead of $\alpha_G$.
		\item%
			For a copointal $p$, a \emph{co-Eilenberg-Moore object} $\coEMp{p}=(\coEMp{p}, u,\upsilon)$ 
			is the one-dimensional universal comodule of $p$
			in the same sense as \cref{def:coEMcomonoid}.
			A co-Eilenberg-Moore object of a copointal $p$, $(\coEMp{p},u,\upsilon)$,
			is \emph{strong} if $\upsilon$ is an opcartesian cell, and it is an \emph{$\zero{M}$-co-Eilenberg-Moore object} for 
			a class $\zero{M}$ of vertical arrows of $\dbl{D}$
			if $u$ belongs to $\zero{M}$.
			We say $\dbl{D}$ has \emph{(strong/$\zero{M}$-) co-Eilenberg-Moore objects for copointals} if 
			every copointal has a (strong/$\zero{M}$-) co-Eilenberg-Moore object.
		\item%
			A vertical arrow $f$ is \emph{copointalic} if the comodule \cref{eq:comonoidic} of $f^*f_!$ seen as a copointal 
			is a co-Eilenberg-Moore object for this copointal.
			Note that if every co-Eilenberg-Moore object is strong, then each leg of a co-Eilenberg-Moore object is copointalic.
			\vspace{-2ex}
	\end{enumerate}
\end{definition}
In general, the two notions of co-Eilenberg-Moore objects are not equivalent.
However, when $\dbl{D}$ is close enough to double categories of relations, they essentially coincide.
This observation was made in \cite{LWW10} in the context of cartesian bicategories,
and we revisit it in the context of double categories.

\begin{remark}
	\label{rem:IdSubterminal}
	If an object $A$ in a cartesian equipment $\dbl{D}$ is discrete,
	then in particular, the diagonal $\Delta\colon A \to A\times A$ is an inclusion. Since the horizontal identity on $A\times A$ is
	isomorphic to $\Id_A\times\Id_A$ in a canonical way, this shows $\Id_A\land\Id_A\cong\Id_A$ in the hom-category $\biH{\dbl{D}}(A, A)$,
	which means $\Id_A$ is a subterminal object in this cartesian category.
	Therefore, in a discrete cartesian equipment, there exists at most one copoint for each horizontal endoarrow.
	Moreover, for each copointal $p\colon A \sto A$, the projections $p\land\Id_A\to p$
	and $\Id_A\land p\to p$ are invertible in $\biH{\dbl{D}}(A, A)$. 
	In other words, we have the following two cartesian cells.
	\[
		\begin{tikzcd}
			A
			\sar[r, "p"]
			\ar[d, "\Delta"']
			\doublecell[rd]{\cart}
				&
				A
				\ar[d, "\Delta"]
			\\
			A\times A
			\sar[r, "p\times\Id_A"']
				&
				A\times A
		\end{tikzcd}
		\hspace{2ex}
		,
		\hspace{2ex}
		\begin{tikzcd}
			A
			\sar[r, "p"]
			\ar[d, "\Delta"']
			\doublecell[rd]{\cart}
				&
				A
				\ar[d, "\Delta"]
			\\
			A\times A
			\sar[r, "\Id_A\times p"']
				&
				A\times A
		\end{tikzcd}
		\vspace{-3ex}
	\]
\end{remark}

\begin{proposition}[{cf.\ \cite[Theorem 3.16]{LWW10}}]
	\label{prop:odotISland}
	Given a discrete cartesian equipment $\dbl{D}$
	and copointals $p,q\colon A \sto A$ in $\dbl{D}$,
	the span $p\to[leftarrow]p\odot q\to q$ in $\biH{\dbl{D}}(A,A)$ induced by the copoints of $p$ and $q$ is the product diagram.
\end{proposition}
\begin{proof}
	Consider the following cartesian and opcartesian cells.
	\[
		\begin{tikzcd}
			&
			\cdot
			\sar[rr, "p"]
			\ar[ld, "\Delta"']
			\doublecell[rd]{\cart}
				&
					&
					\cdot
					\ar[ld, "\Delta"description]
					\ar[rd, "\Delta"description]
						&
							&
							\cdot
							\sal[ll, "q"']
							\ar[rd, "\Delta"]
							\doublecell[ld]{\cart}
								&
		\\
		\cdot
		\sar[rr, "p\times\Id"]
		\ar[rd, equal]
		\doublecell[rrrd]{\opcart}
			&
				&
				\cdot
				\sar[rr]
				\ar[rd, "\id\times\Delta"description]
				\doublecell[rr, shift left=2.5ex]{\opcart}
				\doublecell[rr, shift right=3ex, yshift=1.5ex]{\cart}
					&
						&
						\cdot
						\ar[ld, "\Delta\times\id"description]
							&
								&
								\cdot
								\sal[ll, "\Id\times q"']
								\ar[ld, equal]
								\doublecell[llld]{\opcart}
		\\
			&
			\cdot
			\sar[rr, "p\times\Delta_!"']
			\ar[d, "\id\times\Delta"']
			\doublecell[rrd]{\cart}
				&
					&
					\cdot
					\ar[d, equal]
						&
							&
							\cdot
							\sal[ll, "\Delta^*\times q"]
							\ar[d, "\Delta\times\id"]
							\doublecell[lld]{\cart}
								&
		\\
			&
			\cdot
			\sar[rr, "p\times\Id\times\Id"']
				&
					&
					\cdot
						&
							&
							\cdot
							\sal[ll, "\Id\times\Id\times q"]
								&
	\end{tikzcd}
	\]
	The two cartesian cells on the top row are the same as the cells in \cref{rem:IdSubterminal}.
	The opcartesian cells on the middle row and the cartesian cells on the bottom row are
	obtained by applying to the conjoint and companion of $\Delta$ the pseudo-double functor $-\times-$.
	\cref{lem:Sandwich} shows $p\land \Id\land q\cong p\odot q$.
	Since $p$ is copointed, the projection $p\land\Id\land q\to p\land q$ is invertible. Now it suffices to show that
	projections $p\land\Id\land q \to p$ and $p\land\Id\land q\to q$ are the same as $p\odot q\to p$ and $p\odot q \to q$
	induced from copoints of $p$ and $q$, through the above invertible horizontal cell.

	We only check for $p\land\Id\land q \to p$, while the other is shown similarly.
	Since the local projections for $p\land \Id \land q$
	are obtained by postcomposing projection cells to the cartesian cell defining this product in $\biH{\dbl{D}}$,
	the naturality of the projection cells shows that
	the local projection $p\land\Id\land q\to q$ is obtained as the composite of the cells in the following diagram.
	Here, by $\pi_1$, we mean the first projection for the structure of the cartesian double category of $\dbl{D}$,
	and by $\Delta_3$, we mean the ternary diagonal $A \to A \times A \times A$.
	The two cells $\zeta$ and $\xi$ are cartesian cells in the above diagram dividing the whole cartesian cell.
	Note that since $\Delta$ is an inclusion, $\Delta\times\id$ and $\id\times\Delta$ are inclusions as well, and hence
	the opcartesian cells in the above diagram are all cartesian at the same time.
	\[
		\begin{tikzcd}[column sep=small]
			A
			\sar[rr, "p"]
			\ar[d, "\Delta_3"']
			\doublecell[rrd]{\xi}
				&
					&
					A
					\sar[rr, "q"]
					\ar[d, "\Delta_3"description]
					\doublecell[rrd]{\zeta}
						&
							&
							A
							\ar[d, "\Delta_3"]
			\\
			A\times A \times A
			\sar[rr, "p\times\Id\times\Id"']
			\ar[rd, "\pi_1"']
			\doublecell[rrrd]{\pi_1}
				&
					&
					A\times A \times A
					\sar[rr, "\Id\times\Id\times q"']
					\ar[rd, "\pi_1"']
					\doublecell[rr, shift right=4ex]{\pi_1}
						&
							&
							A\times A \times A
							\ar[ld, "\pi_1"]
			\\
				&
				A
				\sar[rr, "p"']
					&
						&
						A
							&
		\end{tikzcd}
	\]
	Since $\xi$ is cartesian, it gives rise to the isomorphism $p\cong p\land\Id\land\Id$, hence $\xi\fatsemi\pi_1$ is the vertical identity.
	On the other hand, $\zeta\fatsemi\pi_1$ on the right of the above diagram gives the copoint. Therefore, this projection is the same as
	$p\odot q\to p$ induced from the copoint of $p$, since copoints are unique as we mentioned in \cref{rem:IdSubterminal}.
\end{proof}

\begin{corollary}
	\label{cor:copointalcomonoid}
	In a discrete cartesian equipment $\dbl{D}$,
	for each copointal $p\colon A \sto A$,
	there exist unique $\delta\colon p \Rightarrow p\odot p$ and $\varepsilon\colon p \Rightarrow \Id_A$
	such that $(p,\delta,\varepsilon)$ is a comonoid.
\end{corollary}
\begin{proof} 
	By \cref{prop:odotISland}, the cartesian product and the monoidal product $\odot$ on $\biH{\dbl{D}}(A,A)$ coincide for copointals.
	Since every object has a unique comonoid structure in any cartesian monoidal structure,
	for each horizontal endoarrow $p\colon A \sto A$,
	underlying a comonad is a property that is equivalent to being copointed.
\end{proof}
\begin{corollary}
	\label{cor:copointalcomodule}
	In a discrete cartesian equipment $\dbl{D}$,
	for each comonoid $p\colon A \sto A$,
	any pair $(f,\nu)$ of a vertical arrow and a cell of the following form
	\[
		\begin{tikzcd}[column sep=small]
				&
				B
				\ar[ld, "f"']
				\ar[rd, "f"]
					&
			\\
			A
			\sar[rr, "p"']
			\doublecell[rr, shift left=3ex]{\nu}
				&
					&
					A
		\end{tikzcd}
	\]
	is a comodule for $p$.
	In particular, the two notions of a co-Eilenberg-Moore object for copointed arrows and comonoids coincide.
\end{corollary}
\begin{proof}
	Since horizontal identities are subterminal objects in hom-categories by \cref{rem:IdSubterminal},
	the second equation in \cref{defn:comodule} for $(f,\nu)$ to be a comodule automatically holds.
	Observe that through the opcartesian cell defining the comonoid $f^*f_!$,
	$\nu$ corresponds to a horizontal cell $\bar\nu\colon f^*f_!\Rightarrow p$.
	The first equation in \cref{defn:comodule} is equivalent to saying that
	$\bar\nu$ respects the comultiplications in the monoidal category $\biH{\dbl{D}}(A,A)$.
	However, by \cref{prop:odotISland,cor:copointalcomonoid}, this is equivalent to saying that $\bar\nu$ respects the codiagonals, which is trivially true.
\end{proof}

Moving on to the characterisation theorem in terms of co-Eilenberg-Moore objects,
we need the following notion.

\begin{definition}
	Let $\dbl{D}$ be an equipment and $\zero{M}$ be a class of vertical arrows of $\dbl{D}$.
	We say $\dbl{D}$ admits a \emph{unary $\zero{M}$-comprehension scheme w.r.t. copointals}
	if $\dbl{D}$ has strong $\zero{M}$-co-Eilenberg-Moore objects for copointals,
	and every vertical arrow $f\colon A\to B$ in $\zero{M}$ is copointalic.
	Similarly, we say $\dbl{D}$ admits a \emph{unary $\zero{M}$-comprehension scheme w.r.t. comonoids}
	if $\dbl{D}$ has strong $\zero{M}$-co-Eilenberg-Moore objects for comonoids,
	and every vertical arrow $f\colon A\to B$ in $\zero{M}$ is comonoidic.
\end{definition}

The next few pages will be devoted to examining the connection between tabulators and co-Eilenberg-Moore objects,
following the discussions in \cite{Ale18}.

\begin{remark}
	\label{rem:coEMtab}
	For a copointal $p=(A,p,\alpha_p)$ in a unit-pure equipment $\dbl{D}$,
	a cell in the following form
	\[
	\begin{tikzcd}[column sep=small]
		&
		B 
		\ar[ld, "f"']
		\ar[rd, "g"]
			&
	\\
		A
		\sar[rr, "p"']
		\doublecell[rr, shift left=2.5ex]{\nu}
			&
				&
				A
	\end{tikzcd}	
	\]
	must be a comodule of $p$ 
	since postcomposing $\nu$ to $\alpha_p$ gives the identity cell by the unit-pureness.
	In particular, a co-Eilenberg-Moore object of $p$ is the same as a tabulator of $p$.
\end{remark}

\begin{proposition}[cf. {\cite[Corollary 5.1.9]{Ale18}}]\label{prop:coEMtoTabAl}
	Let $\dbl{D}$ be a unit-pure cartesian equipment.
	If $\dbl{D}$ has co-Eilenberg-Moore objects, 
	then, $\dbl{D}$ has tabulators.
	If the co-Eilenberg-Moore objects are strong or $\zero{M}$-co-Eilenberg-Moore objects,
	then, the tabulators are strong or $\zero{M}$-tabulators, respectively.
\end{proposition}
\begin{proof}
	Take a horizontal arrow $p\colon A\sto B$ of $\dbl{D}$ 
	and let $\wh{p}=(A\times B,\wh{p},\alpha_{\wh{p}})$ be a copointal defined by 
	the restriction as follows.
	\[
		\begin{tikzcd}[column sep=small]
			A\times B
			\sar[rr, "\wh{p}"]
			\ar[d, "\Delta\times\id_B"']
			\doublecell[rrd]{\cart}
				&&
				A\times B
				\ar[d, "\id_A\times\Delta"]
			\\
			A\times A\times B
			\sar[rr, "\Id_A\times p\times\Id_B"']
			\ar[rd, "\pi_{13}"']
			\doublecell[rr, shift right=3.5ex]{\id\times!\times\id}
				&&
			A\times B\times B
			\ar[ld, "\pi_{13}"]
			\\
			&
			A\times B
		\end{tikzcd}
		\hspace{2ex}
		=
		\hspace{2ex}
		\begin{tikzcd}[column sep=small]
			A\times B
			\sar[rr, "\wh{p}"]
			\ar[dr,equal]
			\doublecell[rr, shift right=3ex]{\alpha_{\wh{p}}}
				&&
				A\times B
				\ar[dl,equal]
			\\
				&
				A\times B
					&	
		\end{tikzcd}
	\]
	We have the following sequence of correspondences of cells.
	\[
	\begin{tikzcd}[column sep=small]
		&
		X 
		\ar[ld, "{\langle f, g\rangle}"']
		\ar[rd, "{\langle f, g\rangle}"]
			&
	\\
		A\times B
		\sar[rr, "\wh{p}"']
		\doublecell[rr, shift left=2.5ex]{\xi_1}
			&
				&
				A\times B	
	\end{tikzcd}
	\vline
	\,
	\vline
	\begin{tikzcd}[column sep=small]
		&
		X 
		\ar[ld, "{\langle f,f,g\rangle}"']
		\ar[rd, "{\langle f,f,g\rangle}"]
			&
	\\
		A\times A\times B
		\sar[rr, "\Id_A\times p\times\Id_B"']
		\doublecell[rr, shift left=2.5ex]{\xi_2}
			&
				&
				A\times B\times B	
	\end{tikzcd}
	\vline
	\,
	\vline
	\begin{tikzcd}[column sep=small]
		&
		X 
		\ar[ld, "f"']
		\ar[rd, "g"]
			&
	\\
		A
		\sar[rr, "p"']
		\doublecell[rr, shift left=2.5ex]{\xi_3}
			&
				&
				B
	\end{tikzcd}
	\]
	The first correspondence is obtained by postcomposing the cartesian cell defining $\wh{p}$.
	For the second correspondence, we use the horizontal universal property of products and the unit-pureness of $\dbl{D}$.
	Therefore, the cell $\xi_3$ exhibits $X$ as a tabulator of $p$ if and only if the cell $\xi_1$ exhibits $X$ as a co-Eilenberg-Moore object
	of $\wh{p}$.
	In particular, if $\wh{p}$ has an $\zero{M}$-co-Eilenberg-Moore object, an $\zero{M}$-tabulator.
	
	If $\wh{p}$ has a strong co-Eilenberg-Moore object,  
	we have the following diagram.
	\[
		\begin{tikzcd}[column sep=small, row sep=small]
			&&&&&
			\coEM{\wh{p}}
			\ar[llld, "u"']
			\ar[rrrd, "u"]
			\\
				&&
				A\times B
				\sar[rrrrrr, "\wh{p}"]
				\ar[drr, "\Delta_A\times\id_B"]
				\ar[dll, "\pi_1"']
					&&
					\!
					\doublecell[rr,shift right=2ex]{\cart}
					\doublecell[rr,shift left=3ex]{\opcart}
						&&
					\!
							&&
							A\times B
							\ar[dll, "\id_A\times\Delta_B"']
							\ar[drr, "\pi_2"]
								&&
			\\
			A
			\ar[drr, "\Delta_A"description]
			\sar[rrrr, ""]
			\ar[rrrr,
					shift left=2ex,
					phantom, "\opcart"{description, inner sep=0mm},
					start anchor={[xshift=0ex, yshift=0ex]center},
					end anchor={[xshift=-2ex, yshift=0ex]center},
					]
			\ar[rrrr,
					shift right=2ex,
					phantom, "\cart"{description, inner sep=0mm},
					start anchor={[xshift=0ex, yshift=0ex]center},
					end anchor={[xshift=-2ex, yshift=0ex]center},
					]
			\ar[rrrrdd,equal,bend right = 30pt]
				&&
					&&
					A\times A \times B
					\sar[rr, "\Id_A\times p\times \Id_B"]
					\ar[dll, "\pi_{12}"description]
					\doublecell[rrd]{\opcart}
						&&
						A\times B \times B
						\ar[drr, "\pi_{23}"description]
						\sar[rrrr, ""near end]
					\ar[rrrr,
							shift left=2ex,
							phantom, "\opcart"{description, inner sep=0mm},
							start anchor={[xshift=2ex, yshift=0ex]center},
							end anchor={[xshift=0ex, yshift=0ex]center},
							]
					\ar[rrrr,
							shift right=2ex,
							phantom, "\cart"{description, inner sep=0mm},
							start anchor={[xshift=2ex, yshift=0ex]center},
							end anchor={[xshift=0ex, yshift=0ex]center},
							]
							&&
							\!
								&&
								B
								\ar[dll, "\Delta_B"description]
								\ar[ddllll, equal,bend left = 30pt]
			\\
				&&
				A\times A
				\sar[rrrrrr]
				\ar[drr, "\pi_2"']
					&&\!\doublecell[rrd]{\cart}
						&\!	
						&
						\!
						&\!&
							B \times B
							\ar[dll, "\pi_1"]
								&&
			\\
			&&&&
			A
			\sar[rr, "p"']
			&&
			B
		\end{tikzcd}
	\]
	Here, the two diamond cells on both sides satisfy the Beck-Chevalley condition
	since these are obtained by the products of other cells satisfying the Beck-Chevalley condition.
	The hexagon at the bottom is obtained by applying the product $-\times-$ row by row to the following opcartesian and cartesian cells.
	The cells on the left are the vertical identities on $p$ while
	those on the middle and the right are the canonical cells for the companion and conjoint of $!\colon A\to 1$.
	\[
		\begin{tikzcd}[column sep=small, row sep=small]
			A
			\ar[r, "p"]
			\ar[d,equal]
			&
			B 
			\ar[d,equal]
			\\
			A
			\ar[r, "p"]
			\ar[d,equal]
			&
			B 
			\ar[d,equal]
			\\
			A
			\ar[r, "p"]
			&
			B 
		\end{tikzcd}
		\hspace{2ex}
		,
		\hspace{2ex}
			\begin{tikzcd}[column sep=tiny, row sep=small]
					&
					A
					\ar[ld,equal]
					\ar[rd,"!"]
						&
				\\
					A
					\ar[rd,"!"']
					\sar[rr,""]
					&
					&
					1
					\ar[ld,equal]
				\\
					&
					1
						&	
			\end{tikzcd}
			\hspace{2ex}
			,
			\hspace{2ex}
			\begin{tikzcd}[column sep=tiny, row sep=small]
				&
				B
				\ar[rd,equal]
				\ar[ld,"!"']
					&
			\\
				1
				\ar[rd,equal]
				\sar[rr,""]
				&
				&
				B
				\ar[ld,"!"]
			\\
				&
				1
					&	
		\end{tikzcd}
		\]
	Note that for these specific cells, the Beck-Chevalley condition holds without 
	any assumption on the equipment $\dbl{D}$ except for the cartesian structure.
	Back to the diagram, the bottom half is a cartesian cell by \cref{lem:Sandwich} and 
	also a horizontal cell.
	Thus, it is the identity cell on $p$, which makes the whole diagram an opcartesian cell.
	Since it is the very cell corresponding to the canonical cell for the co-Eilenberg-Moore object $\coEM{\wh{p}}$
	through the aforementioned correspondence, every tabulator is strong.
\end{proof}

\begin{remark} 
	The relationship between tabulators and co-Eilenberg-Moore objects is also observed in \cite[4.2.3]{Law15}.
	However, the argument there is more focused on left comodules and a universal property in the horizontal bicategory, that is, comonadicity.
	We will compare the comonoidicity of vertical arrows with the comonadicity of their companions in \Cref{lem:MapBisSpan}
	in a more restricted context.
\end{remark}

\begin{lemma}\label{lem:coEMpresfib}
	Let $\dbl{D}$ be a cartesian equipment with Beck-Chevalley pullbacks.
	If a vertical arrow $f\colon A\to B$ is copointalic,
	then $f$ is a fibration.
\end{lemma}
\begin{proof}
	Applying \cref{prop:BCModularLaw} to $R=!_!!^*$ and $S=\Id_B$,
	we have the following equality of cells on the left,
	where the top horizontal arrow is $f_!$ since the cartesian cell defining $f_!$ is decomposed to two cartesian cells as shown on the right.
	\begin{equation}
		\label{eqn:ModularlawForDefinable}
		\begin{tikzcd}
			A
			\sar[r,"f_!"]
			\ar[d, "\mmbox{\langle\id_A, f\rangle}"']
			\doublecell{\cart}
				&
				B
				\ar[d, "\Delta"]
			\\
			A\times B
			\sar[r, "!_!!^*\times \Id_B"]
			\ar[d, "f\times\id_B"']
			\doublecell{\opcart}
				&
				B\times B
				\ar[d, equal]
			\\
			B\times B
			\sar[r, "f^*!_!!^*\times \Id_B"']
				&
				B\times B
		\end{tikzcd}
		\hspace{1ex}
		=
		\hspace{1ex}
		\begin{tikzcd}
			A
			\sar[r,"f_!"]
			\ar[d, "f"']
			\doublecell{\opcart}
				&
				B
				\ar[d, equal]
			\\
			B
			\sar[r,"f^*f_!"']
			\ar[d, "\Delta"']
			\doublecell{\cart}
				&
				B
				\ar[d, "\Delta"]
			\\
			B\times B
			\sar[r, "f^*!_!!^*\times \Id_B"']
				&
			B\times B
		\end{tikzcd}
		,
		\hspace{3ex}
		\begin{tikzcd}[column sep=small]
			&
			A
			\ar[ld, equal]
			\ar[rd, "f"]
				&
		\\
			A
			\doublecell[rr, shift left=3ex]{\opcart}
			\sar[rr, "f_!"']
			\ar[d, "\mmbox{\langle\id_A, f\rangle}"']
				&
				&
				B
			\ar[d, "\Delta"]
		\\
			A\times B
			\doublecell[rr, shift left=3ex]{\cart} 
			\doublecell[rr, shift right=3.5ex]{\cart}
			\sar[rr, "!_!!^*\times \Id_B"']
			\ar[dr, "!\times\id_B"']
				&
				&
				B\times B
			\ar[dl, "!\times\id_B"]
			\\
				&
				1\times B
				&
		\end{tikzcd}
	\end{equation}

	Now let us observe the following sequence of bijective correspondences of cells and vertical arrows.
	The first correspondence
	is obtained by considering the cartesian cell exhibiting the horizontal composite $f^*!_!!^*$ as the restriction $f^*!_!(\id_B,!)$,
	while the third one is given by the universal property of the cartesian cell 
	at the bottom of the middle diagram above.
	The second correspondence follows from the universality of the product in the category $\dbl{D}_1$.
	Here, note that since $\dbl{D}$ is unit-pure, the vertical arrow $g$ comes with a unique arrow $\Id_X\to\Id_A$ in $\dbl{D}_1$.
	The last is precisely the universality of the (strong) tabulator, obtained by the assumption on $f$ considering \cref{rem:coEMtab}.
	\begin{equation*}
		\begin{tikzcd}[column sep=small]
				&
				X 
				\ar[ld, "g"']
				\ar[rd, "!"]
					&
			\\
			B 
			\doublecell[rr, shift left=2.5ex]{\alpha_1}
			\sar[rr, "f^*!_!"']
				&
					&
					1
		\end{tikzcd}
		\vline
		\,
		\vline
		\begin{tikzcd}[column sep=small]
				&
				X 
				\ar[ld, "g"']
				\ar[rd, "g"]
					&
			\\
			B 
			\doublecell[rr, shift left=2.5ex]{\alpha_2}
			\sar[rr, "f^*!_!!^*"']
				&
					&
					1
		\end{tikzcd}
		\vline
		\,
		\vline
		\begin{tikzcd}[column sep=small]
				&
				X
				\ar[ld, "{\langle g,g\rangle}"']
				\ar[rd, "{\langle g,g\rangle}"]
					&
			\\
				B \times B
				\doublecell[rr, shift left=2.5ex]{\alpha_3}
				\sar[rr, "f^*!_!!^*\times \Id_B"']
					&
					&
				B \times B
		\end{tikzcd}
		\vline
		\,
		\vline
		\begin{tikzcd}[column sep=small]
				&
				X
				\ar[ld, "g"']
				\ar[rd, "g"]
					&
			\\
			B 
			\doublecell[rr, shift left=2.5ex]{\alpha_4}
			\sar[rr, "f^*f_!"']
				&
					&
					B
		\end{tikzcd}
		\vline
		\,
		\vline
		\begin{tikzcd}[column sep=small]
				&
				X
				\ar[ldd, "g"', bend right]
				\ar[d, "\alpha_5"]
			\\
				&
				A
				\ar[ld, "f"']
			\\
			B 
				&
		\end{tikzcd}
	\end{equation*}
	Tracing back the correspondence, the correspondence $\alpha_5\mapsto\alpha_1$ is obtained by postcomposing the following cell, and hence this cell exhibits
	$A$ as a tabulator of $f^*!_!$ whose left leg is $f$.
	\[
		\begin{tikzcd}[column sep=small]
				&
				A
				\ar[ld, "f"']
				\ar[rd, "f"]
					&
			\\
			B
			\ar[d, "\Delta"']
			\doublecell[rrd]{\cart}
			\sar[rr, "f^*f_!"{description, near end}]
			\doublecell[rr, shift left=2.5ex]{\opcart}
				&
					&
					B
					\ar[d, "\Delta"]
			\\
			B\times B
			\ar[d, "\pi_1"]
			\doublecell[rrd]{\pi_1}
			\sar[rr, "f^*!_!!^*\times\Id_B"]
				&
					&
					B\times B
					\ar[d, "\pi_1"]
			\\
			B
			\ar[d, equal]
			\doublecell[rrd]{\cart}
			\sar[rr, "f^*!_!!^*"]
				&
					&
					B
					\ar[d, "!"]
			\\
			B
			\sar[rr, "f^*!_!"']
				&
					&
					1
		\end{tikzcd}
		\vspace{-4ex}
	\]
\end{proof}

\begin{theorem}
	The following are equivalent for a unit-pure cartesian equipment $\dbl{D}$
	and a stable system $\zero{M}$ on $\oneV{\dbl{D}}$.
	\begin{itemize}
		\item[i)] $\dbl{D}$ admits an $\zero{M}$-comprehension scheme,
		\item[ii)] $\dbl{D}$ admits a unary $\zero{M}$-comprehension scheme w.r.t. copointals.
	\end{itemize}
	Moreover, if $\dbl{D}$ is discrete, they are also equivalent to;
	\begin{enumerate}
		\item[iii)] $\dbl{D}$ admits a unary $\zero{M}$-comprehension scheme w.r.t. comonoids.
	\end{enumerate}
\end{theorem}
\begin{proof}
	$i)\Rightarrow ii)$ follows from \cref{rem:coEMtab} and the fact that
	if a vertical arrow $f\colon A\to B$ is in $\zero{M}$, then so is $f\fatsemi\Delta_B$.
	This is because $\Delta_B$ belongs to $\zero{M}=\Fib(\dbl{D})$ in a unit-pure 
	\acl{DCR}. 
	For $ii)\Rightarrow i)$, suppose the second condition holds.
	$\dbl{D}$ has strong $\zero{M}$-tabulators by \cref{prop:coEMtoTabAl},
	especially, $\dbl{D}$ has Beck-Chevalley pullbacks.
	By \cref{lem:coEMpresfib}, every vertical arrow in $\zero{M}$ is a fibration.
	Thus, $\dbl{D}$ admits a left-sided $\zero{M}$-comprehension scheme
	by \cref{rem:fibandcompr}.
	The equivalence of $ii)$ and $iii)$ in the discrete case
	follows from \cref{prop:odotISland}, \cref{cor:copointalcomonoid}, and \cref{cor:copointalcomodule}.
\end{proof}

\begin{corollary}[cf. {\cite[Theorem 5.2]{LWW10} and \cite[Theorem 5.3.2]{Ale18}}]\label{cor:Spanchacompa}
	\label{cor:CharacteriseSpan}
	The following are equivalent for a double category $\dbl{D}$,
	where $\Mor$ denotes the class of all vertical arrows.
	\begin{enumerate}
		\item $\dbl{D}$ is equivalent to $\Span(\one{C})$ for a category $\one{C}$ with finite limits.
		\item $\dbl{D}$ is a unit-pure cartesian equipment and admits a unary $\Mor$-comprehension scheme w.r.t. copointals;
		i.e., it is a unit-pure cartesian equipment with strong co-Eilenberg-Moore objects for copointals and
		every vertical arrow is copointalic.
		\item $\dbl{D}$ is a unit-pure discrete cartesian equipment and admits a unary $\Mor$-comprehension scheme w.r.t. comonoids;
		i.e., it is a unit-pure cartesian equipment with strong co-Eilenberg-Moore objects for comonoids and every vertical arrow is comonoidic.
	\end{enumerate}
\end{corollary}
\begin{proof} 
	Applying \Cref{thm:unitpuredbl} to the case of $\zero{M}=\Mor$,
	a double category $\dbl{D}$ is equivalent to $\Span(\one{C})$ for a category $\one{C}$ with finite limits 
	if and only if it is a unit-pure cartesian equipment with a $\Mor$-comprehension scheme.
	With this and the previous theorem, the first two conditions are equivalent.
	Since $\Span(\one{C})$ is always a discrete double category, 
	the last condition is equivalent to the first two.
\end{proof}

\begin{remark}\label{rem:SpanchaAleiferi}
	Aleiferi states in her PhD thesis \cite{Ale18} that the double categories of spans are 
	characterised by the same conditions as $ii)$ in the above corollary except for the last condition:
	every vertical arrow is copointalic.
	However, without it, the double category of spans cannot be characterised.
	Indeed, every unit-pure double category of relations, say $\Rel{}(\Set)$, 
	also has strong co-Eilenberg-Moore objects for copointals.
\end{remark}

We explain how the characterisation \cite[Theorem 5.2]{LWW10} of spans as a cartesian bicategory is partially reconstructed through our characterisation.
Firstly, let us admit the following fact obtained by assembling some results from \cite{LWW10}.
Note that, in a Cauchy unit-pure cartesian equipment $\dbl{D}$, the notion of a discrete object in $\dbl{D}$ as defined in \cref{defn:discrete}
and that in a cartesian bicategory $\biH{\dbl{D}}$ in the sense of \cite[Definition 3.10]{LWW10} coincide.
\begin{fact}[{\cite[Proposition 2.3, Corollary 3.5, Theorem 3.14, Lemma 4.5]{LWW10}}]
	\label{fact:LWW10}
	Suppose we are given a cartesian bicategory $\bi{B}$, and every map in $\bi{B}$ is comonadic.
	Then, maps are right-cancellable by maps \cite[Lemma 4.5]{LWW10},
	and
	there exists the map double category $\Map(\bi{B})$ (see \cref{def:mapdoublecategory} and the proceeding discussion)
	that is unit-pure, cartesian, and discrete.
\end{fact}
Through this fact, \cref{cor:CharacteriseSpan} results in the main theorem of \cite{LWW10}.
\begin{theorem}[{\cite[Theorem 5.2]{LWW10}}]
	A cartesian bicategory $\bi{B}$ is biequivalent to the bicategory of spans $\biH{\Span(\one{C})}$ for some finitely complete category $\one{C}$,
	if and only if
	every comonad has a co-Eilenberg-Moore object and every map is comonadic.
\end{theorem}
\begin{proof}
	The \textit{only if} part is easily checked; observe that a comonad in $\biH{\Span(\one{C})}$ is the same as a span whose two legs are the same.
	This is verified directly, but one can also check it by utilising \cref{cor:CharacteriseSpan} and \cref{cor:copointalcomonoid} and observing
	a copointal is exactly a span whose legs are the same. See the introduction of \cite{LWW10} for the remainder.

	For the \textit{if} part, suppose that every comonad has a co-Eilenberg-Moore object and every map is comonadic in $\bi{B}$.
	It suffices to show that the map double category $\Map(\bi{B})$ obtained in \cref{fact:LWW10} is equivalent to $\Span(\one{C})$ for some finitely complete category
	$\one{C}$. The following lemma shows that $\Map(\bi{B})$ satisfies $iii)$ of \cref{cor:CharacteriseSpan}, which completes the proof.
\end{proof}
\begin{lemma}
	\label{lem:MapBisSpan}
	Let $\dbl{D}$ be a Cauchy and unit-pure cartesian equipment.
	Then every comonoid in $\dbl{D}$ has a strong co-Eilenberg-Moore object if
	the comonoid, as a comonad in $\biH{\dbl{D}}$, has a co-Eilenberg-Moore object in the bicategory $\biH{\dbl{D}}$.
	Moreover, every vertical arrow $f$ is comonoidic if the map $f_!$ is comonadic in $\biH{\dbl{D}}$.
\end{lemma}
\begin{proof}
	Suppose $G$ is a comonoid on an object $A$ in $\dbl{D}$, or equivalently, a comonad in $\biH{\dbl{D}}$.
	Suppose, moreover, that $G$ has a co-Eilenberg-Moore object $(X,u,\mu)$ as the following diagram.
	Since $u$ is in particular a left adjoint in $\biH{\dbl{D}}$, we have the companion $f$ of $u$.
	Utilising the unit-pureness of $\dbl{D}$,
	it is observed that the triangle cell $\bar\mu$ obtained by the following composite is a comodule for $G$ in the sense of \cref{defn:comodule}.
	As in \cref{prop:repradj}, for each companion, we write $\alpha$ and $\beta$ for the canonical cartesian and opcartesian cells.
	\vspace{-2ex}
	\[
		\begin{tikzcd}[column sep=small]
				&
				X
				\ar[ld, "f"']
				\ar[rd, "f"]
					&
			\\
			A
			\sar[rr, "G"']
			\doublecell[rr, shift left=3ex]{\bar\mu}
				&
					&
					A
		\end{tikzcd}
		=
		\begin{tikzcd}
				&
				X
				\ar[ld, equal]
				\ar[rd, "f"]
					&
			\\
			X
			\sar[rr, "u"]
			\doublecell[rrd]{\mu}
			\doublecell[rr, shift left=4ex]{\beta}
			\ar[d, equal]
				&
					&
					A
					\ar[d, equal]
			\\
			X
			\sar[r, "u"]
			\ar[d, "f"']
			\doublecell[r, shift right=2.2ex, near start]{\alpha}
				&
				A
				\sar[r, "G"]
				\ar[ld, equal]
					&
					A
					\ar[d, equal]
					\doublecell[lld]{\verteq}
			\\
			A
			\sar[rr, "G"']
				&
					&
					A
		\end{tikzcd}
	\]
	Note that $\bar\mu$ is opcartesian because there is an invertible horizontal cell $G\cong f^*f_!$ since $f_!=u$ is comonadic.
	On the other hand, any comodule $(Y,g,\nu)$ of $G$ gives rise to a left comodule $(Y,g_!,\nu_!)$ of the comonad $G$ in $\biH{\dbl{D}}$ as follows.
	\vspace{-1ex}
	\[
		\begin{tikzcd}
			Y
			\sar[rr, "g_!"]
			\doublecell[rrd]{\nu_!}
			\ar[d, equal]
				&
					&
					A
					\ar[d, equal]
			\\
			Y
			\sar[r, "g_!"]
				&
				A
				\sar[r, "G"]
					&
					A
		\end{tikzcd}
		\hspace{2ex}
		:=
		\hspace{2ex}
		\begin{tikzcd}
				&
				Y
				\sar[r, "g_!"]
				\ar[d, "g"description]
				\ar[rd, "g"description]
				\ar[ld, equal]
				\doublecell[r, shift right=2ex, near end]{\alpha}
					&
					A
					\ar[d, equal]
			\\
			Y
			\sar[r, "g_!"']
			\ar[phantom, r, shift left=2ex, "\beta"{description, xshift=1ex, inner sep=0mm}]{}
				&
				A
				\sar[r, "G"']
				\doublecell[r, shift left=2ex, near start]{\nu}
					&
					A
		\end{tikzcd}
	\]
	Therefore, given a comodule $(Y,g,\nu)$ of $G$, there exists a horizontal arrow $v\colon Y \sto X$ equipped with
	the cartesian and opcartesian cells satisfying the first equation of the following.
	Such a horizontal arrow $v$ is unique up to invertible horizontal cells. 
	Since maps are right-cancellable (\cite[Lemma 4.4]{LWW10}), $v$ is also a map.
	Moreover, since $\dbl{D}$ is Cauchy and unit-pure,
	there exists a unique $h\colon Y \to X$ representing $v$ as its companion.
	Again by unit-pureness, we have $h\fatsemi f=g$ and the cartesian and opcartesian cells in the middle of the following equation
	are obtained by composing the canonical cells $\alpha$ and $\beta$ defining the companions of $h$ and $f$. This shows that the last equation follows.
	\vspace{-2ex}
	\[
		\begin{tikzcd}[column sep=small]
				&
				Y
				\ar[ld, "g"']
				\ar[rd, "g"]
					&
			\\
			A
			\sar[rr, "G"']
			\doublecell[rr, shift left=3ex]{\nu}
				&
					&
					A
		\end{tikzcd}
		=
		\begin{tikzcd}
				&
					&
					Y
					\ar[lld, equal]
					\ar[rd, "g"]
						&
			\\
			Y
			\doublecell[rrr, shift left=3ex]{\opcart}
			\sar[r, "v"']
			\ar[d, equal]
			\doublecell[rd]{\verteq}
				&
				X
				\sar[rr, "u"]
				\doublecell[rrd]{\mu}
				\ar[d, equal]
					&
						&
						A
						\ar[d, equal]
			\\
			Y
			\sar[r, "v"]
			\ar[rd, "g"']
			\doublecell[rr, shift right=3ex]{\cart}
				&
				X
				\sar[r, "u"]
					&
					A
					\sar[r, "G"]
					\ar[ld, equal]
					\doublecell[d]{\verteq}
						&
						A
						\ar[ld, equal]
			\\
				&	
				A
				\sar[r, "G"']
					&
					A
						&
		\end{tikzcd}
		=
		\begin{tikzcd}[column sep=small]
				&
				Y
				\ar[d, "h"]
				\ar[ldd, "g"', bend right]
				\ar[rdd, "g", bend left]
					&
			\\
				&
				X
				\ar[ld, "f"']
				\ar[rd, "f"]
					&
			\\
			A
			\sar[rr, "G"']
			\doublecell[rr, shift left=3ex]{\bar\mu}
				&
					&
					A
		\end{tikzcd}
	\]
	Moreover, such a vertical arrow is unique since $\dbl{D}$ is unit-pure and $v$ is unique up to horizontal invertible cells.
	This shows that $(X,f,\bar\mu)$ is a co-Eilenberg-Moore object of the comonoid $G$.

	Let $f$ be a vertical arrow whose companion $f_!$ is comonadic in $\biH{\dbl{D}}$. The same discussion shows that
	there is a strong co-Eilenberg-Moore object of $f^*f_!$ whose leg is isomorphic to $f$ in $\biV{\dbl{D}}$, and the leg is precisely $f$ because $\dbl{D}$ is unit-pure.
	This shows $f$ is comonoidic.
\end{proof}

\section{Future Work}\label{sec:futurework}
\addtocontents{toc}{\protect\setcounter{tocdepth}{1}}

\subsection*{Applications to categorical logic}
One of the leading motivations for studying double categories
of relations is the desire to ground some logic upon them.
Classically, categorical logic is based on the theory of \textit{hyperdoctrines}, or more generally, \textit{fibrations}.
On the other hand, recall that our double categories of relations are constructed through the $\dbl{F}\mathrm{r}$-construction \cref{lem:FrConstruction} from a monoidal bifibration.
In general, we are expecting that from any monoidal fibration $P\colon\one{E}\to\one{B}$ with $\one{B}$ cartesian monoidal,
we can construct a \textit{virtual double category} $\Bil(P)$, which generalises $\dbl{F}\mathrm{r}$-construction and thus our double category of relations,
in the following way.
\begin{itemize}
	\item %
		The vertical category $\oneV{\Bil(P)}$ is the base category $\one{B}$.
	\item %
		A horizontal arrow $e\colon A\sto B$ is an object $e\in\one{E}_{A\times B}$, which can be seen as a \textit{bilateral term}.
	\item %
		A cell
		\[
			\begin{tikzcd}
				A_0
				\sar[r,"p_1"]
				\ar[d,"f"']
				\doublecell[rrrd]{\alpha}
					&
					A_1
					\sar[r,"p_2"]
						&
						\cdots
						\sar[r]
							&
							A_n
							\ar[d,"g"]
				\\
				B_0
				\sar[rrr,"q"']
					&
						&
							&
							B_1
			\end{tikzcd}
		\]
		is defined as an arrow $\alpha\colon\land\vec{p}\to q$ in $\one{E}$
		over $\lambda\vec{x}. \langle f(x_0),g(x_n)\rangle\colon\prod_{i=0,\ldots,n}A_i\to B_0\times B_1$ in $\one{B}$,
		where $\land\vec{p}:=\bigotimes_{i=1,\ldots,n}\pi_{i}^*p_i\in\one{E}_{\prod_j A_j}$ is defined as the $n$-times monoidal product in $\one{E}_{\prod_j A_j}$
		and $\pi_{i}\colon\prod_jA_j\to A_{i-1}\times A_i$
		is the projection for each $i=1,\ldots,n$.
\end{itemize}

Since a suitable class of stable orthogonal factorisation systems, as bifibrations of the form $\one{M}\to\one{C}$ in \cref{prop:OFSfromStableSystem},
provides a categorical treatment for regular theories,
the double category of relations can naturally be seen as a framework to take a model of regular theories.
In the opposite direction, we can construct an internal language for virtual double categories, as proposed in \cite{Nas24}.
Considering such a language of virtual double categories could bring us a new treatment of categorical logic that is more general than the one based on (ordinary) categories.
We expect that there will be some advantages in translating several notions of logical completion of hyperdoctrines, including the ex/reg completion and tripos-to-topos construction, into the double-categorical setting.
Indeed, the Cauchisation in \cref{subsec:cauchyreg} is a logical completion corresponding to constructing the free Cauchy-complete doctrine in \cite{Pas16}.
Other classes of hyperdoctrines and their logical completions would be translated into double categories similarly.

\subsection*{Functorial correspondence between \acsp{DCR} and \acsp{SOFS}}
Although we have explained how double categories and stable orthogonal factorisation systems are related,
the correspondence is not yet functorial. To obtain the result as equivalence of some 2-categories,
we have to determine the 2-category of factorisation systems. From the perspective of the general construction $\Bil$ above,
the construction from a factorisation system should start from some sub 2-category of the 2-category of monoidal fibrations.
In such 2-category, morphisms between \acsp{OFS} need not preserve the left class, and hence it must be different from some (2-)categories of factorisation systems
considered in literature like \cite{Ste24}.
On the other hand, since any pseudo double functor preserves Beck-Chevalley squares, it preserves final morphisms; i.e., the left class.
We expect that lax double functors corresponds to morphisms of monoidal fibrations, while pseudo double functors corresponds to morphisms of monoidal bifibrations,
which also preserves left class.

\subsection*{\Acl{VDCAR}}
Given an \textit{\ac{AWFS}} \cite{BG16} $(L, R)$ on a finitely complete category $\one{C}$,
we can also consider \textit{algebraic relations} as follows.
Recall from \cite[Theorem 9]{BG16} that
a pair $(R\mhyphen\Algdbl, V)$ that satisfies the following conditions serves as the characterisation of an \ac{AWFS} on $\one{C}$.
\begin{itemize}
	\item
		$R\mhyphen\Algdbl$ is a strict double category, and $V\colon R\mhyphen\Algdbl\to\Sq(\one{C})$ is a strict double functor.
		By $\Sq(\one{C})$, 
		we mean the strict double category induced from the following cocategory object in $\Catbi$ through taking the powers
		on $\one{C}$ when we write $\left[n\right]$ for the chain of length $n$.
		\[
			\begin{tikzcd}
				\left[2\right]
				\al[r, "\boundary_1"description]
				&
				\left[1\right]
				\al[r, "\boundary_1"description,  shift left=2ex]
				\ar[r, "\sigma"description]
				\al[r, "\boundary_0"description, shift right=2ex]
				&
				\left[0\right]
			\end{tikzcd}
		\]
	\item
		$R\mhyphen\Algdbl_0=\one{C}$, $V_0=\id_\one{C}$, and $V_1\colon R\mhyphen\Algdbl_1\to\one{C}^{\left[1\right]}$ is faithful.
	\item
		$V_1$ is strictly monadic. We write $R\mhyphen\Algone$ for $R\mhyphen\Algdbl_1$.
	\item
		$V_1$ is a discrete pullback-fibration;
		for every $\mathbf{g}\in R\mhyphen\Algone$, $f\in\one{C}^{\left[1\right]}$, and pullback square $(k,h)\colon f\to V(\mathbf{g})$,
		there exists a unique pair $(\mathbf{f}, \varphi\colon \mathbf{f}\to\mathbf{g})$ satisfying
		$V_1(\varphi)=(k, h)$, and this arrow $\varphi$ is cartesian
		for $\tgt\colon R\mhyphen\Algone\to\one{C}$.
	\item[-]
		$\tgt\colon R\mhyphen\Algone\to\one{C}$ is a fibration and $V_1\colon \tgt\to\cod^\one{C}$ is a morphism of fibrations over $\one{C}$.
	\item[-]
		For each $C\in\one{C}$, $V_1$ restricts to a strictly monadic functor $R/C\mhyphen\one{Alg}\to \one{C}/C$ on fibres for
		$\tgt$ and $\cod^\one{C}$.
\end{itemize}
(The last two are redundant for the characterisation.)
The point here is that an \ac{AWFS} is mostly determined by data of fibration, which is a generalisation of \cref{prop:OFSfromStableSystem}.

Given an AWFS $(R\mhyphen\Algdbl, V)$, we define an \textit{algebraic relation} $\mathbf{p}\colon A\sto B$ between two objects $A, B\in\one{C}$ as
an algebra $\mathbf{p}$ whose target $\tgt(\mathbf{p})$ is $A\times B$.
In other words, it is just a \textit{bilateral term} of the (monoidal) fibration $\tgt\colon R\mhyphen\Algone\to\one{C}$.
Applying the $\Bil$ construction above to this monoidal fibration, we are expecting that an \ac{AWFS}
induces a virtual double category, and this leads us to consider
the notion of \ac{VDCAR} in \cref{theTable}.

A \ac{VDCAR} that corresponds to an orthogonal factorisation system should be called
a \ac{VDCR}. \cite{HJ03} and \cite{BK17} characterise fibrations obtained from the right classes of orthogonal factorisation systems,
and in \cite{BK17}, they said the resulting factorisation system \textit{comprehensive factorisation system}.
Utilising these results and generalising our \cref{thm:CompSchemeToFactorisation} to the virtual settings,
we would complete a double categorical framework for comprehensive factorisation systems that encompasses all the examples in \cite{BK17}.

From the same point of view, the term
`double category of relations' should be used for \acp{VDCR} such that they are double categories.
Then the question arises here whether our definition of \ac{DCR} is consistent with the double category of relations in this sense (written as \ac{DCR}' in the table).
Nevertheless, we can find a clue for this issue that supports the negative expectation.
According to \cite[Theorem  14.4]{Shu08}, the construction $\dbl{F}\mathrm{r}$ is applicable not only for strong \acs{BC} monoidal fibrations, but also for \textit{internally closed} \textit{weakly \acs{BC}} monoidal fibrations.
This implies that if there exists any orthogonal factorisation system whose accompanying fibration $\one{M}$, as in \cref{prop:OFSfromStableSystem}, is not
strong \acs{BC} but admits the construction $\dbl{F}\mathrm{r}$, then \ac{DCR}' and our \ac{DCR} do not coincide.

On the other hand, one can define a \textit{stable} \ac{AWFS} as follows. An \ac{AWFS} is called stable if,
for its category of left coalgebras $L\mhyphen\one{Coalg}$, the forgetful functor $L\mhyphen\one{Coalg}\to\one{C}^{\left[2\right]}$ is
pullback-discrete. We expect that a stable \ac{AWFS} induces a \ac{DCAR}.

\begin{table}[h]
\centering
\begin{tabular}{|c|c|}
	\hline
		\acp{AWFS}
		on finitely complete categories
	& \Acp{VDCAR}
	\\
	\hline
	\hline
	\begin{tikzpicture}[scale=0.7]
		\node[myrectanglefs] (Z)
			at (1, 7)
			{\ac{AWFS}};
		\node[myrectanglefs, text=gray] (Y)
			at (0, 5)
			{?};
		\node[myrectanglefs] (X1)
			at (-1, 3)
			{Stable \ac{AWFS}};
		\node[myrectanglefs] (X2)
			at (4, 6)
			{\ac{OFS}};
		\node[myrectanglefs, text=gray] (A')
			at (3, 4)
			{?};
		\node[myrectanglefs] (A)
			at (2, 2)
			{\ac{SOFS}};
		\node[myrectanglefs, text width=15mm] (C1)
			at (1, -0)
			{left-proper \ac{SOFS}};
		\node[myrectanglefs, text width=15mm] (C2)
			at (6, 0)
			{right-proper \ac{SOFS}};
		\node[myrectanglefs, text width=20mm] (C'1)
			at (0, -2)
			{anti-right-proper \ac{SOFS}};
		\node[myrectanglefs, text width=15mm] (D')
			at (4, -4)
			{regular \ac{OFS}};
		\node[myrectanglefs, text width=15mm] (D)
			at (5, -2)
			{proper \ac{SOFS}};
		\node[myrectanglefs, text width=15mm] (F)
			at (-1, -4)
			{$(\Iso,\Mor)$};
		\draw[-] (Z) -- (Y);
		\draw[-] (Z) to [out=210, in=100] (X1);
		\draw[-] (Z) -- (X2);
		\draw[dashed] (Y) -- (X1);
		\draw[-] (A') -- (A);
		\draw[-] (X2) -- (A');
		\draw[-] (X1) -- (A);
		\draw[-] (Y) -- (A');
		\draw[-] (C1) -- (C'1);
		\draw[-] (C1) -- (D);
		\draw[-] (D) -- (D');
		\draw[-] (A) -- (C1);
		\draw[-] (A) -- (C2);
		\draw[-] (C2) -- (D);
		\draw[-] (C'1) -- (D');
		\draw[-] (C'1) -- (F);
	\end{tikzpicture}
	&
	\begin{tikzpicture}[scale=0.7]
		\node[myrectangledbl, text=gray] (Z)
			at (1, 7)
			{\ac{VDCAR}};
		\node[myrectangledbl, text=gray] (Y)
			at (0, 5)
			{\ac{DCAR}};
		\node[myrectangledbl, text=gray] (X1)
			at (-1, 3)
			{?};
		\node[myrectangledbl, text=gray] (X2)
			at (4, 6)
			{\ac{VDCR}};
		\node[myrectangledbl, text=gray] (A')
			at (3, 4)
			{\ac{DCR}'};
		\node[myrectangledbl] (A)
			at (2, 2)
			{\ac{DCR} in this paper \\
			\ref{thm:MainThm}};
		\node[myrectangledbl] (C1)
			at (1, 0)
			{unit-pure \ac{DCR}\\
			\ref{defn:unitpure}};
		\node[myrectangledbl] (C2)
			at (6, 0)
			{locally preordered \ac{DCR}\\ \ref{lem:DagLocPos}};
		\node[myrectangledbl] (C'1)
			at (0, -2)
			{unit-pure Cauchy \ac{DCR} \\ \ref{thm:Cauchydouble}};
		\node[myrectangledbl] (D')
			at (4, -4)
			{\ac{DCR}\,on\,regular\,categories \\
			\ref{thm:Regcat}, \ref{thm:RegcatCW} \\\cite{Lam22, CW87}};
		\node[myrectangledbl] (D)
			at (5, -2)
			{locally posetal \ac{DCR}\\ \ref{cor:locpos}};
		\node[myrectangledbl] (F)
			at (-1, -4)
			{\acs{DC} of spans \\ \ref{cor:Spanchacompa} \\\cite{LWW10,Ale18}};
		\draw[-] (Z) -- (Y);
		\draw[-] (Z) to [out=210, in=100] (X1);
		\draw[-] (Z) -- (X2);
		\draw[dashed] (Y) -- (X1);
		\draw[-] (A') -- (A);
		\draw[-] (X2) -- (A');
		\draw[-] (X1) -- (A);
		\draw[-] (Y) -- (A');
		\draw[-] (C1) -- (C'1);
		\draw[-] (C1) -- (D);
		\draw[-] (D) -- (D');
		\draw[-] (A) -- (C1);
		\draw[-] (A) -- (C2);
		\draw[-] (C2) -- (D);
		\draw[-] (C'1) -- (F);
		\draw[-] (C'1) -- (D');
	\end{tikzpicture}
	\\
	\hline
\end{tabular}
\caption{Conjectural correspondence between \acfp{AWFS} and \acfp{VDCAR}}
\label{theTable}
\end{table}

\addtocontents{toc}{\protect\setcounter{tocdepth}{2}}

\bibliographystyle{halpha-abbrv}
\bibliography{bibliography}

\begin{thebibliography}{CKWW07}
\expandafter\ifx\csname url\endcsname\relax
  \def\url#1{\texttt{#1}}\fi
\expandafter\ifx\csname doi\endcsname\relax
  \def\doi#1{\burlalt{doi:#1}{http://dx.doi.org/#1}}\fi
\expandafter\ifx\csname urlprefix\endcsname\relax\def\urlprefix{URL }\fi
\expandafter\ifx\csname href\endcsname\relax
  \def\href#1#2{#2}\fi
\expandafter\ifx\csname burlalt\endcsname\relax
  \def\burlalt#1#2{\href{#2}{#1}}\fi

\bibitem[AEBSV01]{AESV01}
J.~Ad\'amek, R.~El~Bashir, M.~Sobral, and J.~Velebil.
\newblock On functors which are lax epimorphisms.
\newblock {\em Theory Appl. Categ.}, 8:509--521, 2001.

\bibitem[Ale18]{Ale18}
E.~Aleiferi.
\newblock {\em Cartesian double categories with an emphasis on characterizing spans}.
\newblock PhD thesis, Dalhousie University, Halifax, Nova Scotia, 2018, \burlalt{1809.06940}{http://arxiv.org/abs/1809.06940}.

\bibitem[BD86]{BD86}
F.~Borceux and D.~Dejean.
\newblock Cauchy completion in category theory.
\newblock {\em Cahiers Topologie G\'{e}om. Diff\'{e}rentielle Cat\'{e}g.}, 27(2):133--146, 1986.

\bibitem[B{\'{e}}n67]{Ben67}
J.~B{\'{e}}nabou.
\newblock Introduction to bicategories.
\newblock In {\em Reports of the {M}idwest {C}ategory {S}eminar}, volume No. 47 of {\em Lecture Notes in Math.}, pages 1--77. Springer, Berlin-New York, 1967.

\bibitem[Bet96]{Bet96}
R.~Betti.
\newblock Formal theory of internal categories.
\newblock {\em Matematiche (Catania)}, 51:35--52, 1996.

\bibitem[BG16]{BG16}
J.~Bourke and R.~Garner.
\newblock Algebraic weak factorisation systems {I}: {A}ccessible {AWFS}.
\newblock {\em J. Pure Appl. Algebra}, 220(1):108--147, 2016.
\newblock \doi{10.1016/j.jpaa.2015.06.002}.

\bibitem[BK17]{BK17}
C.~Berger and R.~M. Kaufmann.
\newblock Comprehensive factorisation systems.
\newblock {\em Tbilisi Math. J.}, 10(3):255--277, 2017.
\newblock \doi{10.1515/tmj-2017-0112}.

\bibitem[Bou87]{Bou87}
D.~Bourn.
\newblock The shift functor and the comprehensive factorization for internal groupoids.
\newblock {\em Cahiers Topologie G\'eom. Diff\'erentielle Cat\'eg.}, 28(3):197--226, 1987.

\bibitem[CJKP97]{CJKP97}
A.~Carboni, G.~Janelidze, G.~M. Kelly, and R.~Par\'{e}.
\newblock On localization and stabilization for factorization systems.
\newblock {\em Appl. Categ. Structures}, 5(1):1--58, 1997.
\newblock \doi{10.1023/A:1008620404444}.

\bibitem[CKS84]{CKS84}
A.~Carboni, S.~Kasangian, and R.~Street.
\newblock Bicategories of spans and relations.
\newblock {\em J. Pure Appl. Algebra}, 33(3):259--267, 1984.
\newblock \doi{10.1016/0022-4049(84)90061-6}.

\bibitem[CKWW07]{CKWW07}
A.~Carboni, G.~M. Kelly, R.~F.~C. Walters, and R.~J. Wood.
\newblock Cartesian bicategories {II}.
\newblock {\em Theory Appl. Categ.}, 19:93--124, 2007.

\bibitem[CS10]{CS10}
G.~S.~H. Cruttwell and M.~A. Shulman.
\newblock A unified framework for generalized multicategories.
\newblock {\em Theory Appl. Categ.}, 24:No. 21, 580--655, 2010.

\bibitem[CW87]{CW87}
A.~Carboni and R.~F.~C. Walters.
\newblock Cartesian bicategories. {I}.
\newblock {\em J. Pure Appl. Algebra}, 49(1-2):11--32, 1987.
\newblock \doi{10.1016/0022-4049(87)90121-6}.

\bibitem[Day77]{Day77}
B.~J. Day.
\newblock Density presentations of functors.
\newblock {\em Bull. Austral. Math. Soc.}, 16(3):427--448, 1977.
\newblock \doi{10.1017/S0004972700023509}.

\bibitem[EBV02]{EV02}
R.~El~Bashir and J.~r. Velebil.
\newblock Simultaneously reflective and coreflective subcategories of presheaves.
\newblock {\em Theory Appl. Categ.}, 10:No. 16, 410--423, 2002.

\bibitem[Ehr63]{Ehr63}
C.~Ehresmann.
\newblock Cat\'{e}gories doubles et cat\'{e}gories structur\'{e}es.
\newblock {\em C. R. Acad. Sci. Paris}, 256:1198--1201, 1963.

\bibitem[FK72]{FK72}
P.~J. Freyd and G.~M. Kelly.
\newblock Categories of continuous functors. {I}.
\newblock {\em J. Pure Appl. Algebra}, 2:169--191, 1972.
\newblock \doi{10.1016/0022-4049(72)90001-1}.

\bibitem[FS90]{FS90}
P.~J. Freyd and A.~Scedrov.
\newblock {\em Categories, allegories}, volume~39 of {\em North-Holland Mathematical Library}.
\newblock North-Holland Publishing Co., Amsterdam, 1990.

\bibitem[GP99]{GP99}
M.~Grandis and R.~Pare.
\newblock Limits in double categories.
\newblock {\em Cahiers Topologie G\'{e}om. Diff\'{e}rentielle Cat\'{e}g.}, 40(3):162--220, 1999.

\bibitem[GP04]{GP04}
M.~Grandis and R.~Pare.
\newblock Adjoint for double categories. {A}ddenda to: ``{L}imits in double categories''.
\newblock {\em Cahiers Topologie G\'{e}om. Diff\'{e}rentielle Cat\'{e}g.}, 45(3):193--240, 2004.

\bibitem[Gra20]{Gra20}
M.~Grandis.
\newblock {\em Higher Dimensional Categories}.
\newblock World Scientific Publishing Co. Pte. Ltd., Hackensack, NJ, 2020.
\newblock From double to multiple categories.

\bibitem[HJ03]{HJ03}
J.~Hughes and B.~Jacobs.
\newblock Factorization systems and fibrations: Toward a fibred {B}irkhoff variety theorem.
\newblock {\em Electronic Notes in Theoretical Computer Science}, 69:156--182, 2003.
\newblock \doi{https://doi.org/10.1016/S1571-0661(04)80564-4}.
\newblock CTCS'02, Category Theory and Computer Science.

\bibitem[HNT20]{HNT20}
S.~N. Hosseini, A.~R. S.~A. Nasab, and W.~Tholen.
\newblock Fraction, restriction, and range categories from stable systems of morphisms.
\newblock {\em J. Pure Appl. Algebra}, 224(9):106361, 28, 2020.
\newblock \doi{10.1016/j.jpaa.2020.106361}.

\bibitem[HNTY22]{HNST22}
S.~N. Hosseini, A.~R. S.~A. Nasab, W.~Tholen, and L.~Yeganeh.
\newblock Quotients of span categories that are allegories and the representation of regular categories.
\newblock {\em Appl. Categ. Structures}, 30(6):1177--1201, 2022.
\newblock \doi{10.1007/s10485-022-09687-9}.

\bibitem[IK86]{IK86}
G.~B. Im and G.~M. Kelly.
\newblock On classes of morphisms closed under limits.
\newblock {\em J. Korean Math. Soc.}, 23(1):1--18, 1986.

\bibitem[Joh02]{Joh02}
P.~T. Johnstone.
\newblock {\em Sketches of an Elephant: a Topos Theory Compendium. {V}ol. 1}, volume~43 of {\em Oxford Logic Guides}.
\newblock The Clarendon Press, Oxford University Press, New York, 2002.

\bibitem[JW00]{JW00}
R.~Jayewardene and O.~Wyler.
\newblock Categories of relations and functional relations.
\newblock volume~8, pages 279--305. 2000.
\newblock \doi{10.1023/A:1008651524610}.
\newblock Papers in honour of Bernhard Banaschewski (Cape Town, 1996).

\bibitem[Kaw73]{Kaw73}
Y.~Kawahara.
\newblock Relations in categories with pullbacks.
\newblock {\em Mem. Fac. Sci. Kyushu Univ. Ser. A}, 27:149--173, 1973.
\newblock \doi{10.2206/kyushumfs.27.149}.

\bibitem[Kel82]{Kel82}
G.~M. Kelly.
\newblock {\em Basic concepts of enriched category theory}, volume~64 of {\em London Mathematical Society Lecture Note Series}.
\newblock Cambridge University Press, Cambridge-New York, 1982.

\bibitem[Kel91]{Kel91}
G.~M. Kelly.
\newblock A note on relations relative to a factorization system.
\newblock In {\em Category theory ({C}omo, 1990)}, volume 1488 of {\em Lecture Notes in Math.}, pages 249--261. Springer, Berlin, 1991.
\newblock \doi{10.1007/BFb0084224}.

\bibitem[Kle70]{Kl70}
A.~Klein.
\newblock Relations in categories.
\newblock {\em Illinois J. Math.}, 14:536--550, 1970.
\newblock \urlprefix\url{http://projecteuclid.org/euclid.ijm/1256052950}.

\bibitem[Kou14]{Kou14}
S.~R. Koudenburg.
\newblock On pointwise {K}an extensions in double categories.
\newblock {\em Theory Appl. Categ.}, 29:No. 27, 781--818, 2014.

\bibitem[Kou22]{Kou22}
S.~R. Koudenburg.
\newblock Augmented virtual double categories, 2022, \burlalt{1910.11189}{http://arxiv.org/abs/1910.11189}.

\bibitem[Lam22]{Lam22}
M.~Lambert.
\newblock Double categories of relations.
\newblock {\em Theory Appl. Categ.}, 38:Paper No. 33, 1249--1283, 2022.
\newblock \doi{10.1002/num.22822}.

\bibitem[Law72]{Law72}
F.~W. Lawvere.
\newblock Teoria delle categorie sopra un topos di base.
\newblock \url{https://github.com/mattearnshaw/lawvere/blob/master/pdfs/1972-perugia-lecture-notes.pdf}, 1972.
\newblock lecture notes from Perugia.

\bibitem[Law15]{Law15}
F.~Lawler.
\newblock {\em Fibrations of Predicates and Bicategories of Relations}.
\newblock PhD thesis, Trinity College, Dublin, 2015.
\newblock Preprint at \url{https://arxiv.org/abs/1502.08017}.

\bibitem[LR20]{LR20}
F.~Loregian and E.~Riehl.
\newblock Categorical notions of fibration.
\newblock {\em Expo. Math.}, 38(4):496--514, 2020.
\newblock \doi{10.1016/j.exmath.2019.02.004}.

\bibitem[LWW10]{LWW10}
S.~Lack, R.~F.~C. Walters, and R.~J. Wood.
\newblock Bicategories of spans as {C}artesian bicategories.
\newblock {\em Theory Appl. Categ.}, 24:No. 1, 1--24, 2010.

\bibitem[Mil00]{Mil00}
S.~Milius.
\newblock Relations in categories, 2000.
\newblock master's thesis, available at \url{https://www8.cs.fau.de/ext/milius/thesis/thesis_a4.pdf}.

\bibitem[Mye18]{Mye18}
D.~J. Myers.
\newblock String diagrams for double categories and equipments, 2018, \burlalt{1612.02762}{http://arxiv.org/abs/1612.02762}.

\bibitem[Nas24]{Nas24}
H.~Nasu.
\newblock An internal logic of virtual double categories, 2024, \burlalt{2410.06792}{http://arxiv.org/abs/2410.06792}.
\newblock \urlprefix\url{https://arxiv.org/abs/2410.06792}.

\bibitem[Nie12]{Nie12}
S.~Niefield.
\newblock Span, cospan, and other double categories.
\newblock {\em Theory Appl. Categ.}, 26:No. 26, 729--742, 2012.

\bibitem[{nLa}23]{nLa}
{nLab authors}.
\newblock sliced adjoint functors -- section.
\newblock \url{https://ncatlab.org/nlab/show/sliced+adjoint+functors+--+section}, Aug. 2023.
\newblock \href{https://ncatlab.org/nlab/revision/sliced+adjoint+functors+--+section/9}{Revision 9}.

\bibitem[Par21]{Par21}
R.~Par\'{e}.
\newblock Morphisms of rings.
\newblock In {\em Joachim {L}ambek: the interplay of mathematics, logic, and linguistics}, volume~20 of {\em Outst. Contrib. Log.}, pages 271--298. Springer, Cham, [2021] \copyright 2021.
\newblock \doi{10.1007/978-3-030-66545-6\_8}.

\bibitem[Pas16]{Pas16}
F.~Pasquali.
\newblock Remarks on the tripos to topos construction: comprehension, extensionality, quotients and functional-completeness.
\newblock {\em Appl. Categ. Structures}, 24(2):105--119, 2016.
\newblock \doi{10.1007/s10485-014-9388-1}.

\bibitem[Pav95]{Pav95}
D.~Pavlovi\'{c}.
\newblock Maps. {I}. {R}elative to a factorisation system.
\newblock {\em J. Pure Appl. Algebra}, 99(1):9--34, 1995.
\newblock \doi{10.1016/0022-4049(94)00054-M}.

\bibitem[Ros99]{Ros99}
G.~Rosolini.
\newblock A note on {C}auchy completeness for preorders.
\newblock {\em Riv. Mat. Univ. Parma (6)}, 2*:89--99, 1999.

\bibitem[SC75]{Suc75}
R.~Succi-Cruciani.
\newblock La teoria delle relazioni nello studio di categorie regolari e di categorie esatte.
\newblock {\em Riv. Mat. Univ. Parma (4)}, 1:143--158, 1975.

\bibitem[Sch15]{Sch15}
P.~Schultz.
\newblock Regular and exact (virtual) double categories, 2015, \burlalt{1505.00712}{http://arxiv.org/abs/1505.00712}.

\bibitem[Shu08]{Shu08}
M.~Shulman.
\newblock Framed bicategories and monoidal fibrations.
\newblock {\em Theory Appl. Categ.}, 20:No. 18, 650--738, 2008.

\bibitem[Shu10]{Shu10}
M.~A. Shulman.
\newblock Constructing symmetric monoidal bicategories, 2010, \burlalt{1004.0993}{http://arxiv.org/abs/1004.0993}.

\bibitem[{\v{S}}t{\v{e}}24]{Ste24}
M.~{\v{S}}t{\v{e}}p\'{a}n.
\newblock Factorization systems and double categories.
\newblock {\em Theory Appl. Categ.}, 41:Paper No. 18, 551--592, 2024.

\bibitem[SW73]{SW73}
R.~Street and R.~F.~C. Walters.
\newblock The comprehensive factorization of a functor.
\newblock {\em Bull. Amer. Math. Soc.}, 79:936--941, 1973.
\newblock \doi{10.1090/S0002-9904-1973-13268-9}.

\bibitem[TY21]{TY21}
W.~Tholen and L.~Yeganeh.
\newblock The comprehensive factorization of {B}urroni's {$\Bbb T$}-functors.
\newblock {\em Theory Appl. Categ.}, 36:Paper No. 8, 206--249, 2021.

\bibitem[Vas14]{Vas14}
C.~Vasilakopoulou.
\newblock Generalization of algebraic operations via enrichment, 2014, \burlalt{1411.3038}{http://arxiv.org/abs/1411.3038}.

\bibitem[Ver11]{Ver11}
D.~Verity.
\newblock Enriched categories, internal categories and change of base.
\newblock {\em Repr. Theory Appl. Categ.}, (20):1--266, 2011.

\bibitem[WW08]{WW08}
R.~F.~C. Walters and R.~J. Wood.
\newblock Frobenius objects in {C}artesian bicategories.
\newblock {\em Theory Appl. Categ.}, 20:No. 3, 25--47, 2008.

\end{thebibliography}

\end{document}